\documentclass[12pt]{article} 
\usepackage[sectionbib]{natbib}
\usepackage{mathrsfs}
\usepackage{fullpage}
\usepackage{amsfonts}
\usepackage{graphicx}
\usepackage{mathrsfs}
\usepackage{amsmath}
\usepackage{amssymb}
\usepackage{epstopdf}
\usepackage{float}
\usepackage{rotating,color}
\usepackage{xcolor}
\usepackage{array,epsfig,fancyheadings,rotating}
\usepackage[]{hyperref}  
\usepackage{sectsty, secdot}
\usepackage{amsmath}
\usepackage{amssymb}
\usepackage{amsfonts}
\usepackage{multirow}
\usepackage{amsthm}

\newtheorem{theorem}{Theorem}
\newtheorem{lemma}{Lemma}

\theoremstyle{definition}

\newtheorem{remark}{Remark}

\usepackage{xr}

\allowdisplaybreaks[4]

\newcommand{\be}{\begin{equation}}
\newcommand{\ee}{\end{equation}}
\newcommand{\beaa}{\begin{eqnarray*}}
\newcommand{\eeaa}{\end{eqnarray*}}
\newcommand{\bea}{\begin{eqnarray}}
\newcommand{\eea}{\end{eqnarray}}
\newcommand{\lbl}{\label}

\def\proof {{\noindent\bf Proof.}\quad}
\def\var{\mathrm {var}}

\def\cov{\mathrm {cov}}
\def\cor{\mathrm {cor}}
\def\diag{\mathrm {diag}}

\newcommand{\bm}{\boldsymbol}

\def\tr{\mathrm {tr}}

\def\A{{\bf A}}

\def\ve{\varepsilon}
\def\z{{\bm z}}

\def\I{{\bf I}}

\def\X{{\bm X}}

\def\Y{{\bm Y}}

\def\tr{\mathrm {tr}}

\def\bms{{\bm\Sigma}}

\def\cd{\mathop{\rightarrow}\limits^{d}}

\def\bmv{\bm \varepsilon}

\begin{document}



\centerline{\large\bf Rank Based Tests for High Dimensional White Noise}
\vspace{.4cm}
\centerline{Dachuan Chen, Fengyi Song and Long Feng}
\vspace{.4cm}
\centerline{\it Nankai University}


\begin{quotation}
\noindent {\it Abstract:}
The development of high-dimensional white noise test is important in both statistical theories and applications, where the dimension of the time series can be comparable to or exceed the length of the time series. This paper proposes several distribution-free tests using the rank based statistics for testing the high-dimensional white noise, which are robust to the heavy tails and do not quire the finite-order moment assumptions for the sample distributions. Three families of rank based tests are analyzed in this paper, including the simple linear rank statistics, non-degenerate U-statistics and degenerate U-statistics. The asymptotic null distributions and rate optimality are established for each family of these tests. Among these tests, the test based on degenerate U-statistics can also detect the non-linear and non-monotone relationships in the autocorrelations. Moreover, this is the first result on the asymptotic distributions of rank correlation statistics which allowing for the cross-sectional dependence in high dimensional data.

\vspace{9pt}
\par
\end{quotation}\par

\def\thefigure{\arabic{figure}}
\def\thetable{\arabic{table}}

\renewcommand{\theequation}{\thesection.\arabic{equation}}

\section{Introduction}

The hypothesis test for white noise is a critical methodology in statistical inference and modeling. It is necessary in diagnosis checking for the linear regression and time series modeling. There has been a vast increase in the amount of high-dimensional data available in recent years, which has received increasing attention from statisticians. The existence of such high-dimensional data is widespread, including the areas of genomics, neuroscience, finance, economics and so on. This brings additional difficulties for the problem of diagnosis checking, which means that in the theoretical development of the test for high-dimensional white noise, the dimension of the time series can be comparable to or exceed the length of the time series.

For the white noise tests designed for univariate time series, many commonly used methodologies are well documented in \citet{li2004diagnostic}. The alternative hypothesis of these tests can be grouped into two different classes: (i) specified alternative in form of some explicit parametric model; (ii) completely unspecified alternative, which means that the departure from white noise can be arbitrary. It is well known that likelihood based tests are more powerful than the omnibus tests under the first class of the alternatives, see, e.g. \citet{chang2017testing}. Under the second class of alternatives, the Box-Pierce portmanteau test and its variations are most popular because of its ease of use in practice {\color{black}and motivate the white noise tests for multivariate time series such as \citet{hosking1980multivariate} and \citet{li1981distribution}. Specially, the Ljung–Box test is a type of statistical test of whether any of a group of autocorrelations of a time series are different from zero. Instead of testing randomness at each distinct lag, it tests the ``overall" randomness based on a number of lags.} These tests enjoy the theoretical benefits of asymptotically distribution-free and $\chi^2$-distributed properties under null hypothesis, see, e.g. \citet{li2004diagnostic} and \citet{lutkepohl2005new}.

There are also some white noise tests constructed for the multivariate time series which assuming that the dimension of the times series is smaller than the length of the time series in asymptotics, see, e.g., \citet{hosking1980multivariate} and \citet{li1981distribution}. However, the existing literature suggests that these tests suffer from the slow convergence to their asymptotic null distributions, see \citet{li2019testing}. This fact calls for the more efficient testing methodologies for multivariate time series, or even high-dimensional time series.

Several omnibus tests for high-dimensional white noise have been developed in recent years, see, e.g., \citet{chang2017testing}, \citet{li2019testing}, \citet{tsay2020testing} and \citet{feng2022whitenoise}. Among these existing theories, the tests proposed in \citet{chang2017testing}, \citet{li2019testing} and \citet{feng2022whitenoise} are distribution-dependent, while the test in \citet{tsay2020testing} is distribution-free. \citet{chang2017testing} developed a max-type test for this purpose based on the maximum absolute auto-correlations and cross-correlations of the component series. \citet{li2019testing} proposed a sum-type test for high dimensional white noise by summing up the squared singular values of the first several lagged sample auto-covariance matrices. In general, the max-type test can only work well under the sparse alternatives where only a few elements in the auto-correlations are nonzero. In contrast, the sum-type test can only work well under the dense alternatives. To test the high dimensional white noise, \citet{feng2022whitenoise} show the asymptotic independence between the max-type test statistic and a new sum-type test statistic. Based on this theoretical result, this paper constructed the Fisher's combination test which is robust to both sparse and dense alternatives. As a distribution-free approach, \citet{tsay2020testing} developed the high-dimensional white noise test based on the Spearman's rank correlation and the theory of extreme values.


More accurately, in this paper we consider the following hypothesis testing problem. Let $\bmv_t$ be a $p$-dimensional weakly stationary time series with mean zero. We want to test the following hypothesis:
\begin{align}
H_0: \{\bmv_t\} ~{\text is ~white~ noise ~versus}~ H_1: \{\bmv_t\} ~is ~not~ white~ noise
\end{align}
{\color{black} In this paper, we said a time series $x_1,\cdots,x_T$ are white noise if they are all independent and identically distributed.} So, under the null hypothesis, $\bmv_{t+k}$ is independent of $\bmv_t$ for all $k>0$. Here the dimension of the time series $p$ is comparable to or even larger than the sample size $n$.

In this paper, we develop the rank based tests for testing the high-dimensional white noise, which are distribution-free. The proposed tests are robust to the heavy tails and do not require the finite-order moment assumptions or any tail assumptions for the sample distribution. There are three families of rank based tests investigated in this paper, including the simple linear rank statistics, non-degenerate U-statistics and degenerate U-statistics, with the examples of Spearman's rho, Kendall's tau, Hoeffding's $D$, Blum-Kiefer-Rosenblatt's $R$ and Bergsma-Dassios-Yanagimoto's $\tau^{*}$. Among these tests, simple linear rank statistics and non-degenerate U-statistics can only work well with the linear or monotone relationships in autocorrelations. In contrast, the degenerate U-statistics can also work well with the non-linear and non-monotone relationships in autocorrelations. As the theoretical results of this paper, we have established the asymptotic null distribution, the power analysis and the rate optimality in terms of power for each family of the rank based test statistics.

Because this paper shows one possible application of the rank correlation statistics in high-dimensional data analysis, we here provide a brief literature review for the rank correlation statistics and point out the theoretical contribution of this paper. \citet{han2017distribution} proposed the rank based tests based on the simple linear rank statistics and non-degenerate U-statistics for testing the mutual independence among all elements in the high-dimensional random vectors. \citet{drton2020high} proposed the hypothesis test based on the degenerate U-statistics with the same purpose as \citet{han2017distribution}. As mentioned earlier, \citet{tsay2020testing} applied the Spearman's rank correlation to the test of high-dimensional white noise. However, the asymptotic distributions of the rank correlation statistics in these three existing literature are all derived based on the assumption of cross-sectional independence in high-dimensional data. Therefore, as the theoretical contribution of our result, this is the first paper in existing literature which established the asymptotic distribution of the rank correlation statistics without assuming the cross-sectional independence.

The main contributions of this paper are summarized as follows.

\begin{enumerate}
\item We develop the rank based tests for testing the high dimensional white noise, which are distribution free. Our test are robust to the heavy tails and do not require the finite-order moment assumptions.
\item Besides the simple linear rank statistics and non-degenerate U-statistics, we also develop the tests for the degenerate U-statistics, which are very useful to detect the non-linear and non-monotone relationships in autocorrelations. Limiting null distributions and the rate optimality in terms of power of these three families of tests are established in this paper.
\item In the existing literature concerning the asymptotic distribution of rank correlation statistics, this paper is the first one on this topic which allowing for the cross-sectional dependence in the high-dimensional data. In contrast, the other existing results are all based on the assumption of cross-sectional independence of the data, see, e.g., \citet{han2017distribution}, \citet{drton2020high} and \citet{tsay2020testing}.
\end{enumerate}

This paper is organized as follows. Section 2 proposes the theoretical results about three families of distribution-free test statistics, including the simple linear rank statistics, non-degenerate U-statistics and degenerate U-statistics. The limiting null distributions of these tests are derived and their rate-optimality in terms of power is also analyzed. Section 3 shows the empirical sizes and the power comparison of the proposed test statistics based on Monte Carlo simulation. Section 4 concludes this paper and discusses several possible directions for the research in the future. All mathematical proofs of the theoretical results in this paper are collected in supplementary material. {\color{black} In the supplementary material, we also consider high dimensional white noise test based on Chatterjee’s rank Correlation \citep{ch2021} and $L$-statistics with the above three-type rank based correlations \citep{cjs2022}.}

\section{Rank based tests}

In this section, we state the theoretical results for three families of rank based methodologies for testing the high-dimensional white noise, including simple linear rank statistics, non-degenerate U-statistics and degenerate U-statistics.

\subsection{Simple linear rank statistics}
{\color{black}First, we restate the definition of relative ranks in Han et al. (2017).} Consider the dependence between $\{(\ve_{1,i},\ve_{k+1,j}),\cdots,(\ve_{n-k,i},\ve_{n,j})\}$ for any two entries $i,j\in\{1, \ldots, p\}$.   Let $Q_{n-k, t}^{i}(k)$ be the rank of $\ve_{t, i}$ in $\left\{\ve_{1, i}, \ldots, \ve_{n-k, i}\right\}$ and let $\tilde{Q}_{n-k, t+k}^{j}(k)$ be the rank of $\ve_{t+k, j}$ in $\{\ve_{k+1,j},\cdots,\ve_{n,j}\}$. Let $R_{n-k, t+k}^{ij}(k)$ be the relative rank of $\ve_{t+k, j}$ compared to $\ve_{t, i}$; that is, $R_{n-k, t+k}^{ij}(k) \equiv \tilde{Q}_{n-k, t'+k}^{j}(k)$ subject to the constraint that $Q_{n-k, t'}^{i}(k)=t$ for $t=1,\cdots,n-k$.

The first family includes tests based on simple linear rank statistics of the form
$$
V_{ij}(k) \equiv (n-k)^{1/2}\sum_{t=1}^{n-k} c_{n-k, t} g\left\{R_{n-k, t+k}^{ij}(k) /(n-k+1)\right\} \quad(i,j \in\{1, \ldots, p\})
$$
where $\left\{c_{n-k, t}, t=1, \ldots, n-k\right\}$ form an array of constants called the regression constants and $g(\cdot)$ is a Lipschitz function called the score function. We assume $\sum_{t=1}^{n-k} c_{n-k,t}^{2}>0$ to avoid triviality. It is immediately clear that Spearman's rho belongs to the family of simple linear rank statistics. To accommodate tests of high-dimensional white noise, we further pose the alignment assumption
$$
c_{n-k, t}=(n-k)^{-1} f\{t /(n-k+1)\}
$$
where $f(\cdot)$ is a Lipschitz function. Under this assumption, the simple linear rank statistic is a general measure of the agreement between the ranks of two sequences. {\color{black}The Spearman’s rho belongs to the family of simple linear rank statistics with $g(x)=f(x)=x-\frac{1}{2}$. }

Under $H_{0},$ the distribution of $V_{ij}(k)$ is irrelevant to the specific distribution of $\ve_t$ for all $i,j \in\{1, \ldots, p\} .$ Accordingly, the mean and variance of $V_{ij}(k)$ are calculable without knowing the true distribution. Let $E_{H_{0}}(\cdot)$ and $\operatorname{var}_{H_{0}}(\cdot)$ be the expectation and variance of a certain statistic under $H_{0} .$ We have
\begin{align}
&E_{H_{0}}\left(V_{ij}(k)\right)= (n-k)^{1/2} \bar{g}_{n-k} \sum_{t=1}^{n-k} c_{n-k, t}, \\
&\sigma_V^2=\operatorname{var}_{H_{0}}\left(V_{ij}(k)\right)=\frac{n-k}{n-k-1} \sum_{t=1}^{n-k}\left[g\{i /(n-k+1)\}-\bar{g}_{n-k}\right]^{2} \sum_{t=1}^{n-k}\left(c_{n-k, t}-\bar{c}_{n-k}\right)^{2} \label{sv}
\end{align}
where $\bar{g}_{n-k} \equiv (n-k)^{-1} \sum_{t=1}^{n-k} g\{t /(n-k+1)\}$ is the sample mean of $g\left\{R_{n-k, t}^{ij} /(n-k+1)\right\}(t=1, \ldots, n-k)$ and $\bar{c}_{n-k}=(n-k)^{-1} \sum_{t=1}^{n-k}c_{n-k, t}$. Based
on $\left\{V_{ij}(k), 1 \leqslant i,j \leqslant p, 1\le k\le K\right\},$ we propose the following statistic for testing $H_{0}:$
$$
V_{n} \equiv \max_{1\le k\le K} \max _{1\le i, j\le p}\left|V_{ij}(k)-E_{H_{0}}\left(V_{ij}(k)\right)\right|
$$
{\color{black}Note that we can allow $K$ to grow as $n$ increase.}
Let $N=Kp^2$. We define the following assumption for any matrix $\bms$:
\begin{itemize}
\item[(A1)]Let $\boldsymbol{\Sigma}=\left(\sigma_{i j}\right)_{1 \leq i, j \leq N}$. For some $\varrho \in(0,1),$ assume $\left|\sigma_{i j}\right| \leq \varrho$ for all $1 \leq i<j \leq N$ and
$N \geq 2 .$ Suppose $\left\{\delta_{N} ; N \geq 1\right\}$ and $\left\{\varsigma_{N} ; N \geq 1\right\}$ are positive constants with $\delta_{N}=o(1 / \log N)$
and $\varsigma=\varsigma_{N} \rightarrow 0$ as $N \rightarrow \infty .$ For $1 \leq i \leq N,$ define $B_{N, i}=\left\{1 \leq j \leq N ;\left|\sigma_{i j}\right| \geq \delta_{N}\right\}$
and $C_{N}=\left\{1 \leq i \leq N ;\left|B_{N, i}\right| \geq N^{\varsigma}\right\} .$ We assume that $\left|C_{N}\right| / N \rightarrow 0$ as $N \rightarrow \infty$.
\end{itemize}

Here we define $\{\nu_1,\cdots,\nu_{Kp^2}\}=\{V_{ij}(k)/\sigma_V\}_{1\le i,j\le p, 1\le k \le K}$. Define  $\sigma^V_{ij}=\cor(\nu_i,\nu_j)$ and $\bms_V=(\sigma_{ij}^V)_{1\le i,j\le N}$.

To derive the limiting null distribution of simple linear rank statistics, we need the following conditions.

\begin{itemize}
\item[(C1)] The regression constants $\left\{c_{n-k, 1}, \ldots, c_{n-k, n-k}\right\}$ satisfying
\begin{align*}
&\max _{1 \leqslant i \leqslant n-k}\left|c_{n-k, i}-\bar{c}_{n-k}\right|^{2} \leqslant \frac{C_{1}^{2}}{n-k} \sum_{i=1}^{n-k}\left(c_{n-k, i}-\bar{c}_{n-k}\right)^{2},\\
&\left|\sum_{i=1}^{n-k}\left(c_{n-k, i}-\bar{c}_{n-k}\right)^{3}\right|^{2} \leqslant \frac{C_{2}^{2}}{n-k}\left\{\sum_{i=1}^{n-k}\left(c_{n-k, i}-\bar{c}_{n-k}\right)^{2}\right\}^{3}
\end{align*}
where $\bar{c}_{n-k} \equiv \sum_{i=1}^{n-k} c_{n-k, i}$ represents the sample mean of the regression constants and $C_{1}$ and $C_{2}$ are two constants.
\item[(C2)] The score function $g(\cdot)$ is differentiable with bounded Lipschitz constant.
\item[(C3)] The correlation matrix $\bms_V$ satisfies Assumption (A1).
\end{itemize}

{\color{black}
{\bf Remark:} The assumption (A1) is the same as the condition (2.2) in Feng et al. (2022a), which demands the number of variables that are strongly-correlated with many other variables should not be too much. If the eigenvalues of $\bms_V$ are all bounded, we have $\max_{1\le i\le N} \sum_{j=1}^N \sigma_{ij}^{V2}\le C$ for some constant $C>0$. Then, let $\delta_N=(\log N)^{-2}$ for $N\ge e^e$, so for each $1 \leqslant i \leqslant N, \delta_N^2 \cdot\left|B_{N, i}\right| \leqslant \sum_{j=1}^N \sigma_{i j}^{V2} \leqslant C$. Hence, $\left|B_{N, i}\right| \leqslant C \cdot(\log N)^2<N^\kappa$ where $\kappa=\kappa_N:=5(\log \log N) / \log N$ for large $N$ . As a result, $\left|C_N\right|=0$ and condition (C3) holds. Condition (C1) is commonly used to deviate the asymptotical normality of the simple linear rank statistics, see Hájek et al. (1999) and Kallenberg (1982). If $f$ is a linear function, Condition (C1) will hold directly.

}

Next, we state the theoretical result about the limiting null distribution of simple linear rank statistics.

\begin{theorem}\label{th11}
Suppose (C1)-(C3) hold. Then, under $H_{0}$, for any $y \in \mathbb{R}$, we have
$$
\left|\operatorname{P}\left(V_{n}^{2} / \sigma_{V}^{2}-2\log(Kp^2)+\log \log (Kp^2) \leqslant y\right)-\exp \left\{-\pi^{-1 / 2} \exp (-y / 2)\right\}\right|=o(1)
$$
where $\sigma_{V}^{2}=\var_{H_0}(V_{ij}(k))$  if $N=o(n^\epsilon)$ as $n\to \infty$ for some positive constant $\epsilon$.
\end{theorem}

We propose the following size-$\alpha$ test $T^V_\alpha$ of $H_0$:
\begin{equation}
T_\alpha^V\doteq I\left(V_{n}^{2} / \sigma_{V}^{2}-2\log(Kp^2)+\log \log (Kp^2)\ge q_\alpha\right), \label{DEF_Test_TV}
\end{equation}
where $q_\alpha=-\log (\pi)-2\log\log(1-\alpha)^{-1}$.

To specify the alternative hypothesis, we introduce a notation for a set of vectors which satisfying some specific condition. Define $N=Kp^2$.
Let $\mathcal{U}(c)$ be a set of vectors indexed by a constant $c$:
$$
\mathcal{U}(c) \equiv\left\{M=(m_{l})_{1 \le l \le N} \in \mathbb{R}^{N} \mid \max_{1\le l\le N} m_{l} \geqslant c(\log N)^{1 / 2}\right\}.
$$
Based on the above definition, we know that $\mathcal{U}(c)$ is the set of vectors of which at least one element has magnitude greater than $c(\log N)^{1 / 2}$ for some large enough constant $c>0$.

Next, we specify the sparse local alternative based on $\mathcal{U}(c)$. We define the random vector $\hat{V}=\left[\hat{V}_{ij}(k)\right] \in \mathbb{R}^{N}$ by
$$
\hat{V}_{ij}(k)=\sigma_{V}^{-1}\left\{V_{ij}(k)-E_{H_{0}}\left(V_{ij}(k)\right)\right\},\quad(1 \leqslant i,j \leqslant p; 1 \leqslant k \leqslant K)
$$
where $\sigma_{V}$ is defined in (\ref{sv}) and $\left\{V_{ij}(k), 1 \leqslant i,j \leqslant p,1\le k\le K\right\}$ are the simple linear rank statistics. Let the population version of $\hat{V}$ be $V \equiv E(\hat{V})$. We study the power of tests against the alternative
$$
H_{\mathrm{a}}^{V}(c) \equiv\{F(\bmv): V\{F(\bmv)\} \in \mathcal{U}(c)\}
$$
where $F(\bmv)$ is the joint distribution function of $\bmv$ and we write $V\{F(\bmv)\}$ to emphasize that $V=E(\hat{V})=\int \hat{V} \mathrm{d} F(\bmv)$ is a function of $F(\bmv)$.

The following theorem now describe the conditions under which the power of the test based on simple linear rank statistics converges to one as $n$ and $p$ going to infinity, under the sparse local alternative $H_{\mathrm{a}}^{V}$.

\begin{theorem}\label{th12}
Assume Conditions (C1)-(C3) hold. And assume that $\sigma_{V}^{2}=A_{1}\{1+o(1)\}$ and $\max \{|f(0)|,|g(0)|\} \leqslant A_{2}$ for some positive constants $A_{1}$ and $A_{2}$. Further assume that $f(\cdot)$ and $g(\cdot)$ have bounded Lipschitz constants. Then, for some large scalar $B_{1}$ depending only on $A_{1}, A_{2}$ and the Lipschitz constants of $f(\cdot)$ and $g(\cdot)$
$$
\inf _{F(\bmv) \in H_{\mathrm{a}}^{V}\left(B_{1}\right)} \operatorname{pr}\left(T^V_{\alpha}=1\right)=1-o(1)
$$
where the infimum is taken over all distributions $F(\bmv)$ such that $V\{F(\bmv)\} \in \mathcal{U}\left(B_{1}\right)$.
\end{theorem}

Define $r_{ij}(k)$ is the correlation between $\ve_{t,i}$ and $\ve_{t+k,j}$. To investigate the rate optimality of the test based on simple linear rank statistics, we need the following assumption for the distribution:
\begin{itemize}
\item[(A2)] When $\bmv$ is Gaussian, suppose that for large $n$ and $p$, $c V_{ij}(k) \leqslant r_{ij}(k) \leqslant C V_{ij}(k)$ for $1\le i,j\le p, 1\le k\le K$ {\color{black}with probability tending to one}, where $c$ and $C$ are two constants.
\end{itemize}

For each $n$, define $\mathcal{T}_{\alpha}$ to be the set of all measurable size-$\alpha$ tests. In other words, $\mathcal{T}_{\alpha}:=\{T_{\alpha}: \operatorname{pr}(T_{\alpha}=1 | H_0) \le \alpha\}$.

Finally, the rate optimality result can be stated by the following theorem. Recall that $T_\alpha^V$ defined in (\ref{DEF_Test_TV}) can correctly reject the null hypothesis provided that at least one element in $V$ has magnitude greater than $c(\log N)^{1 / 2}$ for some constant $c$. In the following theorem, we show that the rate of the signal gap $(\log N)^{1 / 2}$ cannot be further relaxed.

\begin{theorem}\label{th13}
Suppose that the simple linear rank statistics $\left\{V_{ij}(k), 1 \leqslant i,j \leqslant p,1\le k\le K\right\}$ satisfy all the conditions in Theorems \ref{th11} and \ref{th12}. Suppose also that Assumption (A2) holds. Then,  the corresponding size-$\alpha$ test $T^V_{\alpha}$ is rate-optimal. In other
words, there exist two constants $D_{1}<D_{2}$ such that:
\begin{itemize}
\item[(i)] $\sup _{F(\bmv) \in H_{\mathrm{a}}^{V}\left(D_{2}\right)} \operatorname{pr}\left(T_{\alpha}=0\right)=o(1)$;
\item[(ii)] for any $\beta>0$ satisfying $\alpha+\beta<1,$ for large $n$ and $p$ we have
$$
\inf _{T_{\alpha} \in \mathcal{T}_{\alpha}} \sup _{F(\bmv) \in H_{\mathrm{a}}^{V}\left(D_{1}\right)} \operatorname{pr}\left(T_{\alpha}=0\right) \geqslant 1-\alpha-\beta
$$
\end{itemize}
\end{theorem}

The above theorem means that any measurable size-$\alpha$ test cannot distinguish between the null hypothesis and the sparse alternative when the coefficient $c$ in $H_{\mathrm{a}}^{V}(c)$ is small enough.

As an example of simple linear rank statistic, we state the high-dimensional white noise test based on the Spearman's rho as follows.

\paragraph{Example 1 (Spearman's rho).}  Recall that $Q_{n-k, t}^{i}(k)$ and $\tilde{Q}_{n-k, t+k}^{j}(k)$ be the ranks of $\ve_{t, i}$ and $\ve_{t+k, j}$ among $\left\{\ve_{1, i}, \ldots, \ve_{n-k, i}\right\}$ and $\{\ve_{k+1,j},\cdots,\ve_{n,j}\}$, respectively. Let $R_{n-k, t+k}^{ij}(k)$ be the relative rank of $\ve_{t+k, j}$ compared to $\ve_{t, i}$; that is, $R_{n-k, t+k}^{ij}(k) \equiv \tilde{Q}_{n-k, t'+k}^{j}(k)$ subject to the constraint that $Q_{n-k, t'}^{i}(k)=t$ for $t=1,\cdots,n-k$. Spearman's rho is defined as
$$
\begin{aligned}
\rho_{i j}(k) &=\frac{\sum_{t=1}^{n-k}\left(Q_{n-k, t}^{i}(k)-\bar{Q}_{n-k}^{i}(k)\right)\left(\tilde{Q}_{n-k, t+k}^{j}(k)-\bar{\tilde{Q}}_{n-k}^{j}(k)\right)}{\left\{\sum_{t=1}^{n-k}\left(Q_{n-k, t}^{i}(k)-\bar{Q}_{n-k}^{i}(k)\right)^{2} \sum_{t=1}^{n-k}\left(\tilde{Q}_{n-k, t+k}^{j}(k)-\bar{\tilde{Q}}_{n-k}^{j}(k)\right)^{2}\right\}^{1 / 2}} \\
&=\frac{12}{(n-k)\left((n-k)^{2}-1\right)} \sum_{t=1}^{n-k}\left(i-\frac{n-k+1}{2}\right)\left(R_{n-k, t+k}^{ij}(k)-\frac{n-k+1}{2}\right) \quad(i,j \in\{1, \ldots, p\})
\end{aligned}
$$
where $\bar{Q}_{n-k}^{i}(k)=\bar{\tilde{Q}}_{n-k}^{j}(k) \equiv(n-k+1) / 2 .$ This is a simple linear rank statistic, and we have
$$
E_{H_{0}}\left(\rho_{ij}(k)\right)=0, \quad \operatorname{var}_{H_{0}}\left(\rho_{ij}(k)\right)=(n-k-1)^{-1} \quad(i,j\in\{1, \ldots, p\})
$$
According to (\ref{DEF_Test_TV}), the corresponding test statistic is
$$
L_{\rho}=I\left\{ \max_{1\le k\le K} \max _{1\le i,j\le p} (n-k) \rho_{ij}(k)^{2}-2 \log(Kp^2)+\log \log (Kp^2) \geqslant q_{\alpha}\right\}
$$
where $q_{\alpha} \equiv-\log ( \pi)-2 \log \log (1-\alpha)^{-1}$.

\subsection{Non-degenerate U-statistics}
The second family includes the tests based on non-degenerate U-statistics of the form {\color{black}(Han et al. 2017)}
\begin{align}
U_n=\max_{1\le k\le K} \max_{1\le i,j\le p}(n-k)^{1/2}\left| {U}_{ij}(k)-E_{H_0}(U_{ij}(k))\right|
\end{align}
where $A_{n-k}^m=(n-k)(n-k-1)\cdots(n-k-m+1)$,
\begin{align}\label{ud}
{U}_{ij}(k)=&\frac{1}{A_{n-k}^m}\sum_{1\le t_1\not=t_2,\cdots,\not=t_m\le n-k}h((\varepsilon_{t_1,i},\varepsilon_{t_1+k,j})^\top,\cdots,(\varepsilon_{t_m,i},\varepsilon_{t_m+k,j})^\top)
\end{align}
Here $U_{ij}(k)$ depends only on  $\{R_{n-k,t}^{ij}(k)\}_{t=k+1}^{n}$.
For our purposes $h$ may always be assumed to be bounded but not necessarily symmetric. The boundedness assumption is mild since correlation
is the object of interest.

Further concepts concerning U-statistics are needed to state the assumption for the derivation of the limiting null distribution. For $m \in \mathbb{Z}^{+}$, we define $[m]=\{1,2,\cdots,m\}$ and write $\mathcal{P}_m$ for the set of all $m!$ permutations of $[m]$. For any kernel $h(\cdot),$ any number $\ell \in[m],$ and any measure $\mathbb{P}_{\boldsymbol{Z}},$ we write
\begin{equation}\label{h10}
h_{\ell}\left(\boldsymbol{z}_{1} \ldots, \boldsymbol{z}_{\ell} ; \mathbb{P}_{\boldsymbol{Z}}\right):=\mathbb{E} h\left(\boldsymbol{z}_{1} \ldots, \boldsymbol{z}_{\ell}, \boldsymbol{Z}_{\ell+1}, \ldots, \boldsymbol{Z}_{m}\right)
\end{equation}
and
\begin{equation}\label{hl}
h^{(\ell)}\left(\boldsymbol{z}_{1}, \ldots, \boldsymbol{z}_{\ell} ; \mathbb{P}_{\boldsymbol{Z}}\right):=h_{\ell}\left(\boldsymbol{z}_{1}, \ldots, \boldsymbol{z}_{\ell} ; \mathbb{P}_{\boldsymbol{Z}}\right)-\mathbb{E} h-\sum_{k=1}^{\ell-1} \sum_{1 \leq i_{1}<\cdots<i_{k} \leq \ell} h^{(k)}\left(\boldsymbol{z}_{i_{1}}, \ldots, \boldsymbol{z}_{i_{k}} ; \mathbb{P}_{\boldsymbol{Z}}\right)
\end{equation}
where $\boldsymbol{Z}_{1}, \ldots, \boldsymbol{Z}_{m}$ are $m$ independent random vectors with distribution $\mathbb{P}_{\boldsymbol{Z}}$ and $\mathbb{E} h:=\mathbb{E} h\left(\boldsymbol{Z}_{1}, \ldots, \boldsymbol{Z}_{m}\right)$. The kernel as well as the corresponding U-statistic is non-degenerate under $\mathbb{P}_{\boldsymbol{Z}}$ if the variance of $h_{1}(\cdot)$ is not zero.

Based on above definitions, we state the following conditions which are needed to derive the limiting null distribution.

\begin{itemize}
\item[(C4)] The kernel function $h(\cdot)$ is bounded and non-degenerate.
\item[(C5)] The correlation matrix of $U_{ij}(k)$--$\bms_{U}$ satisfies Assumption (A1).
\end{itemize}

The following theorem show the asymptotic distribution of the non-degenerate U-statistics under the null hypothesis.

\begin{theorem}\label{th21}
Suppose (C4)-(C5) hold. Then under $H_{0}$, for any $y \in \mathbb{R}$ we have
$$
\left|\operatorname{P}\left(U_{n}^{2} / \sigma_{U}^{2}-2 \log(Kp^2)+\log \log (Kp^2) \leqslant y\right)-\exp \left\{-\pi^{-1 / 2} \exp (-y / 2)\right\}\right|=o(1)
$$
where $\sigma_{U}^{2}=(n-k)\var_{H_0}(U_{ij}(k))$   if $N=o(n^\epsilon)$ as $n\to \infty$ for some positive constant $\epsilon$.
\end{theorem}

We propose the following size-$\alpha$ test $T^U_\alpha$ of $H_0$:
\begin{equation}
T_\alpha^U\doteq I\left(U_{n}^{2} / \sigma_{U}^{2}-2\log(Kp^2)+\log \log (Kp^2)\ge q_\alpha\right)
\end{equation}
where $q_\alpha=-\log (\pi)-2\log\log(1-\alpha)^{-1}$.

To specify the sparse local alternative for the tests based on non-degenerate U-statistics, we first define the random vector $\hat{U}=\left[\hat{U}_{ij}(k)\right] \in \mathbb{R}^{N}$ by
$$
\hat{U}_{ij}(k)=\sigma_{U}^{-1}(n-k)^{1/2}\left\{U_{ij}(k)-E_{H_{0}}\left(U_{ij}(k)\right)\right\},\quad(1 \leqslant i,j \leqslant p; 1 \leqslant k \leqslant K)
$$
where $\sigma_{U}$ is defined in Theorem \ref{th21} and $\left\{U_{ij}(k), 1 \leqslant i,j\leqslant p,1\le k\le K\right\}$ are the non-degenerate U-statistics. Let the population version of $\hat{U}$ be $U \equiv E(\hat{U}) .$ We study the power of tests against the alternative
$$
H_{\mathrm{a}}^{U}(c) \equiv\{F(\bmv): U\{F(\bmv)\} \in \mathcal{U}(c)\}
$$
where $F(\bmv)$ is the joint distribution function of $\bmv$ and we write $U\{F(\bmv)\}$ to emphasize that $U=E(\hat{U})=\int \hat{U} \mathrm{d} F(\bmv)$ is a function of $F(\bmv)$.

The following theorem states the conditions which are required to establish the convergence of the power of ${T}^U_{\alpha}$ to one as $n$ and $p$ going to infinity under the sparse alternative.

\begin{theorem}\label{th22}
 Suppose that the kernel function $h(\cdot)$ in (\ref{ud}) is bounded with $|h(\cdot)| \leqslant A_{3}$ and
$$
m^{2} \operatorname{var}_{H_{0}}\left[E_{H_{0}}\left\{h\left((X_{11},X_{12})^\top, \ldots, (X_{m1},X_{m2})^\top\right) \mid (X_{11},X_{12})^\top\right\}\right]=A_{4} \{1+o(1)\}
$$
for some positive constants $A_{3}$ and $A_{4} .$ Then, for some large scalar $B_{2}$ depending only on $A_{3}, A_{4}$ and $m$,
$$
\inf _{F(\bmv) \in H_{\mathrm{a}}^{U}\left(B_{2}\right)} P\left({T}_{\alpha}=1\right)=1-o(1)
$$
where the infimum is taken over all distributions $F(\bmv)$ such that $U\{F(\bmv)\} \in \mathcal{U}\left(B_{2}\right) .$
\end{theorem}

To study the rate optimality in terms of power for the tests based on non-degenerate U-statistics, we need the following assumption for the distribution:
\begin{itemize}
\item[(A3)] When $\bmv$ is Gaussian, suppose that for non-degenerate U-statistics $U_{ij}(k)$ and large $n$ and $p$, $c U_{ij}(k) \leqslant r_{ij}(k) \leqslant C U_{ij}(k)$ for $1\le i,j\le p, 1\le k\le K$ {\color{black}with probability tending to one}, where $c$ and $C$ are two constants.
\end{itemize}

The rate optimality result and related conditions for the tests based on the non-degenerate U-statistics can be shown as follows, which implies that the rate of the signal gap $(\log N)^{1 / 2}$ cannot be further relaxed.

\begin{theorem}\label{th23}
Suppose that Non-degenerate U-statistics $\left\{U_{ij}(k), 1 \leqslant i,j \leqslant p,1\le k\le K\right\}$ satisfy all the conditions in Theorems \ref{th21} and \ref{th22}. Suppose also that Assumption (A3) holds. Then,  the corresponding size-$\alpha$ test $T^U_{\alpha}$ is rate-optimal. In other
words, there exist two constants $D_{3}<D_{4}$ such that:
\begin{itemize}
\item[(i)] $\sup _{F(\bmv) \in H_{\mathrm{a}}^{U}\left(D_{4}\right)} \operatorname{pr}\left(T_{\alpha}=0\right)=o(1)$;
\item[(ii)] for any $\beta>0$ satisfying $\alpha+\beta<1,$ for large $n$ and $p$ we have
$$
\inf _{T_{\alpha} \in \mathcal{T}_{\alpha}} \sup _{F(\bmv) \in H_{\mathrm{a}}^{U}\left(D_{3}\right)} \operatorname{pr}\left(T_{\alpha}=0\right) \geqslant 1-\alpha-\beta.
$$
\end{itemize}

\end{theorem}

As an example of non-degenerate U-statistics, we state the high-dimensional white noise test based on the Kendall's tau as follows.

\paragraph{Example 2 (Kendall's tau).} Kendall's tau is defined, for $i,j \in\{1, \ldots, p\},$ by
$$
\begin{aligned}
\tau_{ij}(k)& = \frac{2}{(n-k)(n-k-1)}\sum_{1 \le l<l^{\prime} \le n-k} \operatorname{sign}\left(\ve_{l^{\prime}, i}-\ve_{l, i}\right) \operatorname{sign}\left(\ve_{l^{\prime}+k, j}-\ve_{l+k, j}\right) \\
&=\frac{2}{(n-k)(n-k-1)} \sum_{1 \le l<l^{\prime} \le n-k} \operatorname{sign}\left(R_{n-k, l^{\prime}+k}^{ij}(k) - R_{n-k, l+k}^{ij}(k)\right)
\end{aligned}
$$
where the sign function $\operatorname{sign}(\cdot)$ is defined as $\operatorname{sign}(x)=x /|x|$ with the convention $0 / 0=0$. This statistic is a function of the relative ranks $\left\{R_{n-k, t+k}^{ij}(k), t=1, \ldots, n-k\right\}$ and is also a $U$-statistic with bounded kernel $h\left(x_{1,\{1,2\}}, x_{2,\{1,2\}}\right) \equiv \operatorname{sign}\left(x_{1,1}-x_{2,1}\right) \operatorname{sign}\left(x_{1,2}-x_{2,2}\right)$. Accordingly, Kendall's
tau is a rank-type $U$-statistic. Moreover,
$$
E_{H_{0}}\left(\tau_{ij}(k)\right)=0, \quad \operatorname{var}_{H_{0}}\left(\tau_{ij}(k)\right)=\frac{2(2 (n-k)+5)}{9 (n-k)(n-k-1)} \quad(i,j \in\{1, \ldots, p\})
$$
According to (8), the proposed test statistic based on Kendall's tau is
$$
L_{\tau}=I\left\{  \max_{1\le k\le K} \max _{1\le i,j\le p} \frac{9 (n-k)(n-k-1)}{2(2 (n-k)+5)} \tau_{ij}(k)^{2}-2\log(Kp^2)+\log \log (Kp^2) \geqslant q_{\alpha}\right\}
$$

\subsection{Degenerate U-statistics}

The third family includes the tests based on degenerate U-statistics, which are very useful to detect the non-linear and non-monotone relationships in the autocorrelations. We use the term completely degenerate to indicate that the variances of $h_{1}(\cdot), \ldots, h_{m-1}(\cdot)$ are all zero. Finally, let $\mathbb{P}_{0}$ be the uniform distribution on $[0,1],$ and write $\mathbb{P}_{0} \otimes \mathbb{P}_{0}$ for its product measure, the uniform distribution on $[0,1]^{2}$.


In order to derive the limiting null distribution and establish the theoretical results related to the power of the tests based on the degenerate U-statistics, we need the following assumption concerning the kernel function $h$.

\begin{itemize}
\item[(C6)] The kernel $h$ is rank-based, symmetric, and has the following three properties:
\begin{itemize}
\item[(i)] $h$ is bounded.
\item[(ii)] $h$ is mean-zero and degenerate under independent continuous margins, i.e., $\mathbb{E}\left\{h_{1}\left(\boldsymbol{Z}_{1} ; \mathbb{P}_{0} \otimes\right.\right.$ $\left.\left.\mathbb{P}_{0}\right)\right\}^{2}=0$ as $\boldsymbol{Z}_{1} \sim \mathbb{P}_{0} \otimes \mathbb{P}_{0}$
\item[(iii)] $h_{2}\left(\boldsymbol{z}_{1}, \boldsymbol{z}_{2} ; \mathbb{P}_{0} \otimes \mathbb{P}_{0}\right)$ has uniformly bounded eigenfunctions, that is, it admits the expansion
$$
h_{2}\left(\boldsymbol{z}_{1}, \boldsymbol{z}_{2} ; \mathbb{P}_{0} \otimes \mathbb{P}_{0}\right)=\sum_{v=1}^{\infty} \lambda_{v} \phi_{v}\left(\boldsymbol{z}_{1}\right) \phi_{v}\left(\boldsymbol{z}_{2}\right)
$$
where $\left\{\lambda_{v}\right\}$ and $\left\{\phi_{v}\right\}$ are the eigenvalues and eigenfunctions satisfying the integral equation
$$
\begin{aligned}
\mathbb{E} h_{2}\left(\boldsymbol{z}_{1}, \boldsymbol{Z}_{2}\right) \phi\left(\boldsymbol{Z}_{2}\right) &=\lambda \phi\left(\boldsymbol{z}_{1}\right) \text { for all } \boldsymbol{z}_{1} \in \mathbb{R}^{2} \\
\text {with } \boldsymbol{Z}_{2} \sim \mathbb{P}_{0} \otimes \mathbb{P}_{0}, \lambda_{1} \geq \lambda_{2} \geq \cdots \geq 0, \Lambda &=\sum_{v=1}^{\infty} \lambda_{v} \in(0, \infty), \text { and } \sup _{v}\left\|\phi_{v}\right\|_{\infty}<\infty.
\end{aligned}
$$
\end{itemize}
\end{itemize}

The first requirement about boundedness property can be easily verified for the rank correlations which are commonly used, for example, Spearman's rho, Kendall's tau and many others. The other two requirements are much more specific, but can be satisfied by some typical rank correlation measures as long as their consistency properties are known. Moreover, it is easy to see that the assumption $\Lambda>0$ implies $\lambda_{1}>0,$ so that $h_{2}(\cdot)$ is not a constant function.

We also need the following condition to derive the limiting null distribution for the degenerate U-statistics. We first make several definitions which will be used in the following condition. Define a quantity $\theta$, which is any absolute constant such that
$$
\begin{aligned} \label{DEF_theta}
\theta < \sup \left\{ q \in [0, 1/3) : \sum_{v > [n^{(1-3q)/5}]} \lambda_v = O(n^{-q}) \right\}
\end{aligned}
$$
if infinitely many eigenvalues $\lambda_v$ are nonzero, and $\theta = 1/3$ otherwise. Define $\omega_{l,v}=(n-k_l)^{-1/2}\sum_{t=1}^{n-k_l} \phi_v(Z_{t,l})$ for $l=1,\cdots,N$, $v=1,\cdots,M$ where $M=[n^{(1-3\theta)/5}]$ and $Z_{t,l}$ is the corresponding $Z$ of $U_{ij}(k)$ in Condition (C6). Let $b_{lv, rs}=\cov(\omega_{l,v},\omega_{r,s})$ for  $1\le l,r\le N, 1\le v,s\le M$.
    Let $\bm \omega_l=(\omega_{l,1},\cdots,\omega_{l,M})$ and $\Xi_l=\bms_l\bms_l^\top$ where $\bms_{l}\in \mathbb{R}^{M\times (N-1)M}$ is the covariance matrix between $\bm \omega_l$ with $\bm \omega_{r}, r\in \{1,\cdots,N\}\setminus\{l\}$.

\begin{itemize}
\item[(C7)]  There exists a constant $\delta\in (0,1)$ satisfying $\lambda_{max}(\Xi_l)\le \delta$ for all $1\le l\le N$. Suppose $\left\{\delta_{N} ; N \geq 1\right\}$ and $\left\{\varsigma_{N} ; N \geq 1\right\}$ are positive constants with $\delta_{N}=o(1 /\log N)$
and $\varsigma=\varsigma_{N} \rightarrow 0$ as $N \rightarrow \infty .$ Let $\Xi_{ij}=\cov(\bm \omega_i,\bm \omega_j)$. For $1 \leq i \leq N,$ define $B_{N, i}=\left\{1 \leq j \leq N \mid \lambda_{max}(\Xi_{ij}\Xi_{ij}^\top)\geq \delta_{N}^{2+2c}\right\}$ for some constant $c>0$
and $C_{N}=\left\{1 \leq i \leq N ;\left|B_{N, i}\right| \geq N^{\varsigma}\right\} .$ We assume that $\left|C_{N}\right| / N \rightarrow 0$ as $N \rightarrow \infty$.

\end{itemize}

In the following theorem, we show the limiting null distribution and related conditions for the degenerate U-statistics.

\begin{theorem}\label{th31}
Under conditions (C6)-(C7).
Then  for any absolute constant $y \in \mathbb{R}$ that
$$
\begin{aligned}
& \mathbb{P}\left\{\max_{1\le k\le K} \max _{1\le i,j\le p} \frac{n-k-1}{\lambda_{1}\left(\begin{array}{c}
m \\
2
\end{array}\right)} {U}_{ij}(k)-2 \log (Kp^2)-\left(\mu_{1}-2\right) \log \log (Kp^2)+\frac{\Lambda}{\lambda_{1}} \leq y\right\} \\
=& \exp \left\{-\frac{  \kappa}{\Gamma\left(\mu_{1} / 2\right)} \exp \left(-\frac{y}{2}\right)\right\}+o(1)
\end{aligned}
$$
for $\log N=o(n^\theta)$ as $n\to \infty$. Here $\mu_1$ is the multiplicity of the largest eigenvalue $\lambda_1$ in the sequence $\{\lambda_1, \lambda_2, \ldots\}$, $\kappa:= \prod_{v=\mu_1+1}^{\infty} (1-\lambda_v/\lambda_1)^{-1/2}$ and $\Gamma(z):=\int_{0}^{\infty} x^{z-1}e^{-x}dx$ is the gamma function.
\end{theorem}

We propose the following size-$\alpha$ test $T_\alpha^D$ for degenerate U-statistics:
\begin{equation} \label{DEF_T_D_alpha}
T_\alpha^D=I\left(\max_{1\le k\le K} \max _{1\le i,j\le p} \frac{n-k-1}{\lambda_{1}\left(\begin{array}{c}
m \\
2
\end{array}\right)} {U}_{ij}(k)-2 \log (Kp^2)-\left(\mu_{1}-2\right) \log \log (Kp^2)+\frac{\Lambda}{\lambda_{1}}\ge \tilde q_\alpha\right)
\end{equation}
where $\tilde q_\alpha$ is the $1-\alpha$ quantile of the Gumbel distribution function $ \exp \left\{-{  \kappa}/{\Gamma\left(\mu_{1} / 2\right)} \exp \left(-{y}/{2}\right)\right\}$, i.e.,
$$
\tilde q_\alpha=-\log\left(\frac{\Gamma^2\left(\mu_{1} / 2\right)}{  \kappa^2}\right)-2\log\log(1-\alpha)^{-1}.
$$
It is easy to show that $\mathbb{P}_{H_{0}}\left(T_{\alpha}^D=1\right)=\alpha+o(1)$

We study the power of the proposed test based on the degenerate U-statistics from now on. It is necessary to introduce a new distribution family which is also useful to specify the alternative.
Recall the definition of $h^{(1)}(\cdot)$ in (\ref{hl}) . For any kernel function $h(\cdot)$ and constants $\gamma>0$ and $n,p\in \mathbb{Z}^{+},$ define a general $p$-dimensional (not necessarily continuous) distribution family as follows:
$$
\mathcal{D}(\gamma, np ; h):=\left\{F(\boldsymbol{X}): \boldsymbol{X} \in \mathbb{R}^{np}, \operatorname{Var}_{ij k}\left\{h^{(1)}\left(\cdot ; \mathbb{P}_{ij k}\right)\right\} \leq \gamma \mathbb{E}_{ij k} h \text { for all } 1\le i,j\le p, 1\le k\le K\right\}
$$
where $F(\boldsymbol{X})$ is the distribution (law) of $\boldsymbol{X},$ and $\mathbb{P}_{ij k}, \mathbb{E}_{ij k}(\cdot),$ and $\operatorname{Var}_{ij k}(\cdot)$ stand for the probability measure, expectation, and variance operated on the bivariate distribution of $\left(\varepsilon_{ti}, \varepsilon_{t+k,j}\right)^{\top},$ respectively. The family $\mathcal{D}(\gamma, np ; h)$ intrinsically characterizes the slope of the non-negative function $\operatorname{Var}_{ij k}\left\{h^{(1)}\left(\cdot ; \mathbb{P}_{ij k}\right)\right\}$ with regard to the dependence between $\varepsilon_{ti}$ and $\varepsilon_{t+k,j},$ characterized by the non-negative correlation measure $\mathbb{E}_{ij k} h$. Under the null hypothesis, we have
$$
\operatorname{Var}_{ij k}\left\{h^{(1)}\left(\cdot ; \mathbb{P}_{ij k}\right)\right\} = \mathbb{E}_{ij k} h =0
$$
provided that Condition (C6) holds for $h(\cdot)$. Therefore, as the dependence between $\varepsilon_{ti}$ and $\varepsilon_{t+k,j}$ increasing, it can be expected that the variance $\operatorname{Var}_{ij k}\left\{h^{(1)}\left(\cdot ; \mathbb{P}_{ij k}\right)\right\}$ will depart away from zero with the same or a slower rate compared to $\mathbb{E}_{ij k} h$.

In the following theorem, we show that the power of the proposed test $T^D_{\alpha}$ converges to one as $n$ and $p$ increasing to infinity under a newly specified sparse alternative.

\begin{theorem}\label{th32}
Given any $\gamma>0$ and a kernel $h(\cdot)$ satisfying Condition (C6), there exists some sufficiently large $B_3$ depending on $\gamma$ such that
$$
\inf _{F(\bmv)\in \mathcal{D}(\gamma, np ; h) \cap H_{\mathrm{a}}^{U}\left(B_{3}\right)} \mathbb{P}\left(T^D_{\alpha}=1\right)=1-o(1)
$$
\end{theorem}

The establishment of the rate optimality of the tests based on degenerate U-statistics requires the following assumption for the distribution:
\begin{itemize}
\item[(A4)] When $\bmv$ is Gaussian, suppose that for degenerate U-statistics $U_{ij}(k)$ and large $n$ and $p$, $c U_{ij}(k) \leqslant r_{ij}(k) \leqslant C U_{ij}(k)$ for $1\le i,j\le p, 1\le k\le K$ {\color{black}with probability tending to one}, where $c$ and $C$ are two constants.
\end{itemize}

Under the new type of sparse local alternative, we could show the rate optimality in terms of power for the proposed test in the following theorem.

\begin{theorem}\label{th33}
Suppose that Degenerate U-statistics $\left\{U_{ij}(k), 1 \leqslant i,j \leqslant p,1\le k\le K\right\}$ satisfy all the conditions in Theorems \ref{th31} and \ref{th32}. Suppose also that Assumption (A4) holds. Then,  the corresponding size-$\alpha$ test $T^D_{\alpha}$ is rate-optimal. In other
words, there exist two constants $D_{5}<D_{6}$ such that:
\begin{itemize}
\item[(i)] $\sup _{F(\bmv) \in \mathcal{D}(\gamma, np ; h) \cap H_{\mathrm{a}}^{U}\left(D_{6}\right)} \operatorname{pr}\left(T_{\alpha}=0\right)=o(1)$;
\item[(ii)] for any $\beta>0$ satisfying $\alpha+\beta<1,$ for large $n$ and $p$ we have
$$
\inf _{T_{\alpha} \in \mathcal{T}_{\alpha}} \sup _{F(\bmv) \in \mathcal{D}(\gamma, np ; h) \cap H_{\mathrm{a}}^{U}\left(D_{5}\right)} \operatorname{pr}\left(T_{\alpha}=0\right) \geqslant 1-\alpha-\beta
$$
\end{itemize}

\end{theorem}

Three examples belonging to the family of degenerate U-statistics are provided to test the high dimensional white noise as follows.

\paragraph{Example 3 (Hoeffding's $D$).} The Hoeffding's $D$ statistic is a rank-based U-statistic of order 5, which is based on the symmetric kernel
$$
\begin{array}{l}
h_{D}\left(z_{1}, \ldots, z_{5}\right):=\frac{1}{16} \sum_{\left(i_{1}, \ldots, i_{5}\right) \in \mathcal{P}_{5}} \\
{\left[\left\{I\left(z_{i_{1}, 1} \leq z_{i_{5}, 1}\right)-I\left(z_{i_{2}, 1} \leq z_{i_{5}, 1}\right)\right\}\left\{I\left(z_{i_{3}, 1} \leq z_{i_{5}, 1}\right)-I\left(z_{i_{4}, 1} \leq z_{i_{5}, 1}\right)\right\}\right]} \\
{\left[\left\{I\left(z_{i_{1}, 2} \leq z_{i_{5}, 2}\right)-I\left(z_{i_{2}, 2} \leq z_{i_{5}, 2}\right)\right\}\left\{I\left(z_{i_{3}, 2} \leq z_{i_{5}, 2}\right)-I\left(z_{i_{4}, 2} \leq z_{i_{5}, 2}\right)\right\}\right]}.
\end{array}
$$
Thus, the Hoeffding's $D$ correlation measure is given by $\mathbb{E} h_{D}$. Based on \citet[Proposition 7]{weihs2018symmetric} or \citet[Theorem 4.4]{nandy2016large}, under the measure $\mathbb{P}_{0} \otimes \mathbb{P}_{0}$, the eigenvalues and corresponding eigenfunctions of $h_{D, 2}(\cdot)$ are:
$$
\lambda_{i, j ; D}=3 /\left(\pi^{4} i^{2} j^{2}\right)>0, \quad i, j \in \mathbb{Z}^{+}
$$
and
$$
\phi_{i, j ; D}\left\{\left(z_{1,1}, z_{1,2}\right)^{\top}\right\}=2 \cos \left(\pi i z_{1,1}\right) \cos \left(\pi j z_{1,2}\right), \quad i, j \in \mathbb{Z}^{+},
$$
where $\Lambda_{D}:=\sum_{i, j} \lambda_{i, j ; D}=1 / 12$ and $\sup _{i, j}\left\|\phi_{i, j ; D}\right\|_{\infty} \leq 2$. Therefore, by considering the results in \citet{hoeffding1948non}, the kernel $h_{D}(\cdot)$ satisfies the three properties in Condition (C6).
Based on the result in \citet[p. 547]{hoeffding1948non}, the correlation measure $\mathbb{E} h_{D}$ is non-negative for arbitrary pair of random variables. Moreover, as shown by \citet{hoeffding1948non} and \citet{yanagimoto1970measures}, for a pair of random variables which is absolutely continuous in $\mathbb{R}^{2}$, the sufficient and necessary condition for their independence is that $\mathbb{E} h_{D}=0$. However, this result does not hold when the data is discrete or is continuous but not absolute continuous, e.g. a counter example is given in Remark 1 of \citet{yanagimoto1970measures}.

Define $\{\X_{tijk}\}=\{(\varepsilon_{t,i},\varepsilon_{t+k,j})^\top\}_{1\le t\le n-k}$.
According to (\ref{DEF_T_D_alpha}), the corresponding test is
\begin{align}
\widehat{D}_{ij}(k):=&\left(\begin{array}{c}
n-k \\
5
\end{array}\right)^{-1} \sum_{t_{1}<\cdots<t_{5}} h_{D}\left(\X_{t_{1},ijk}, \ldots, \X_{{t_{5}},ijk}\right)
\end{align}
and
\begin{align}
L_D:=&I\left\{ \max_{1\le k\le K} \max _{1\le i,j\le p} \frac{\pi^{4}(n-k-1)}{30} \widehat{D}_{ij}(k)-2 \log (Kp^2)+\log \log (Kp^2)+\frac{\pi^{4}}{36}>Q_{D, \alpha}\right\}
\end{align}
where $Q_{D, \alpha}:=\log \left\{\kappa_{D}^{2} /\pi\right\}-2 \log \log (1-\alpha)^{-1}$ and
$$
\kappa_{D}:=\left\{2 \prod_{n=2}^{\infty} \frac{\pi / n}{\sin (\pi / n)}\right\}^{1 / 2} \approx 2.467.
$$

\paragraph{Example 4 (Blum-Kiefer-Rosenblatt's $R$).} The Blum-Kiefer-Rosenblatt's $R$ statistic (\citet{blum1961distribution}) is a rank-based U-statistic of order 6, which is based on the symmetric kernel:
$$
\begin{array}{l}
h_{R}\left(z_{1}, \ldots, z_{6}\right):=\frac{1}{32} \sum_{\left(i_{1}, \ldots, i_{6}\right) \in \mathcal{P}_{6}} \\
{\left[\left\{I\left(z_{i_{1}, 1} \leq z_{i_{5}, 1}\right)-I\left(z_{i_{2}, 1} \leq z_{i_{5}, 1}\right)\right\}\left\{I\left(z_{i_{3}, 1} \leq z_{i_{5}, 1}\right)-I\left(z_{i_{4}, 1} \leq z_{i_{5}, 1}\right)\right\}\right]} \\
{\left[\left\{1\left(z_{i_{1}, 2} \leq z_{i_{6}, 2}\right)-I\left(z_{i_{2}, 2} \leq z_{i_{6}, 2}\right)\right\}\left\{I\left(z_{i_{3}, 2} \leq z_{i_{6}, 2}\right)-1\left(z_{i_{4}, 2} \leq z_{i_{6}, 2}\right)\right\}\right]}.
\end{array}
$$
The three properties in Condition (C6) can be easily verified based on the fact that $h_{R, 2}=2 h_{D, 2}$. Similarly, the correlation measure $\mathbb{E} h_{R}$ is non-negative for arbitrary pair of random variables. $\mathbb{E} h_{R}=0$ if and only if the pair of random variables are independent (without requiring the continuity properties), see, e.g. page 490 of \citet{blum1961distribution}.

According to (\ref{DEF_T_D_alpha}), the corresponding test is
$$
\widehat{R}_{ij}(k):=\left(\begin{array}{c}
n-k \\
6
\end{array}\right)^{-1} \sum_{t_{1}<\cdots<t_{6}} h_{R}\left(\X_{t_{1},ijk}, \ldots, \X_{t_{6},ijk}\right)
$$
and
$$
L_R:=I\left\{ \max_{1\le k\le K} \max _{1\le i,j\le p} \frac{\pi^{4}(n-k-1)}{90} \widehat{R}_{ij}(k)-2 \log (Kp^2)+\log \log (Kp^2)+\frac{\pi^{4}}{36}>Q_{R, \alpha}\right\}
$$
where $Q_{R, \alpha}:=Q_{D, \alpha}$

\paragraph{Example 5 (Bergsma-Dassios-Yanagimoto's $\tau^{*}$).} \citet{bergsma2014consistent} introduced a rank
correlation statistic as a U-statistic of order 4 with the symmetric kernel
\begin{align*}
& h_{\tau^{*}}\left(z_{1}, \ldots, z_{4}\right) \\
:=&\frac{1}{16}  \sum_{\left(i_{1}, \ldots, i_{4}\right) \in \mathcal{P}_{4}}\left\{1\left(z_{i_{1}, 1}, z_{i_{3}, 1}<z_{i_{2}, 1}, z_{i_{4}, 1}\right)+I\left(z_{i_{2}, 1}, z_{i_{4}, 1}<z_{i_{1}, 1}, z_{i_{3}, 1}\right)\right.\\
&\left.-I\left(z_{i_{1}, 1}, z_{i_{4}, 1}<z_{i_{2}, 1}, z_{i_{3}, 1}\right)-I\left(z_{i_{2}, 1}, z_{i_{3}, 1}<z_{i_{1}, 1}, z_{i_{4}, 1}\right)\right\}\\
&\left\{I\left(z_{i_{1}, 2}, z_{i_{3}, 2}<z_{i_{2}, 2}, z_{i_{4}, 2}\right)+I\left(z_{i_{2}, 2}, z_{i_{4}, 2}<z_{i_{1}, 2}, z_{i_{3}, 2}\right)\right.\\
&\left.-I\left(z_{i_{1}, 2}, z_{i_{4}, 2}<z_{i_{2}, 2}, z_{i_{3}, 2}\right)-I\left(z_{i_{2}, 2}, z_{i_{3}, 2}<z_{i_{1}, 2}, z_{i_{4}, 2}\right)\right\},
\end{align*}
where $I\left(y_{1}, y_{2}<y_{3}, y_{4}\right):=I\left(y_{1}<y_{3}\right) I\left(y_{1}<y_{4}\right) I\left(y_{2}<y_{3}\right) I\left(y_{2}<y_{4}\right)$. Based on the fact that $h_{\tau^{*}, 2}=3 h_{D, 2}$, all properties in Condition (C6) can be verified for $h_{\tau^{*}}(\cdot)$. As shown by Theorem 1 in \citet{bergsma2014consistent}, for a pair of random variables whose distribution is discrete, absolutely continuous, or a mixture of both, the correlation measure $\mathbb{E} h_{\tau^{*}}$ is non-negative and $\mathbb{E} h_{\tau^{*}}=0$ if and only if the pair is independent.

According to (\ref{DEF_T_D_alpha}), it yields the test
$$
\widehat{\tau}_{ij}^{*}(k):=\left(\begin{array}{c}
n-k \\
4
\end{array}\right)^{-1} \sum_{t_{1}<\cdots<t_{4}} h_{\tau^{*}}\left(\X_{t_{1},ijk}, \ldots, \X_{t_{4},ijk}\right)
$$
and
$$
L_{\tau^{*}}:=I\left\{ \max_{1\le k\le K} \max _{1\le i,j\le p}\frac{\pi^{4}(n-k-1)}{54}  \widehat{\tau}_{ij}^{*}(k)-2 \log (Kp^2)+\log \log (Kp^2)+\frac{\pi^{4}}{36}>Q_{\tau^{*}, \alpha}\right\}
$$
where $Q_{\tau^{*}, \alpha}:=Q_{D, \alpha}$

\section{Simulation}

In this section, we evaluate the empirical sizes and powers of several test statistics based on Monte Carlo simulation. We mainly compare the performance of the following test statistics:
\begin{itemize}
\item $L_r$: the max-type test statistic provided by \citet{chang2017testing};
\item $S_r$: the sum-type test statistic provided by \citet{li2019testing};
\item $L_{\rho}$: the Spearman's rho statistic defined in Example 1;
\item $L_{\tau}$: the Kendall's tau statistic defined in Example 2;
\item $L_D$: the Hoeffding's $D$ statistic defined in Example 3;
\item $L_R$: the Blum-Kiefer-Rosenblatt's $R$ statistic defined in Example 4;
\item $L_{\tau^{*}}$: the Bergsma-Dassios-Yanagimoto's $\tau^{*}$ statistic defined Example 5.
\end{itemize}
\subsection{Empirical sizes}
Let $\bmv_t=\A\z_t$. We consider the following four distribution for $\z_t$: (a) $\z_t\sim N(\bm 0, \I_p)$; (b) $\z_t=\bm w_t^{1/3}$ with  $\bm w_t\sim N(\bm 0, \I_p)$; (c) $\z_t=\bm w_t^{3}$ with  $\bm w_t\sim N(\bm 0, \I_p)$; (d) $\z_t=(z_{t1},\cdots,z_{tp})^\top$ with $z_{ti}\stackrel{i.i.d}{\sim} t(3)/\sqrt{3}$. For the Models (i)-(iv), we consider $\A=\bms^{1/2}$ with $\bms=(0.5^{|i-j|})_{1\le i,j\le p}$ and $\z_t$ follows the settings of (a)-(d). For the Models (v)-(viii), we set $\A=(a_{ij})_{1\le i,j\le p}$ with $a_{ij}\stackrel{i.i.d}{\sim} U(-1,1)$ and $\z_t$ follows the settings of (a)-(d).

{\color{black}Here we use parametric bootstrap method to determine the critical value of $L_r$;}
The empirical sizes of the seven test statistics listed above are reported in Tables \ref{tab:t1} for $K=2$. The simulation results with $K=4, 6$ are in the Supplementary Material. In each table, the results are summarized for Models (i)-(viii) with different combinations of $n$ and $p$, i.e., $n=100, 200$ and $p=30, 60, 120, 240$. It is easy to see that the degenerate U-statistics $L_D$, $L_R$ and $L_{\tau^{*}}$ can control the sizes very well in most of the cases. However, the empirical sizes of the sum-type statistic $S_r$, the simple linear rank statistic $L_{\rho}$ and non-degenerate U-statistic $L_{\tau}$ are a little smaller than the nominal level. {\color{black} The parametric bootstrap method proposed by Chang et al. (2017) can control the empirical sizes of $L_r$ in most cases. In the supplemental material, we also proposed a permutation procedure to control the empirical sizes of the above rank-based test statistics. Additional simulation studies show the good performance of the permutation procedure.}


\begin{table}[!ht]
\begin{center}
\footnotesize
\caption{\label{tab:t1} Sizes of tests with $K=2$ under Model (i)-(viii).}
                    \renewcommand\tabcolsep{1.2pt}
                     \renewcommand{\arraystretch}{1}
                     {
\begin{tabular}{cc|ccccccc|ccccccc}
\hline \hline
$n$&$p$&\multicolumn{7}{c}{  i}&\multicolumn{7}{c}{  ii}\\ \hline
&&$L_r$&$L_{\tau}$&$L_\rho$&$L_{\tau^*}$&$L_D$& $L_R$&$S_r$&$L_r$&$L_{\tau}$&$L_\rho$&$L_{\tau^*}$&$L_D$& $L_R$&$S_r$\\ \hline
100&30&0.043&0.013&0.014&0.054&0.067&0.05&0.02&0.045&0.016&0.015&0.046&0.065&0.04&0.027\\
100&60&0.036&0.01&0.008&0.052&0.073&0.041&0.018&0.039&0.016&0.011&0.047&0.076&0.031&0.015\\
100&120&0.034&0.007&0.006&0.029&0.073&0.023&0.002&0.042&0.009&0.007&0.033&0.068&0.026&0.004\\
100&240&0.044&0.009&0.009&0.034&0.081&0.027&0&0.043&0.009&0.005&0.04&0.082&0.028&0\\
200&30&0.035&0.017&0.016&0.045&0.05&0.04&0.034&0.056&0.016&0.015&0.04&0.056&0.035&0.042\\
200&60&0.038&0.022&0.023&0.045&0.054&0.046&0.025&0.042&0.014&0.016&0.034&0.044&0.03&0.027\\
200&120&0.039&0.012&0.008&0.032&0.056&0.025&0.015&0.046&0.015&0.018&0.04&0.044&0.034&0.008\\
200&240&0.043&0.008&0.01&0.037&0.06&0.032&0.005&0.047&0.018&0.014&0.051&0.067&0.046&0\\ \hline
&&\multicolumn{7}{c}{  iii}&\multicolumn{7}{c}{  iv}\\ \hline
100&30&0.022&0.009&0.007&0.03&0.052&0.026&0.085&0.038&0.011&0.012&0.044&0.061&0.034&0.038\\
100&60&0.024&0.006&0.004&0.034&0.05&0.027&0.089&0.042&0.008&0.005&0.028&0.057&0.022&0.05\\
100&120&0.032&0.01&0.011&0.044&0.087&0.029&0.055&0.061&0.011&0.009&0.055&0.097&0.042&0.057\\
100&240&0.036&0.033&0.01&0.035&0.086&0.03&0.034&0.054&0.009&0.008&0.04&0.088&0.028&0.035\\
200&30&0.051&0.014&0.015&0.041&0.056&0.04&0.08&0.037&0.018&0.019&0.044&0.055&0.046&0.071\\
200&60&0.035&0.015&0.014&0.04&0.054&0.035&0.08&0.043&0.019&0.018&0.048&0.059&0.045&0.075\\
200&120&0.036&0.016&0.011&0.046&0.062&0.039&0.086&0.047&0.016&0.017&0.039&0.066&0.039&0.056\\
200&240&0.042&0.021&0.016&0.041&0.062&0.035&0.034&0.036&0.014&0.012&0.045&0.068&0.041&0.062\\
\hline
&&\multicolumn{7}{c}{ v}&\multicolumn{7}{c}{vi}\\ \hline
100&30&0.044&0.018&0.017&0.049&0.067&0.041&0.031&0.051&0.017&0.012&0.049&0.067&0.043&0.033\\
100&60&0.038&0.015&0.013&0.05&0.082&0.044&0.015&0.043&0.011&0.012&0.042&0.066&0.039&0.016\\
100&120&0.035&0.008&0.006&0.034&0.073&0.026&0.001&0.036&0.011&0.007&0.035&0.076&0.025&0.004\\
100&240&0.044&0.012&0.009&0.045&0.097&0.036&0&0.029&0.009&0.003&0.045&0.091&0.027&0\\
200&30&0.057&0.013&0.012&0.031&0.046&0.031&0.055&0.061&0.018&0.022&0.054&0.063&0.044&0.035\\
200&60&0.037&0.009&0.008&0.038&0.053&0.034&0.038&0.043&0.014&0.013&0.043&0.06&0.043&0.029\\
200&120&0.044&0.018&0.017&0.049&0.075&0.043&0.021&0.039&0.016&0.012&0.046&0.061&0.04&0.016\\
200&240&0.037&0.008&0.009&0.032&0.065&0.034&0.008&0.052&0.014&0.016&0.052&0.068&0.048&0.003\\ \hline
&&\multicolumn{7}{c}{  vii}&\multicolumn{7}{c}{  viii}\\ \hline
100&30&0.041&0.01&0.009&0.037&0.056&0.031&0.09&0.053&0.013&0.015&0.055&0.072&0.05&0.047\\
100&60&0.038&0.011&0.009&0.045&0.07&0.037&0.063&0.041&0.010&0.006&0.041&0.075&0.035&0.057\\
100&120&0.037&0.017&0.008&0.046&0.087&0.036&0.055&0.045&0.008&0.008&0.035&0.075&0.026&0.04\\
100&240&0.053&0.009&0.005&0.034&0.085&0.022&0.023&0.038&0.009&0.005&0.042&0.094&0.028&0.037\\
200&30&0.061&0.024&0.024&0.053&0.063&0.053&0.087&0.048&0.011&0.013&0.032&0.04&0.03&0.067\\
200&60&0.036&0.02&0.017&0.049&0.066&0.046&0.075&0.045&0.015&0.016&0.046&0.058&0.046&0.058\\
200&120&0.048&0.021&0.02&0.05&0.077&0.042&0.064&0.039&0.019&0.016&0.055&0.081&0.047&0.063\\
200&240&0.049&0.01&0.011&0.039&0.056&0.036&0.041&0.054&0.02&0.023&0.041&0.063&0.036&0.056\\
\hline
\hline
\end{tabular}}
\end{center}
\end{table}

\subsection{Power comparison}
We consider the following eight examples as the data generation procedure in order to investigate the powers of different test statistics. Let $\z_t\sim N(\bm 0,\I_p)$. In the following, with slight abuse of notation, we write $f(v)=\left(f\left(v_{1}\right), \ldots, f\left(v_{p}\right)\right)^{\top}$ for any univariate function $f: \mathbb{R} \rightarrow \mathbb{R}$ and $v=\left(v_{1}, \ldots, v_{p}\right)^{\top} \in \mathbb{R}^{p}$. That is,
(I) $\bmv_t=\A\bmv_{t-1}+\z_t$;
(II) $\bmv_t=\sin(\frac{2\pi}{3}\A\bmv_{t-1})+\z_t$;
(III) $\bmv_t=\sin(\frac{\pi}{3}(\A\bmv_{t-1})^{1/3})+\z_t$;
(IV) $\bmv_t=(\A\bmv_{t-1})^{1/3}+\z_t$;
(V) $\bmv_t=\z_t+\A\z_{t-1}$;
(VI) $\bmv_t=\z_t+\sin(\frac{2\pi}{3}\A\z_{t-1})$;
(VII) $\bmv_t=\z_t+\sin(\frac{\pi}{3}(\A\z_{t-1})^{1/3})$;
(VIII) $\bmv_t=\z_t+(\A\z_{t-1})^{1/3}$.

We consider $\A=(a_{ij})_{1\le i,j\le p}$ with $a_{ij}\sim U(-\rho,\rho)$ if $1\le i,j\le k_0$ and $a_{ij}=0$ otherwise.
Models (I) and (V) can be classified as the linear relationship in autocorrelations, while Models (IV) and (VIII) can be classified as the monotone relationship in autocorrelations. The rest of the models are the non-linear and non-monotone relationships in autocorrelation.

Figures \ref{fig1} report the power curves with different $\rho$ for $K=2$. The power curves with $K=4,6$ are in the Supplementary Material. For these three figures, we set $k_0=2$, $n=100$ and $p=30$. It is clear that the parameter $\rho$ controls the level of the autocorrelation. Therefore, as $\rho$ increasing, the power curves show the upward trend as well for most of the models. Moreover, the power curves of the degenerate U-statistics $L_D$, $L_R$ and $L_{\tau^{*}}$ are higher than that of other max-type test statistics, i.e., $L_r$, $L_{\rho}$ and $L_{\tau}$. The sum-type test statistic $S_r$ has the lowest power curve in most cases. It is not surprising because $k_0$ was set as 2 here and sum-type test cannot work well under the sparse alternatives.

Figure \ref{fig4} shows the power curves with different $k_0$. For fixed $p$ and $\rho$, the parameter $k_0$ is used to control the sparsity of the autocorrelations. The higher value of $k_0$ yields the lower level of the sparsity in the autocorrelations. As expected, the power curves of the max-type test statistics, i.e., $L_r$, $L_{\rho}$, $L_{\tau}$, $L_D$, $L_R$ and $L_{\tau^{*}}$ have the downward trend when $k_0$ increasing in most of the models. Moreover, among the six max-type test statistics, the power curves of the degenerate U-statistics $L_D$, $L_R$ and $L_{\tau^{*}}$ are relatively higher than that of the other three test statistics. In contrast, the power curve of the sum-type test statistic $S_r$ has the upward trend as $k_0$ increasing in most models.

Figure \ref{fig5} shows the power curves with different $p$. For fixed $\rho$ and $k_0$, as the parameter $p$ increasing, the signal strength tends to decrease. Therefore, it is not surprising that all power curves of the seven test statistics show the downward trend as $p$ increasing. The power curves of degenerate U-statistics $L_D$, $L_R$ and $L_{\tau^{*}}$ are the highest in the seven test statistics. In contrast, the power curve of the sum-type test statistic $S_r$ is the lowest among the seven statistics.

\begin{figure}[p]
\centering
\includegraphics[width=1\textwidth]{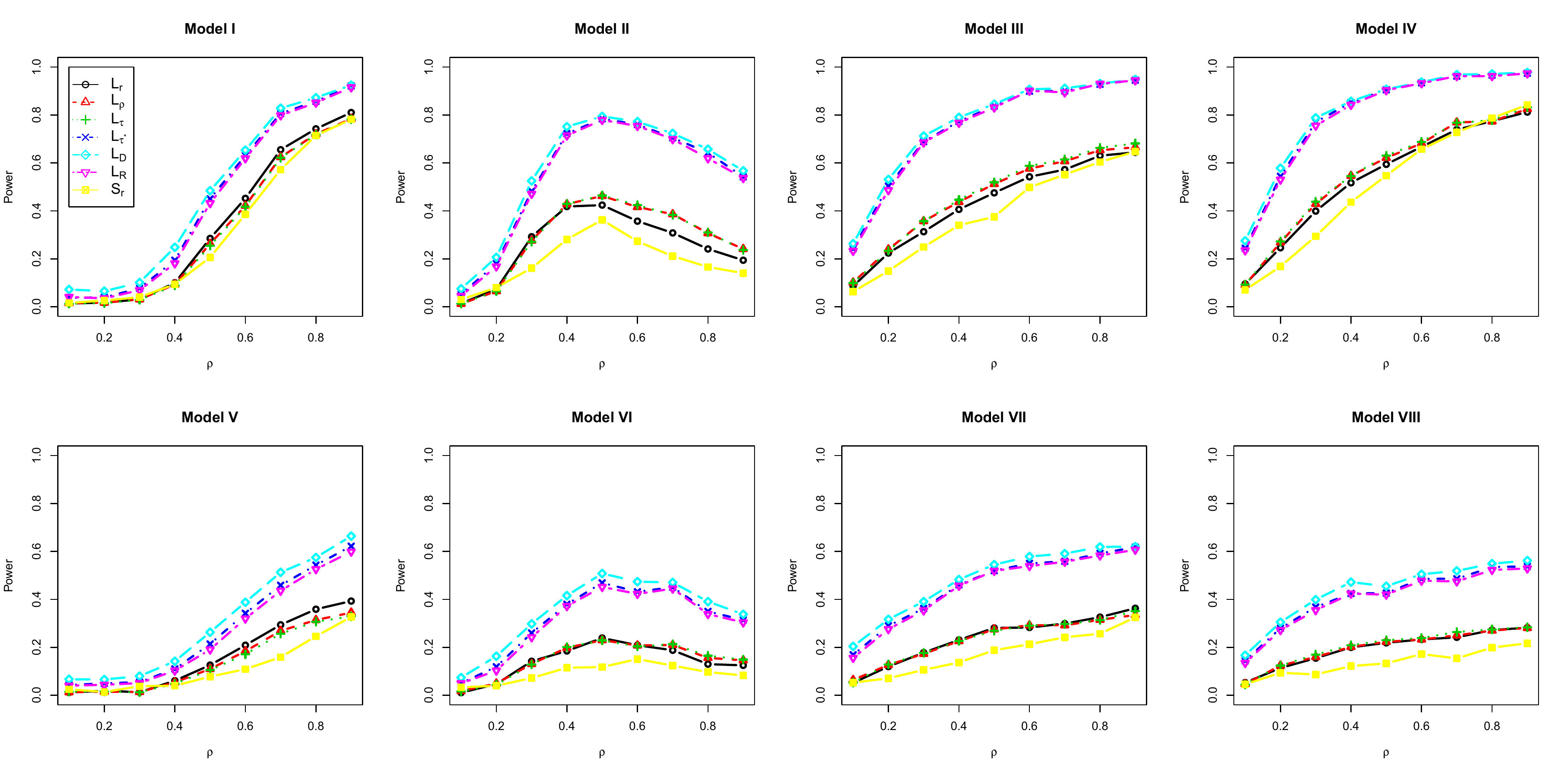}
\caption{\label{fig1}{Power curves of different methods  with different $\rho$ and $k_0=2, n=100,p=30, K=2$.}}
\end{figure}

\begin{figure}[p]
\centering
\includegraphics[width=1\textwidth]{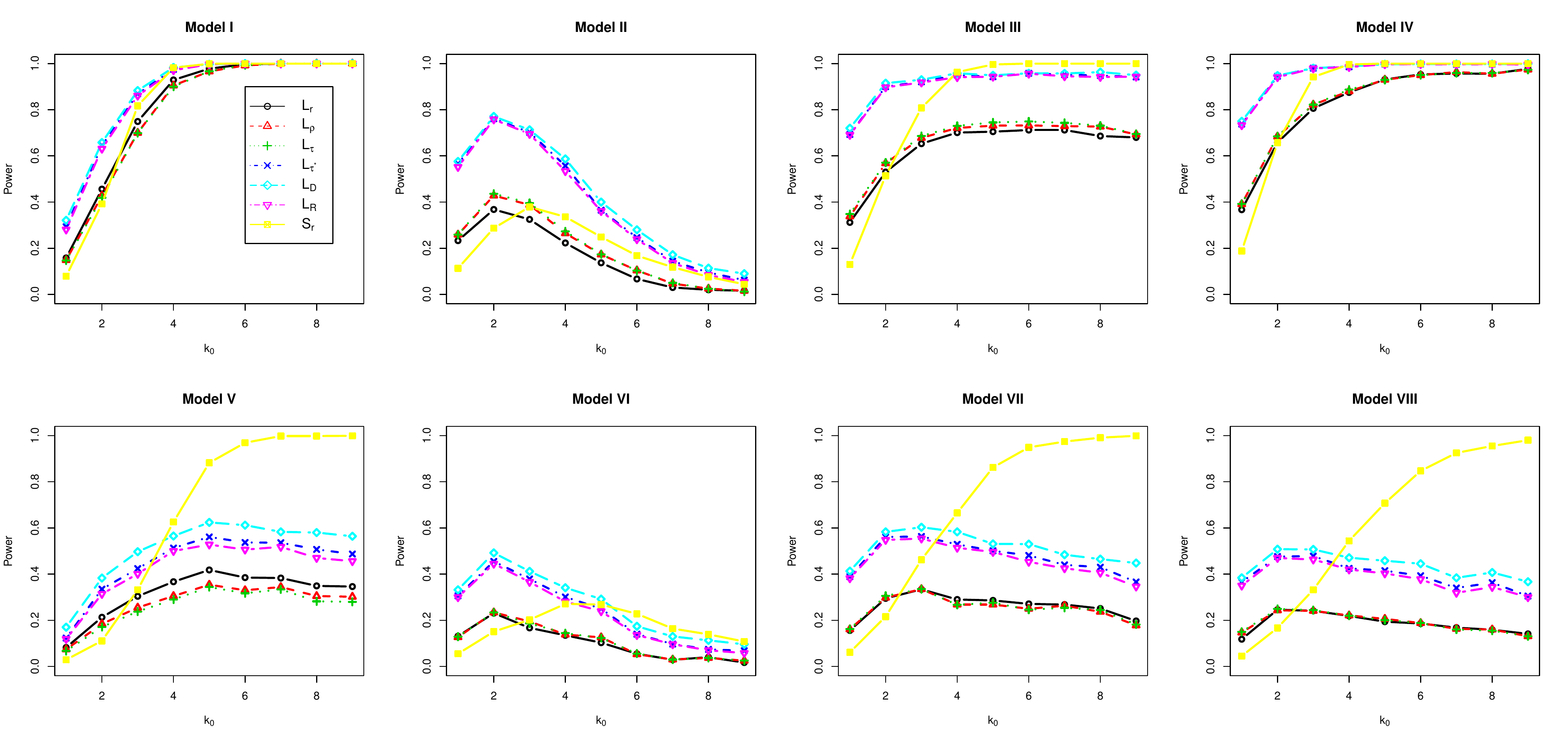}
\caption{\label{fig4}{Power curves  of different methods with different $k_0$ and $\rho=0.6, n=100,p=30, K=2$.}}
\end{figure}

\begin{figure}[ht]
\centering
\includegraphics[width=1\textwidth]{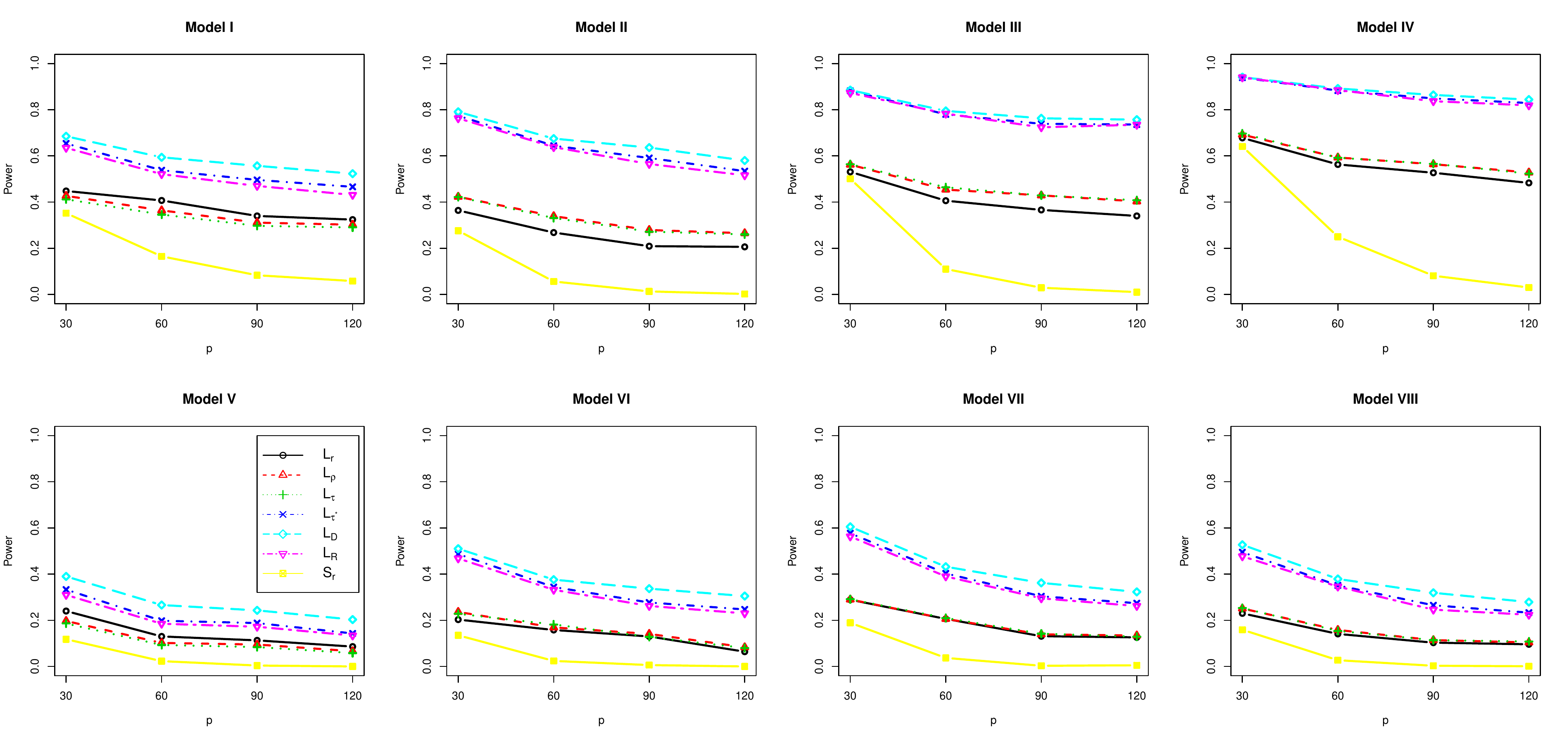}
\caption{\label{fig5}{Power curves  of different methods with different $p$ and $\rho=0.6, n=100,k_0=2, K=2$.}}
\end{figure}

\section{Conclusion}

To test the high-dimensional white noise, we develop the max-type tests based on three families of rank based statistics, including the simple linear rank statistics, non-degenerate U-statistics and degenerate U-statistics. The proposed tests are distribution free and in particular, the degenerate U-statistics can be used to detect the non-linear and non-monotone relationships in autocorrelations. Finally, as the theoretical contribution of this paper, we have relaxed the cross-sectional independence assumption in existing literature when deriving the asymptotic distributions for the rank correlation statistics. {\color{black} From the simulation studies, we found that the power of degenerate U-statistics $L_D, L_R, L_{\tau^*}$ have the best performance. So we suggest the degenerate U-statistics proposed in subsection 2.3 in practice.}

For the future directions related to the high-dimensional white noise test, it is also important to develop the theory for the sum-type tests based on the rank based statistics. The asymptotic independence between the max-type test and sum-type test based on the rank based statistics is also necessary to be established because of its usefulness in constructing some combination test which can be robust to both sparse and dense alternatives.

\section*{Acknowledgement}
The research of Dachuan Chen is supported by the National Natural Science Foundation of China (Grants 12101335 and 12271271), the Natural Science Foundation of Tianjin (Grant 21JCQNJC00020), the Fundamental Research Funds for the Central Universities, Nankai University (Grants 63211088, 63221050, and 63231013) and Wukong Investment Research Funds. Long Feng was partially supported by Shenzhen Wukong Investment Company, the Fundamental Research Funds for the Central Universities under Grant No. ZB22000105 and 63233075, the China National Key R\&D Program (Grant Nos. 2019YFC1908502, 2022YFA1003703, 2022YFA1003802, 2022YFA1003803) and the National Natural Science
Foundation of China Grants (Nos. 12271271, 11925106, 12231011, 11931001 and 11971247). Fengyi Song and Long Feng are co-corresponding authors and equally contributed to this paper.

\linespread{1.0}

\begin{footnotesize}

\bibhang=1.7pc
\bibsep=2pt
\fontsize{9}{14pt plus.8pt minus .6pt}\selectfont
\renewcommand\bibname{\large \bf References}
\expandafter\ifx\csname
natexlab\endcsname\relax\def\natexlab#1{#1}\fi
\expandafter\ifx\csname url\endcsname\relax
  \def\url#1{\texttt{#1}}\fi
\expandafter\ifx\csname urlprefix\endcsname\relax\def\urlprefix{URL}\fi



\vskip .65cm
\noindent
School of Statistics and Data Science, KLMDASR, LEBPS, and LPMC, Nankai University
\vskip 2pt
\noindent
E-mail: dchen@nankai.edu.cn
\vskip 2pt

\noindent
School of Statistics and Data Science, KLMDASR, LEBPS, and LPMC, Nankai University
\vskip 2pt
\noindent
E-mail: sauntbai@163.com
\vskip 2pt

\noindent
School of Statistics and Data Science, KLMDASR, LEBPS, and LPMC, Nankai University
\vskip 2pt
\noindent
E-mail: flnankai@nankai.edu.cn

\end{footnotesize}

\newpage

\centerline{\large\bf Supplement to}
\vspace{0.5cm}
\centerline{\large\bf ``Rank Based Tests for High Dimensional White Noise"}
\vspace{.4cm}
\centerline{Dachuan Chen, Fengyi Song and Long Feng}
\vspace{.4cm}
\centerline{\it Nankai University}

{\color{black}
\section{Chatterjee's rank Correlation}
For a pair of continuous random variables $(X,Y)$, Chatterjee (2021) proposed a new rank correlation, i.e.
\begin{align*}
\xi_n(X,Y)=1-\frac{3\sum_{i=1}^{n-1}|r_{i+1}-r_i|}{n^2-1}
\end{align*}
where $r_i$ is the rank of $Y_{(i)}$. Here we rearrange the data as $(X_{(1)},Y_{(1)}),\cdots,(X_{(n)},Y_{(n)})$ such that $X_{(1)}\le \cdots\le X_{(n)}$.
So, we can also consider the test based on Chatterjee's rank correlation of the form
\begin{align}
\Xi_{ij}(k)=1-\frac{3\sum_{t=1}^{n-k}|R_{n-k,t+k+1}^{ij}(k)-R_{n-k,t+k}^{ij}(k)|}{(n-k+1)^2-1}
\end{align}
By Theorem 2.2 in Chatterjee (2021), $\sqrt{n-k+1}\Xi_{ij}(k)\cd N(0,2/5)$ as $n\to \infty$  under the null hypothesis. Thus, we propose the following statistics for testing $H_0$:
\begin{align}
\Xi_n=\max_{1\le k\le K}\max_{1\le i,j\le p}\sqrt{\frac{5(n-k+1)}{2}}|\Xi_{ij}(k)|
\end{align}

Next, we state the theoretical result about the limiting null distribution of $\Xi_n$.
\begin{theorem}
If $(\varepsilon_{t1},\cdots,\varepsilon_{tp})$ are mutually independent continuous random variables, under $H_0$, for any $y\in \mathbb{R}$, we have
\begin{align*}
P\left(\Xi_n^2-2\log(N)+\log\log(N)\le y\right)=\exp \left\{-\pi^{-1 / 2} \exp (-y / 2)\right\}+o(1)
\end{align*}
as $n,p\to \infty$.
\end{theorem}

Based on Theorem S1, we proposed the following high dimensional white noise test based on Chatterjee's Correlation
\begin{equation}
T_\alpha^\xi\doteq I\left(\Xi_{n}^{2}-2\log(Kp^2)+\log \log (Kp^2)\ge q_\alpha\right), \label{DEF_Test_Tsr}
\end{equation}
where $q_\alpha=-\log (\pi)-2\log\log(1-\alpha)^{-1}$.

Table \ref{tcrs} show the empirical sizes of the proposed test $T_{\alpha}^\xi$ under the same settings as subsection 3.1 in the main text. We observe that the empirical sizes of $T_{\alpha}^\xi$ is a little conservative in most cases. Additionally, we also show the power of the proposed test $T_{\alpha}^\xi$ and the seven tests in the main text in Table \ref{tcrp} with $n=100,p=30,\rho=0.5, k_0=2$. The other settings are all the same as subsection 3.2 in the main text. We found that $T_{\alpha}^\xi$ does not perform very well in most cases, which is consistent with many recent studies (Cao and Bickel, 2020; Shi et al., 2021; Lin and Han, 2023). They showed that independence tests based on Chatterjee’s rank correlation are unfortunately rate-inefficient against various local alternatives.

\begin{table}[!ht]
\begin{center}
\footnotesize
\caption{\label{tcrs} Empirical sizes of $L_{\xi}$ under Models (i)-(viii). ($L_{\xi}$: the max-type test defined in (\ref{DEF_Test_Tsr}).)}
                    \renewcommand\tabcolsep{5pt}
                     \renewcommand{\arraystretch}{1}
                     {
\begin{tabular}{cc|cccccccc}
\hline \hline
&Models&i&ii&iii&iv&v&vi&vii&viii\\ \hline
$n$&$p$&\multicolumn{8}{c}{$K=2$}\\ \hline
100&30&0.022&0.03&0.026&0.029&0.026&0.032&0.022&0.019\\
100&60&0.026&0.019&0.028&0.026&0.027&0.02&0.019&0.024\\
100&120&0.032&0.019&0.023&0.023&0.027&0.025&0.031&0.022\\
100&240&0.031&0.018&0.02&0.027&0.03&0.022&0.023&0.019\\
200&30&0.033&0.025&0.039&0.03&0.039&0.03&0.038&0.035\\
200&60&0.033&0.029&0.031&0.03&0.029&0.032&0.035&0.028\\
200&120&0.031&0.030&0.029&0.029&0.035&0.037&0.029&0.038\\
200&240&0.030&0.035&0.039&0.027&0.025&0.033&0.038&0.027\\\hline
&&\multicolumn{8}{c}{$K=4$}\\ \hline
100&30&0.029&0.022&0.023&0.022&0.025&0.031&0.032&0.016\\
100&60&0.016&0.022&0.021&0.025&0.02&0.031&0.028&0.022\\
100&120&0.018&0.023&0.032&0.015&0.029&0.019&0.029&0.028\\
100&240&0.03&0.023&0.03&0.029&0.018&0.02&0.031&0.031\\
200&30&0.036&0.033&0.038&0.033&0.04&0.032&0.032&0.033\\
200&60&0.04&0.03&0.023&0.032&0.023&0.037&0.028&0.023\\
200&120&0.022&0.033&0.023&0.032&0.02&0.028&0.034&0.023\\
200&240&0.02&0.038&0.035&0.036&0.028&0.024&0.022&0.02\\ \hline
&&\multicolumn{8}{c}{$K=6$}\\ \hline
100&30&0.03&0.027&0.022&0.023&0.027&0.018&0.016&0.028\\
100&60&0.015&0.028&0.028&0.032&0.03&0.03&0.032&0.018\\
100&120&0.021&0.022&0.021&0.022&0.029&0.021&0.018&0.021\\
100&240&0.021&0.027&0.025&0.024&0.016&0.016&0.022&0.029\\
200&30&0.023&0.028&0.02&0.013&0.026&0.027&0.024&0.032\\
200&60&0.012&0.022&0.022&0.012&0.02&0.019&0.017&0.018\\
200&120&0.021&0.018&0.02&0.012&0.02&0.011&0.019&0.018\\
200&240&0.023&0.021&0.012&0.010&0.022&0.012&0.025&0.022\\
\hline
\hline
\end{tabular}}
\end{center}
\end{table}

\begin{table}[!ht]
\begin{center}
\footnotesize
\caption{\label{tcrp} Power of tests with $n=100,p=30,\rho=0.5, k_0=2$ under Models (I)-(VIII). ($L_{\xi}$: the max-type test defined in (\ref{DEF_Test_Tsr}).)}
                    \renewcommand\tabcolsep{5pt}
                     \renewcommand{\arraystretch}{1}
                     {
\begin{tabular}{c|cccccccc}
\hline \hline
&\multicolumn{8}{c}{Methods}\\ \hline
Models&$L_r$&$L_{\tau}$&$L_\rho$&$L_{\tau^*}$&$L_D$& $L_R$&$S_r$&$L_{\xi}$\\ \hline
&\multicolumn{8}{c}{$K=2$}\\ \hline
I&0.22&0.22&0.22&0.38&0.42&0.37&0.22&0.03\\
II&0.41&0.43&0.44&0.81&0.83&0.79&0.43&0.02\\
III&0.45&0.48&0.48&0.79&0.8&0.8&0.41&0.06\\
IV&0.64&0.63&0.64&0.88&0.9&0.88&0.49&0.21\\
V&0.15&0.1&0.12&0.26&0.27&0.25&0.05&0.02\\
VI&0.29&0.26&0.23&0.45&0.46&0.44&0.1&0.01\\
VII&0.38&0.33&0.37&0.54&0.55&0.55&0.16&0.06\\
VIII&0.26&0.24&0.23&0.42&0.49&0.43&0.17&0.02\\ \hline
&\multicolumn{8}{c}{$K=4$}\\ \hline
I&0.21&0.21&0.2&0.37&0.42&0.34&0.03&0.02\\
II&0.42&0.46&0.44&0.7&0.75&0.71&0.12&0.08\\
III&0.44&0.42&0.42&0.76&0.81&0.75&0.16&0.03\\
IV&0.58&0.56&0.56&0.83&0.86&0.83&0.28&0.1\\
V&0.11&0.06&0.06&0.21&0.26&0.19&0&0.04\\
VI&0.17&0.17&0.17&0.37&0.38&0.36&0.03&0\\
VII&0.18&0.18&0.18&0.47&0.51&0.45&0.06&0.01\\
VIII&0.13&0.16&0.16&0.39&0.42&0.38&0.04&0.03\\ \hline
&\multicolumn{8}{c}{$K=6$}\\ \hline
I&0.19&0.17&0.17&0.36&0.42&0.35&0.04&0.05\\
II&0.31&0.39&0.42&0.71&0.72&0.68&0.08&0.07\\
III&0.42&0.45&0.45&0.75&0.78&0.75&0.02&0.01\\
IV&0.49&0.51&0.53&0.87&0.88&0.87&0.12&0.12\\
V&0.1&0.07&0.06&0.13&0.19&0.12&0&0.04\\
VI&0.07&0.11&0.09&0.26&0.3&0.24&0&0.03\\
VII&0.2&0.19&0.2&0.38&0.38&0.37&0.02&0.04\\
VIII&0.14&0.18&0.16&0.31&0.34&0.31&0.01&0.02\\
\hline
\hline
\end{tabular}}
\end{center}
\end{table}

\section{L-statistic}
As shown in the main text, the max-type test statistics performs very well under sparse alternative. Motivated by \citet{cjs2022}, we consider an L-statistic for high dimensional white noise test, which combines the first several largest signals together. That is,
\begin{align}
L_V=\sum_{l=1}^L V_{(l)}, ~~ L_U=\sum_{l=1}^L U_{(l)}
\end{align}
where $V_{(l)}$ and $U_{(l)}$ are the $l$-th largest maximum of $\{|V_{ij}(k)|\}_{1\le i,j\le p, 1\le k\le K}$ and $\{|U_{ij}(k)|\}_{1\le i,j\le p, 1\le k\le K}$, respectively.

It is difficult to establish the limit null distribution of $L_V$ and $L_U$. So we adopt the permutation test to calculate the critical value of each test. We randomly rearrange $\{\varepsilon_1,\cdots,\varepsilon_T\}$ as $\{\varepsilon_{\pi(1)},\cdots,\varepsilon_{\pi(T)}\}$ where $\pi$ is a permutation of $\{1,\cdots,T\}$. The permutation test statistic $\tilde{L}_V$ and $\tilde{L}_U$ are accordingly built from the permutation sample $\{\varepsilon_{\pi(1)},\cdots,\varepsilon_{\pi(T)}\}$. When this procedure is
repeated many times, the permutation critical value $z_{\alpha}^L$ and $z_{\alpha}^U$  are the empirical $1-\alpha$ quantile of the permutation test statistic, respectively. The tests with rejection region $L_V\ge z_{\alpha}^L$ and $L_U\ge z_{\alpha}^U$ are our proposal.

Here we give a simulation study of the $L$-statistics. Let $\tilde L_{\tau}$, $\tilde L_\rho$, $\tilde L_{\tau^*}$, $\tilde L_D$, $\tilde L_R$ denote the corresponding $L$-statistics based on Kendall's tau, Spearman's rho, Bergsma-Dassios-Yanagimoto's $\tau^{*}$, Hoeffding's $D$, Blum-Kiefer-Rosenblatt's $R$, respectively. Table \ref{tab:st1} reports the empirical sizes of the above test statistics with $K=2$ and different $L$ under Model (i). We found that the permutation procedure can control all the empirical sizes of these tests very well. To show the performance of $L$-statistics, we consider a power comparison of $\tilde{L}_{\tau^*}$ with different $L$ because the Bergsma-Dassios-Yanagimoto's $\tau^{*}$ statistic performs very well in most cases in the simulation studies in the main text. We consider the same settings as subsection 3.2 in the main text except that $k_0=1,\cdots,10$ and $\rho=0.98k_0^{2/3}$. From Figure \ref{figs1}, we observe that when the number of non-zero autocorrelations is small, $\tilde{L}_{\tau^*}$ with small $L$s have better performance than large $L$s and vice versa. So the optimal $L$ depends on the sparsity of the autocorrelations. Generally speaking, $\tilde{L}_{\tau^*}$ with $L>1$ outperforms $\tilde{L}_{\tau^*}$ with $L=1$, i.e. $L_{\tau^*}$ in most cases. So how to derive the limit null distribution of rank based $L$-statistics and choose the optimal $L$ for high dimensional white noise test deserves further studies.

\begin{table}[!ht]
\begin{center}
\footnotesize
\caption{\label{tab:st1} Sizes of $L$-statistic tests with $K=2$ under Model (i).}
                    \renewcommand\tabcolsep{4pt}
                     \renewcommand{\arraystretch}{1}
                     {
\begin{tabular}{c|cccccccccc}
\hline \hline
&\multicolumn{10}{c}{ $L$}\\ \hline
&1&2&3&4&5&6&7&8&9&10\\ \hline
&\multicolumn{10}{c}{ $(n,p)=(100,30)$}\\ \hline
$\tilde L_{\tau}$&0.041&0.046&0.059&0.052&0.046&0.058&0.05&0.052&0.059&0.048\\
$\tilde L_\rho$&0.046&0.059&0.056&0.048&0.052&0.059&0.042&0.05&0.05&0.054\\
$\tilde L_{\tau^*}$&0.059&0.046&0.05&0.043&0.047&0.045&0.047&0.047&0.056&0.05\\
$\tilde L_D$&0.046&0.045&0.05&0.05&0.059&0.043&0.05&0.058&0.045&0.043\\
$\tilde L_R$&0.054&0.051&0.059&0.052&0.056&0.048&0.045&0.056&0.054&0.045\\ \hline
&\multicolumn{10}{c}{ $(n,p)=(100,60)$}\\ \hline
$\tilde L_{\tau}$&0.054&0.051&0.043&0.05&0.04&0.054&0.054&0.04&0.046&0.057\\
$\tilde L_\rho$&0.044&0.04&0.057&0.047&0.048&0.051&0.04&0.054&0.055&0.048\\
$\tilde L_{\tau^*}$&0.04&0.06&0.048&0.041&0.052&0.049&0.056&0.045&0.056&0.058\\
$\tilde L_D$&0.04&0.059&0.056&0.049&0.044&0.049&0.058&0.042&0.041&0.058\\
$\tilde L_R$&0.046&0.044&0.056&0.058&0.047&0.04&0.048&0.05&0.043&0.05\\
\hline
&\multicolumn{10}{c}{ $(n,p)=(100,90)$}\\ \hline
$\tilde L_{\tau}$&0.047&0.041&0.053&0.048&0.058&0.049&0.048&0.043&0.041&0.049\\
$\tilde L_\rho$&0.047&0.057&0.042&0.057&0.046&0.058&0.055&0.043&0.059&0.052\\
$\tilde L_{\tau^*}$&0.058&0.053&0.041&0.048&0.053&0.049&0.045&0.053&0.051&0.041\\
$\tilde L_D$&0.054&0.054&0.052&0.05&0.044&0.057&0.047&0.05&0.06&0.06\\
$\tilde L_R$&0.051&0.057&0.06&0.056&0.057&0.054&0.049&0.057&0.058&0.049\\ \hline
&\multicolumn{10}{c}{ $(n,p)=(100,120)$}\\ \hline
$\tilde L_{\tau}$&0.052&0.049&0.049&0.057&0.049&0.043&0.041&0.048&0.045&0.04\\
$\tilde L_\rho$&0.058&0.054&0.044&0.049&0.056&0.058&0.06&0.049&0.049&0.045\\
$\tilde L_{\tau^*}$&0.044&0.05&0.059&0.059&0.049&0.06&0.055&0.05&0.052&0.049\\
$\tilde L_D$&0.058&0.059&0.044&0.044&0.054&0.041&0.05&0.045&0.041&0.044\\
$\tilde L_R$&0.06&0.05&0.053&0.055&0.046&0.049&0.042&0.045&0.047&0.042\\
\hline\hline
\end{tabular}}
\end{center}
\end{table}

\begin{figure}[p]
\centering
\includegraphics[width=1\textwidth]{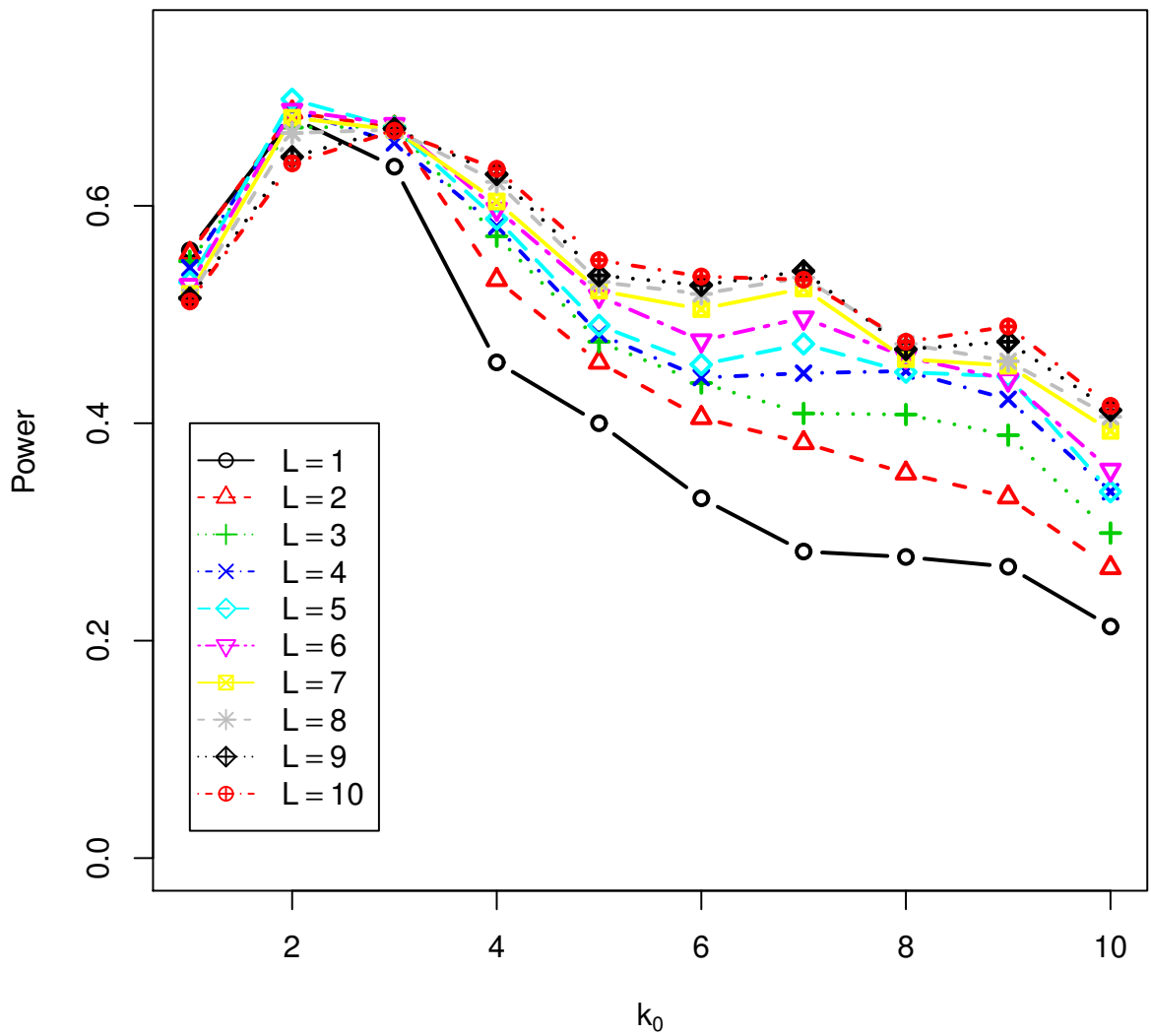}
\caption{\label{figs1}{Power curves of different $L$-statistics with different $L$ and $n=100,p=30, K=2$ under Model I.}}
\end{figure}

\section{Additional Simulation Results of Section 3 in the main document}

The empirical sizes of the seven test statistics listed in the beginning of Section 3 in the main document are reported in Table \ref{tab:t2} and \ref{tab:t3}, respectively. Figures \ref{fig2} and \ref{fig3} report the power curves of the seven test statistics with different $\rho$ for $K=4$ and $6$, respectively.
\begin{table}[!ht]
\begin{center}
\footnotesize
\caption{\label{tab:t2} Sizes of tests with $K=4$ under Model (i)-(viii).}
                    \renewcommand\tabcolsep{1.2pt}
                     \renewcommand{\arraystretch}{1}
                     {
\begin{tabular}{cc|ccccccc|ccccccc}
\hline \hline
$n$&$p$&\multicolumn{7}{c}{  i}&\multicolumn{7}{c}{  ii}\\ \hline
&&$L_r$&$L_{\tau}$&$L_\rho$&$L_{\tau^*}$&$L_D$& $L_R$&$S_r$&$L_r$&$L_{\tau}$&$L_\rho$&$L_{\tau^*}$&$L_D$& $L_R$&$S_r$\\ \hline
100&30 &0.042&0.016&0.016&0.055&0.087&0.047&0.013&0.038&0.007&0.007&0.04&0.074&0.037&0.014\\
100&60 &0.037&0.012&0.008&0.04&0.07&0.036&0.001&0.031&0.006&0.005&0.031&0.073&0.03&0.002\\
100&120&0.036&0.011&0.008&0.036&0.09&0.026&0&0.027&0.009&0.004&0.04&0.075&0.03&0\\
100&240&0.029&0.009&0.001&0.037&0.091&0.018&0&0.024&0.006&0.004&0.028&0.08&0.021&0\\
200&30 &0.047&0.017&0.014&0.047&0.055&0.041&0.022&0.051&0.016&0.015&0.036&0.047&0.034&0.034\\
200&60 &0.031&0.016&0.013&0.036&0.048&0.028&0.012&0.037&0.014&0.01&0.038&0.058&0.035&0.012\\
200&120&0.029&0.013&0.013&0.042&0.057&0.038&0.002&0.033&0.015&0.016&0.047&0.066&0.04&0.001\\
200&240&0.031&0.018&0.014&0.038&0.059&0.036&0&0.028&0.008&0.005&0.045&0.066&0.033&0\\ \hline
&&\multicolumn{7}{c}{  iii}&\multicolumn{7}{c}{  iv}\\ \hline
100&30&0.047&0.011&0.011&0.039&0.06&0.035&0.06&0.052&0.008&0.006&0.04&0.071&0.036&0.036\\
100&60&0.036&0.013&0.009&0.039&0.077&0.034&0.044&0.045&0.009&0.009&0.056&0.092&0.041&0.033\\
100&120&0.038&0.01&0.008&0.031&0.081&0.025&0.028&0.047&0.011&0.009&0.039&0.094&0.03&0.023\\
100&240&0.028&0.004&0.003&0.027&0.082&0.02&0.008&0.051&0.006&0.003&0.037&0.104&0.021&0.015\\
200&30&0.031&0.017&0.018&0.04&0.049&0.036&0.06&0.057&0.017&0.016&0.045&0.057&0.037&0.059\\
200&60&0.042&0.014&0.012&0.038&0.047&0.034&0.056&0.053&0.015&0.012&0.043&0.057&0.032&0.045\\
200&120&0.034&0.014&0.008&0.037&0.057&0.036&0.046&0.038&0.014&0.018&0.041&0.056&0.037&0.035\\
200&240&0.030&0.01&0.01&0.043&0.074&0.033&0.012&0.048&0.012&0.012&0.045&0.072&0.041&0.029\\ \hline
&&\multicolumn{7}{c}{ v}&\multicolumn{7}{c}{vi}\\ \hline
100&30&0.039&0.011&0.01&0.038&0.071&0.041&0.013&0.039&0.013&0.011&0.039&0.072&0.029&0.026\\
100&60&0.041&0.013&0.014&0.055&0.09&0.047&0.003&0.033&0.007&0.005&0.034&0.067&0.03&0.001\\
100&120&0.031&0.016&0.006&0.044&0.086&0.036&0.001&0.032&0.013&0.007&0.037&0.083&0.027&0\\
100&240&0.027&0.006&0.004&0.025&0.08&0.019&0&0.034&0.008&0.005&0.035&0.106&0.029&0\\
200&30&0.053&0.024&0.023&0.049&0.064&0.043&0.027&0.052&0.017&0.016&0.047&0.051&0.042&0.031\\
200&60&0.044&0.013&0.015&0.029&0.053&0.03&0.016&0.039&0.014&0.014&0.043&0.06&0.038&0.017\\
200&120&0.032&0.023&0.021&0.053&0.077&0.049&0.011&0.031&0.015&0.013&0.047&0.066&0.038&0.003\\
200&240&0.029&0.008&0.006&0.029&0.067&0.024&0&0.029&0.011&0.009&0.032&0.053&0.029&0\\ \hline
&&\multicolumn{7}{c}{  vii}&\multicolumn{7}{c}{  viii}\\ \hline
100&30&0.053&0.008&0.007&0.038&0.059&0.034&0.052&0.053&0.01&0.013&0.043&0.07&0.04&0.037\\
100&60&0.042&0.01&0.007&0.042&0.071&0.038&0.041&0.044&0.012&0.01&0.039&0.067&0.026&0.031\\
100&120&0.037&0.012&0.006&0.035&0.083&0.02&0.023&0.031&0.009&0.004&0.028&0.073&0.019&0.016\\
100&240&0.026&0.009&0.003&0.036&0.101&0.025&0.006&0.043&0.005&0.005&0.028&0.108&0.017&0.017\\
200&30&0.043&0.025&0.025&0.052&0.062&0.049&0.073&0.046&0.016&0.014&0.045&0.057&0.044&0.047\\
200&60&0.036&0.016&0.015&0.046&0.057&0.041&0.046&0.033&0.011&0.01&0.054&0.058&0.047&0.024\\
200&120&0.031&0.01&0.008&0.034&0.05&0.032&0.022&0.043&0.015&0.012&0.043&0.066&0.033&0.036\\
200&240&0.028&0.012&0.01&0.04&0.059&0.042&0.005&0.034&0.015&0.024&0.039&0.053&0.026&0.034\\
\hline
\hline
\end{tabular}}
\end{center}
\end{table}

\begin{table}[!ht]
\begin{center}
\footnotesize
\caption{\label{tab:t3} Sizes of tests with $K=6$ under Model (i)-(viii).}
                    \renewcommand\tabcolsep{1.2pt}
                     \renewcommand{\arraystretch}{1}
                     {
\begin{tabular}{cc|ccccccc|ccccccc}
\hline \hline
$n$&$p$&\multicolumn{7}{c}{  i}&\multicolumn{7}{c}{  ii}\\ \hline
&&$L_r$&$L_{\tau}$&$L_\rho$&$L_{\tau^*}$&$L_D$& $L_R$&$S_r$&$L_r$&$L_{\tau}$&$L_\rho$&$L_{\tau^*}$&$L_D$& $L_R$&$S_r$\\ \hline
100&30&0.049&0.016&0.012&0.051&0.082&0.046&0.002&0.055&0.013&0.014&0.041&0.066&0.038&0.004\\
100&60&0.047&0.009&0.006&0.038&0.078&0.032&0.001&0.053&0.012&0.01&0.046&0.099&0.035&0\\
100&120&0.036&0.016&0.011&0.053&0.105&0.038&0&0.043&0.008&0.003&0.028&0.083&0.027&0 \\
100&240&0.034&0.005&0.002&0.023&0.077&0.014&0&0.041&0.007&0.002&0.026&0.095&0.018&0\\
200&30&0.055&0.011&0.013&0.043&0.047&0.04&0.015&0.045&0.015&0.015&0.036&0.046&0.032&0.013\\
200&60&0.056&0.016&0.014&0.037&0.051&0.036&0.002&0.042&0.012&0.01&0.046&0.054&0.046&0.008\\
200&120&0.037&0.009&0.007&0.045&0.062&0.036&0&0.041&0.014&0.012&0.032&0.053&0.032&0\\
200&240&0.039&0.019&0.013&0.043&0.071&0.04&0&0.036&0.01&0.007&0.046&0.078&0.041&0\\\hline
&&\multicolumn{7}{c}{  iii}&\multicolumn{7}{c}{  iv}\\ \hline
100&30&0.058&0.01&0.01&0.041&0.071&0.034&0.043&0.046&0.011&0.007&0.037&0.074&0.024&0.021\\
100&60&0.046&0.008&0.007&0.036&0.071&0.027&0.02&0.048&0.01&0.007&0.053&0.091&0.038&0.015\\
100&120&0.048&0.009&0.005&0.038&0.089&0.029&0.009&0.039&0.011&0.008&0.03&0.081&0.023&0.009\\
100&240&0.039&0.007&0.006&0.023&0.087&0.018&0.003&0.041&0.01&0.007&0.034&0.096&0.021&0.001\\
200&30&0.042&0.009&0.009&0.034&0.05&0.035&0.053&0.038&0.017&0.019&0.044&0.058&0.039&0.04\\
200&60&0.048&0.015&0.012&0.039&0.055&0.031&0.043&0.065&0.014&0.012&0.038&0.055&0.035&0.036\\
200&120&0.051&0.01&0.006&0.037&0.058&0.032&0.012&0.057&0.008&0.007&0.033&0.056&0.03&0.028\\
200&240&0.054&0.021&0.018&0.045&0.076&0.042&0.003&0.049&0.013&0.01&0.031&0.06&0.03&0.014\\\hline
&&\multicolumn{7}{c}{ v}&\multicolumn{7}{c}{vi}\\ \hline
100&30&0.049&0.014&0.01&0.036&0.06&0.032&0.007&0.049&0.018&0.01&0.05&0.082&0.046&0.007\\
100&60&0.038&0.01&0.007&0.043&0.072&0.033&0.002&0.036&0.01&0.01&0.037&0.076&0.029&0\\
100&120&0.032&0.012&0.007&0.033&0.086&0.026&0&0.032&0.009&0.008&0.026&0.069&0.022&0\\
100&240&0.031&0.011&0.006&0.028&0.093&0.021&0&0.033&0.008&0.003&0.035&0.096&0.019&0\\
200&30&0.042&0.019&0.014&0.044&0.065&0.04&0.02&0.047&0.017&0.014&0.038&0.047&0.036&0.021\\
200&60&0.041&0.016&0.014&0.051&0.071&0.043&0.006&0.035&0.014&0.011&0.044&0.06&0.037&0.001\\
200&120&0.032&0.014&0.011&0.034&0.053&0.03&0&0.041&0.006&0.005&0.033&0.049&0.031&0\\
200&240&0.037&0.012&0.007&0.035&0.065&0.029&0&0.037&0.013&0.013&0.04&0.073&0.032&0\\ \hline
&&\multicolumn{7}{c}{  vii}&\multicolumn{7}{c}{  viii}\\ \hline
100&30&0.052&0.016&0.014&0.045&0.079&0.04&0.025&0.051&0.01&0.011&0.042&0.081&0.035&0.025\\
100&60&0.038&0.007&0.008&0.041&0.084&0.035&0.013&0.042&0.01&0.008&0.03&0.069&0.026&0.015\\
100&120&0.036&0.008&0.005&0.03&0.078&0.02&0.005&0.038&0.005&0.001&0.031&0.084&0.025&0.005\\
100&240&0.033&0.008&0.006&0.032&0.103&0.022&0&0.033&0.003&0.002&0.028&0.101&0.018&0.005\\
200&30&0.046&0.022&0.021&0.044&0.059&0.044&0.052&0.047&0.016&0.015&0.048&0.058&0.04&0.027\\
200&60&0.043&0.017&0.013&0.041&0.063&0.038&0.033&0.052&0.016&0.015&0.054&0.068&0.046&0.041\\
200&120&0.038&0.015&0.014&0.032&0.061&0.027&0.016&0.033&0.014&0.012&0.042&0.066&0.034&0.016\\
200&240&0.037&0.013&0.008&0.047&0.07&0.041&0.001&0.037&0.01&0.011&0.039&0.066&0.031&0.013\\
\hline
\hline
\end{tabular}}
\end{center}
\end{table}

\begin{figure}[p]
\centering
\includegraphics[width=1\textwidth]{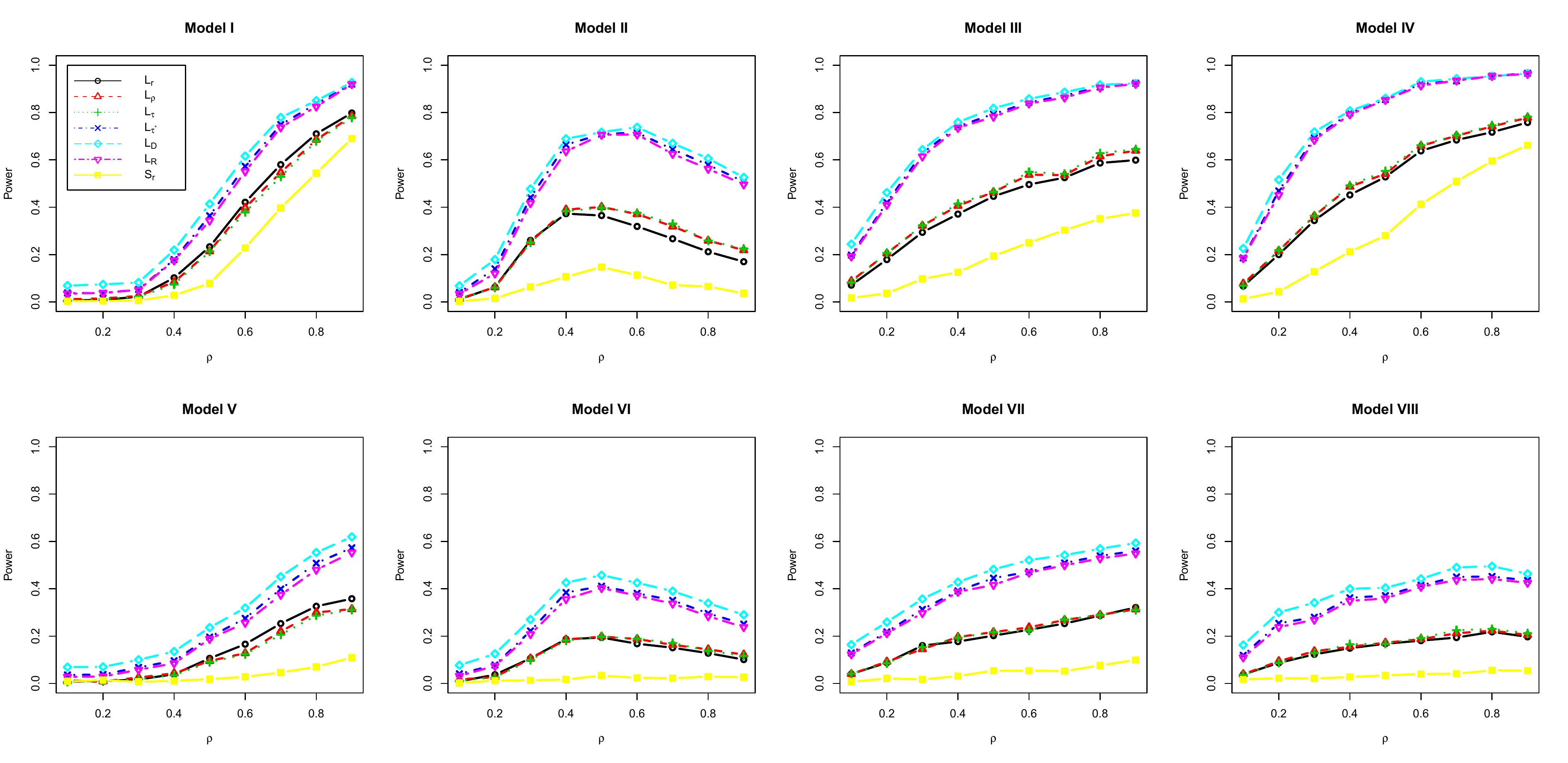}
\caption{\label{fig2}{Power curves of different methods with different $\rho$ and $k_0=2, n=100,p=30, K=4$.}}
\end{figure}

\begin{figure}[p]
\centering
\includegraphics[width=1\textwidth]{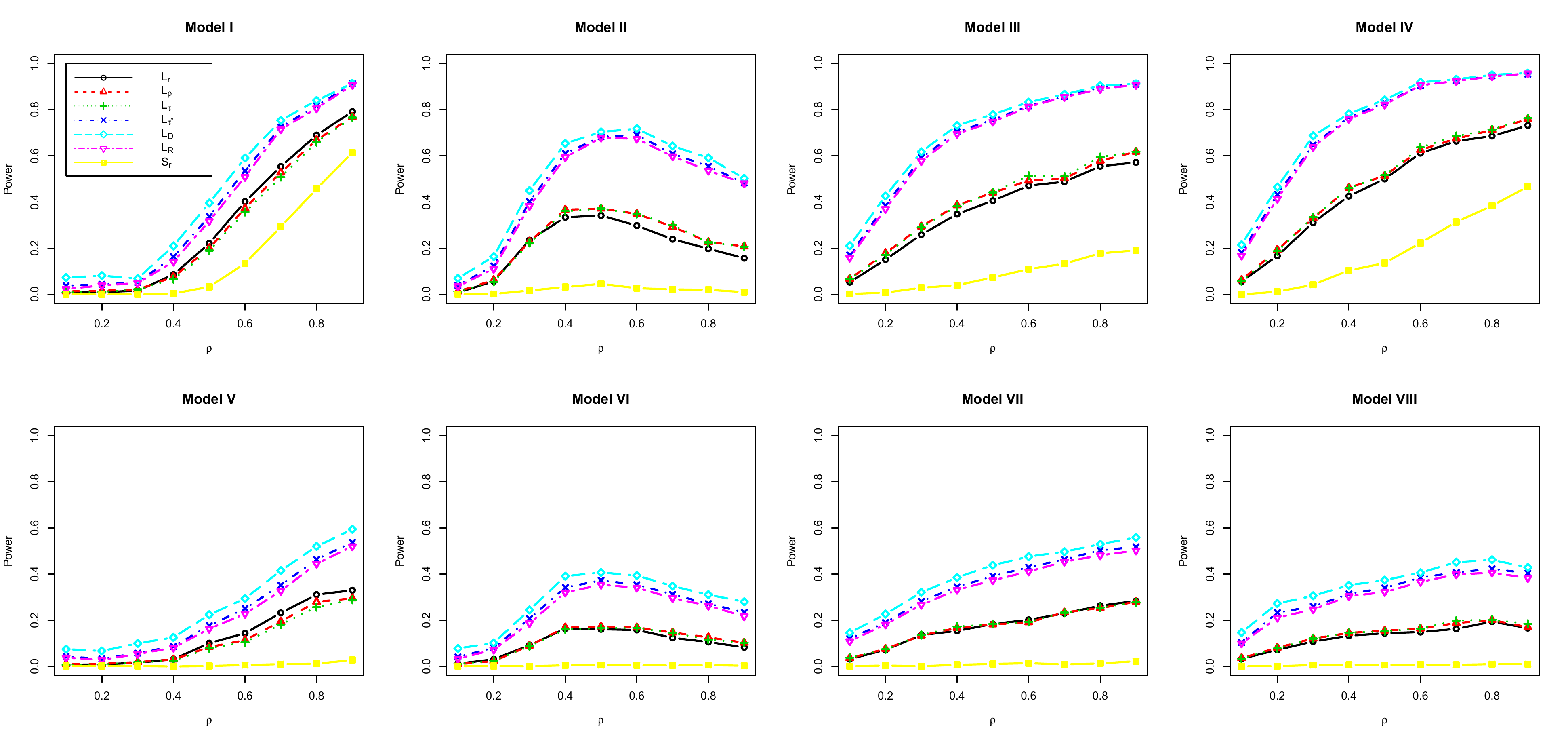}
\caption{\label{fig3}{Power curves of different methods with different $\rho$ and $k_0=2, n=100,p=30, K=6$.}}
\end{figure}

}
\section{Proof of Theorems}

\subsection{Proof of Theorems of Simple Linear Rank Statistics}

\begin{lemma}\label{ri}
Suppose that $X, Y$ are two independent continuous random variables. Let $X_{1}, \ldots, X_{n}$ and $Y_{1}, \ldots, Y_{n}$ be independent observations of $X$ and $Y .$ Let $\left\{Q_{i}^{X}, i=1, \ldots, n\right\}$ and $\left\{Q_{i}^{Y}, i=\right.$ $1, \ldots, n\}$ be the rank of $X_{i}$ and $Y_{i}$ in the samples $\left\{X_{i}\right\}_{i=1}^{n}$ and $\left\{Y_{i}\right\}_{i=1}^{n} .$ Let $\left\{R_{n i}\right\}_{i=1}^{n}$ represent the relative ranks:
$$
R_{n i}=Q_{i^{\prime}}^{Y} \text { subject to } Q_{i^{\prime}}^{X}=i
$$
We then have $\left\{R_{n 1}, \ldots, R_{n n}\right\}$ are uniformly distributed in all permutations of $\{1, \ldots, n\}$ with
$$
\operatorname{pr}\left(R_{n 1}=i_{1}, \ldots, R_{n n}=i_{n}\right)=\frac{1}{n !}
$$
for any permeation $\left\{i_{1}, \ldots, i_{n}\right\}$ of $\{1, \ldots, n\} .$ Here n! represents the factorial of n.
\end{lemma}

\begin{lemma} \label{le62}
{\it(Concentration Inequality For Simple Linear Rank Statistics)}   . Assume the setting and notation in Lemma \ref{ri} . Consider the simple linear rank statistic
$$
V \equiv \sum_{i=1}^{n} c_{n i} g\left(\frac{R_{n i}}{n+1}\right)=\frac{1}{n} \sum_{i=1}^{n} f\left(\frac{Q_{i}^{X}}{n+1}\right) g\left(\frac{Q_{i}^{Y}}{n+1}\right)
$$
where $f(\cdot)$ and $g(\cdot)$ are Lipschitz functions with Lipschitz constant $\Delta<\infty$ and $\max \{|f(0)|,|g(0)|\} \leq$ $A_{2} .$ We have, for any $t>0$
$$
\operatorname{pr}(|V-E V|>t) \leq 2 \exp \left(-C n t^{2}\right)
$$
for some scalar $C$ only depending on $\Delta$ and $A_{2} .$
\end{lemma}

\begin{lemma}\label{st2}
Suppose that the regularity conditions in Theorem \ref{th11} hold. Under the null hypothesis $H_{0}$ holds, we have in the region $x \in\left(0, O\left(n^{1 / 6-\epsilon}\right)\right)$ for some $\epsilon>0$
$$
\operatorname{pr}\left[\frac{V_{ij}(k)-E_{H_{0}}\left(V_{ij}(k)\right)}{\left\{\operatorname{var}_{H_{0}}\left(V_{ij}(k)\right)\right\}^{1 / 2}}>x\right]=\{1-\Phi(x)\}\left\{1+O\left(\frac{1+x^{3}}{n^{1 / 2}}+\frac{x}{n^{1 / 6}}\right)\right\}
$$
\end{lemma}
\subsubsection{Proof of Theorem \ref{th11}}
And without loss of generality, we assume that $\sum_{i=1}^n c_{ni}=0$. Thus, $E_{H_0}(V_{ij}(k))=0$.
Let $\{\nu_i\}_{i=1}^N=\{V_{ij}(k)/\sigma_V\}_{1\le k\le K,1\le i,j\le p}$. Define $z=(2\log(N)-\log\log(N)+y)^{1/2}$.
By Lemma \ref{st2}, we have
\begin{align*}
P(|\nu_i| \geq z)=2 \{1-\Phi(z)\}\{1+o(1)\}\sim \frac{1}{\sqrt{\pi}}\frac{e^{-y/2}}{N}
\end{align*}
Thus,
$$
P\left(\max _{i \in C_{N}}\left|\nu_{i}\right|>z\right) \leq\left|C_{N}\right| \cdot P(|\nu_i| \geq z) \rightarrow 0
$$
$\operatorname{as } p \rightarrow \infty .$ Set $D_{N}:=\left\{1 \leq i \leq N ;\left|B_{N, i}\right|<N^{\varsigma}\right\} .$ By assumption, $\left|D_{N}\right| / N \rightarrow 1$ as $N \rightarrow \infty$
Easily,
$$
\begin{aligned}
P\left(\max _{i \in D_{N}}\left|\nu_{i}\right|>z\right) & \leq P\left(\max _{1 \leq i \leq N}\left|\nu_{i}\right|>z\right) \\
& \leq P\left(\max _{i \in D_{N}}\left|\nu_{i}\right|>z\right)+P\left(\max _{i \in C_{N}}\left|\nu_{i}\right|>z\right)
\end{aligned}
$$
Therefore, to prove Theorem 1, it is enough to show
$$
\lim _{N \rightarrow \infty} P\left(\max _{i \in D_{N}}\left|\nu_{i}\right|>z\right)=1-\exp \left(-\frac{1}{\sqrt{\pi}} e^{-x / 2}\right)
$$
as $N \rightarrow \infty$. Define
$$
\alpha_{t}=\sum^* P\left(\left|\nu_{i_{1}}\right|>z, \cdots,\left|\nu_{i_{t}}\right|>z\right)
$$
for $1 \leq t \leq N,$ where the sum runs over all $i_{1}<\cdots<i_{t}$ and $i_{1} \in D_{N}, \cdots, i_{t} \in D_{N} .$
First, we will prove next that
$$
\lim _{N \rightarrow \infty} \alpha_{t}=\frac{1}{t !} \pi^{-t / 2} e^{-t y / 2}
$$
for each $t \geq 1 .$ Because $g(F_{j}(\varepsilon_{t+k,j}))$ is bounded by a constant $C_g$, thus all the assumptions in Theorem 1.1 in \citet{zaitsev1987gaussian} are satisfied. Thus, we have
\begin{align*}
 &\sum^* P\left(\left|Z_{i_{1}}\right|>z+\epsilon_n(\log(N))^{-1}, \cdots,\left|Z_{i_{t}}\right|>z+\epsilon_n(\log(N))^{-1}\right)\\
 &-\left(\begin{array}{c} |D_N|\\ t\end{array}\right) c_1t^{5/2}\exp\left(-\frac{n^{1/2}\epsilon_n}{c_2t^{3}(\log N)^{1/2}}\right)\\
\le &\sum^* P\left(\left|\nu_{i_{1}}\right|>z, \cdots,\left|\nu_{i_{t}}\right|>z\right)\\
\le &\sum^* P\left(\left|Z_{i_{1}}\right|>z-\epsilon_n(\log(N))^{-1}, \cdots,\left|Z_{i_{t}}\right|>z-\epsilon_n(\log(N))^{-1}\right)\\
&+\left(\begin{array}{c} |D_N|\\ t\end{array}\right)c_1t^{5/2}\exp\left(-\frac{n^{1/2}\epsilon_n}{c_2t^{3}(\log N)^{1/2}}\right)
\end{align*}
where $(Z_{i_1},\cdots,Z_{i_t})$ follows a multivariate normal distribution with mean zero and the same covariance matrix with $(\nu_{i_1},\cdots,\nu_{i_t})$. By the proof of Theorem 2 in \citet{feng2022test}, we have
\begin{align*}
\sum^* P\left(\left|Z_{i_{1}}\right|>z+\epsilon_n(\log(N))^{-1}, \cdots,\left|Z_{i_{t}}\right|>z+\epsilon_n(\log(N))^{-1}\right)\to \frac{1}{t !} \pi^{-t / 2} e^{-t y / 2}\\
\sum^* P\left(\left|Z_{i_{1}}\right|>z-\epsilon_n(\log(N))^{-1}, \cdots,\left|Z_{i_{t}}\right|>z-\epsilon_n(\log(N))^{-1}\right)\to \frac{1}{t !} \pi^{-t / 2} e^{-t y / 2}
\end{align*}
with $\epsilon_n\to 0$ and $N\to \infty$.
Additionally,
\begin{align*}
\left(\begin{array}{c} |D_N|\\ t\end{array}\right)c_1t^{5/2}\exp\left(-\frac{n^{1/2}\epsilon_n}{c_2t^{3}(\log N)^{1/2}}\right)\le C \left(\begin{array}{c} N\\ t\end{array}\right)t^{5/2}\exp\left(-\frac{n^{1/2}\epsilon_n}{c_2t^{3}(\log N)^{1/2}}\right)\to 0
\end{align*}
for $\epsilon_n\to 0$ sufficiently slow. Thus, we have
\begin{align*}
\sum^* P\left(\left|\nu_{i_{1}}\right|>z, \cdots,\left|\nu_{i_{t}}\right|>z\right)\to \frac{1}{t !} \pi^{-t / 2} e^{-t y / 2}.
\end{align*}
Then, by Bonferroni inequality,
$$
\sum_{t=1}^{2 k}(-1)^{t-1} \alpha_{t} \leq P\left(\max _{i \in D_{N}}\left|\nu_{i}\right|>z\right) \leq \sum_{t=1}^{2 k+1}(-1)^{t-1} \alpha_{t}
$$
for any $k \geq 1 .$  let $N \rightarrow \infty$, we have
$$
\begin{aligned}
\sum_{t=1}^{2 k}(-1)^{t-1} \frac{1}{t !}\left(\frac{1}{\sqrt{\pi}} e^{-x / 2}\right)^{t} & \leq \liminf _{N \rightarrow \infty} P\left(\max _{i \in D_{N}}\left|\nu_{i}\right|>z\right) \\
& \leq \limsup _{N \rightarrow \infty} P\left(\max _{i \in D_{N}}\left|\nu_{i}\right|>z\right) \leq \sum_{t=1}^{2 k+1}(-1)^{t-1} \frac{1}{t !}\left(\frac{1}{\sqrt{\pi}} e^{-x / 2}\right)^{t}
\end{aligned}
$$
for each $k \geq 1 .$ By letting $k \rightarrow \infty$ and using the Taylor expansion of the function $1-e^{-x}$, so we obtain the result.

\subsubsection{Proof of Theorem \ref{th12}}
By Lemma \ref{le62}, there exist a constant $c$ such that, for any $t>0$,
\begin{align*}
P\left(|\hat{V}_{ij}(k)-E(\hat{V}_{ij}(k))|>t\right)\le 2e^{-t^2/c}.
\end{align*}
Then,
\begin{align*}
P\left(\max_{1\le i,j\le p,1\le k\le K}|\hat{V}_{ij}(k)-E(\hat{V}_{ij}(k))|>t\right)\le N2e^{-t^2/c}.
\end{align*}
which implies that, with probability at least $1-N^{-1}$,
\begin{align*}
\max_{1\le i,j\le p,1\le k\le K}|\hat{V}_{ij}(k)-E(\hat{V}_{ij}(k))|\le \sqrt{3c\log N}.
\end{align*}
So, for large enough $n$, we have
\begin{align*}
V_{n}^{2} / \sigma_{V}^{2}&=\max_{1\le i,j\le p,1\le k\le K} \widehat{V}_{ij}(k)^{2} \geq \left(\max_{1\le i,j\le p,1\le k\le K}\left|E(\hat V_{ij}(k))\right|-\max_{1\le i,j\le p,1\le k\le K}\left|\widehat{V}_{ij}(k)-E(\hat V_{ij}(k))\right|\right)^{2}\\
 &\geq (2+\varrho)\log N
\end{align*}
for some small positive constant $\varrho$. Accordingly, for any given $q_\alpha$, with probability tending to one,
\begin{align*}
V_{n}^{2} / \sigma_{V}^{2}>2\log N-\log\log N-q_{\alpha}.
\end{align*}
Then we complete the proof. \hfill$\square$

\subsubsection{Proof of Theorem \ref{th13}}

According to Theorem 3 in \citet{feng2022whitenoise} and Assumption (A2), we can easily obtain the result. \hfill$\square$

\subsection{Proof of Theorems of Non-Degenerate U-Statistics}
\begin{lemma}\label{c7}
Suppose that $U$ is a $U$-statistic with degree $m$ and bounded kernel $|h(\cdot)| \leq M$. We then have, for any $t>0$,
$$
P(|U-E U|>t) \leq 2 \exp \left\{-n t^{2} /\left(2 m M^{2}\right)\right\}.
$$
\end{lemma}

\begin{lemma}\label{st22}
Suppose that the boundedness assumption in Theorem \ref{th21} hold. We then have, in a region $x \in\left(0, o\left(n^{1 / 6}\right)\right)$
$$
P\left[\frac{U_{ij}(k)-E\left(U_{ij}(k)\right)}{\left\{\operatorname{var}\left(U_{ij}(k)\right)\right\}^{1 / 2}}>x\right]=\{1-\Phi(x)\}\left\{1+O\left(\frac{1+x^{3}}{n^{1 / 2}}\right)\right\}.
$$
\end{lemma}
\subsubsection{Proof of Theorem \ref{th21}}

First, we consider the following U-statistics with bounded and symmetric kernels, i.e.
\begin{align}
{U}_{ij}(k)=&\frac{1}{C_{n-k}^m}\sum_{1\le t_1<t_2,\cdots,<t_m\le n-k}h((\varepsilon_{t_1,i},\varepsilon_{t_1+k,j})^\top,\cdots,(\varepsilon_{t_m,i},\varepsilon_{t_m+k,j})^\top)
\end{align}
Here we define $\{u_s\}_{s=1}^{N}=\{U_{ij}(k)\}_{1\le i,j\le p, 1\le k\le K}$ and $\{\X_{tijk}\}=\{(\varepsilon_{t,i},\varepsilon_{t+k,j})^\top\}_{1\le t\le n-k}$. So we rewrite $U_{ij}(k)$ in the following forms
\begin{align}\label{uf}
u_s=\frac{1}{C_{n-k_s}^m}\sum_{1\le t_1<t_2,\cdots,<t_m\le n-k_s}h(\X_{t_1,i_s j_s k_s},\cdots,\X_{t_m,i_s j_s k_s})
\end{align}
Without loss of generality, we assume that $E(u_s)=0$. By the condition, we have
\begin{align*}
\mu_q\doteq& E|h(\X_{t_1},\cdots,\X_{t_m})|^q<\infty,\\
\psi_s(x)=&E(h(\X_{t_1},\cdots,\X_{t_m})|\X_{t_1}=x),~~ \sigma_{\psi}^2=\var(\psi_s(\X_{t_1}))>0.
\end{align*}
for any $q\ge 2$.
By Lemma 1 in \citet{malevich1979large}, we can rewrite $u_s$ as follow
\begin{align*}
u_s=S_s+\eta_s, ~~ S_s=\frac{m}{n-k_s}\sum_{i=1}^{n-k_s}\psi_s(\X_i),~~ \eta_s=\sum_{l=2}^m C_m^l u_{s,l}
\end{align*}
where $u_{s,l}$ is a $U$-statistics of the form (\ref{uf}) with kernel $\psi^{(l)}(x_1,\cdots,x_l)$ such that
\begin{align*}
&\mathbf{P}\left\{\mathbf{E}\left[\psi^{(l)}\left(\X_{1}, \cdots, \X_{l}\right) \mid \X_{1}, \cdots, \X_{l-1}\right]=0\right\}=1,\\
&\mathbf{E}\left|\psi^{(l)}\left(\X_{1}, \cdots, \X_{l}\right)\right|^{q} \leqq 2^{m q} \mu_{q}.
\end{align*}
By Lemma 2 in \citet{malevich1979large}, we have
\begin{align*}
E|C_{n-k_s}^lu_{s,l}|^q\le C (n-k_s)^{lq/2}.
\end{align*}
Thus, by the Markov inequality,
\begin{align*}
&P(\max_{1\le s \le N}(n - k_s)^{1/2}|\eta_s|>(\log N)^{-1})\\
&\le NP(|\eta_s|>(\log N)^{-1}(n - k_s)^{-1/2})\\
&\le N(\log N)^{2q}(n-k_s)^{q} E(|\eta_s|^{2q})\\
&\le N(\log N)^{2q}(n-k_s)^{q} E\left(\left|\sum_{l=2}^mC_m^l u_{s,l}\right|^{2q}\right)\\
& =  N(\log N)^{2q}(n-k_s)^{q} m^{2q} E\left(\left|\frac{1}{m}\sum_{l=2}^mC_m^l u_{s,l}\right|^{2q}\right)\\
&\le N(\log N)^{2q}(n-k_s)^{q} m^{2q-1} \left\{ \sum_{l=2}^m (C_m^l)^{2q} E(|u_{s,l}|^{2q}) \right\}\\
&\le C N(\log N)^{2q}(n-k_s)^{q} m^{2q-1} \left\{ \sum_{l=2}^m (C_m^l)^{2q} (n-k_s)^{-lq} \right\} \\
&\le C m^{2q-1} \left\{ \sum_{l=2}^m (C_m^l)^{2q} \right\} N(\log N)^{2q}(n-k_s)^{q}  (n-k_s)^{-2q}\\
& =  C_0(m,q) N(\log N)^{2q} (n-k_s)^{-q} \to 0,
\end{align*}
for some positive integer $q$ by $N=o(n^\epsilon)$.
Thus, by
\begin{align*}\label{u1}
  &\left|\max_{1\le s\le N}(n-k_s)u_s^2-\max_{1\le s\le N}(n-k_s)S_s^2\right|\le 2\max_{1\le s\le N}(n-k_s)^{1/2}|S_s|\max_{1\le s\le N}(n-k_s)^{1/2}|\eta_s|+\max_{1\le s\le N}(n-k_s)\eta_s^2
\end{align*}
we only need to show that
$$\operatorname{P}\left( \max_{1\le s\le N}(n-k_s)S^2_s/ \sigma_{U}^{2}-2\log(Kp^2)+\log \log (Kp^2) \leqslant y\right)\to\exp \left\{-\pi^{-1 / 2} \exp (-y / 2)\right\}
$$

 Here we define $\upsilon_s=(n-k_s)^{1/2}S_s/\sigma_U$
and $z=(2\log(N)-\log\log(N)+y)^{1/2}$.
Since $(n-k_s)^{1/2}\eta_s$ is negligible, the tail behavior of $\upsilon_s$ is the same as that of $(n-k_s)^{1/2}u_s/\sigma_U$. Therefore, by Lemma \ref{st22}, we have
\begin{align*}
P(|\upsilon_i| \geq z)=2\{1-\Phi(z)\}\{1+o(1)\}\sim \frac{1}{\sqrt{\pi}}\frac{e^{-y/2}}{N}
\end{align*}
Thus,
$$
P\left(\max _{i \in C_{N}}\left|\upsilon_{i}\right|>z\right) \leq\left|C_{N}\right| \cdot P(|\upsilon_i| \geq z) \rightarrow 0
$$
$\operatorname{as} p \rightarrow \infty .$ Set $D_{N}:=\left\{1 \leq i \leq N ;\left|B_{N, i}\right|<N^{\varsigma}\right\} .$ By assumption, $\left|D_{N}\right| / N \rightarrow 1$ as $N \rightarrow \infty$
Easily,
$$
\begin{aligned}
P\left(\max _{i \in D_{N}}\left|\upsilon_{i}\right|>z\right) & \leq P\left(\max _{1 \leq i \leq N}\left|\upsilon_{i}\right|>z\right) \\
& \leq P\left(\max _{i \in D_{N}}\left|\upsilon_{i}\right|>z\right)+P\left(\max _{i \in C_{N}}\left|\upsilon_{i}\right|>z\right)
\end{aligned}
$$
Therefore, to prove Theorem \ref{th21}, it is enough to show
$$
\lim _{N \rightarrow \infty} P\left(\max _{i \in D_{N}}\left|\upsilon_{i}\right|>z\right)=1-\exp \left(-\frac{1}{\sqrt{\pi}} e^{-x / 2}\right)
$$
as $N \rightarrow \infty$. Define
$$
\alpha_{t}=\sum^* P\left(\left|\upsilon_{i_{1}}\right|>z, \cdots,\left|\upsilon_{i_{t}}\right|>z\right)
$$
for $1 \leq t \leq N,$ where the sum runs over all $i_{1}<\cdots<i_{t}$ and $i_{1} \in D_{N}, \cdots, i_{t} \in D_{N} .$
First, we will prove next that
$$
\lim _{N \rightarrow \infty} \alpha_{t}=\frac{1}{t !} \pi^{-t / 2} e^{-t y / 2}
$$
for each $t \geq 1 .$ Because $\psi_s(\X_i)$ is bounded, thus all the assumptions in Theorem 1.1 in \citet{zaitsev1987gaussian} are satisfied. Thus, we have
\begin{align*}
 &\sum^* P\left(\left|Z_{i_{1}}\right|>z+\epsilon_n(\log(N))^{-1/2}, \cdots,\left|Z_{i_{t}}\right|>z+\epsilon_n(\log(N))^{-1/2}\right)\\
 &-\left(\begin{array}{c} |D_N|\\ t\end{array}\right) c_1t^{5/2}\exp\left(-\frac{n^{1/2}\epsilon_n}{c_2t^{3}(\log N)^{1/2}}\right)\\
\le &\sum^* P\left(\left|\upsilon_{i_{1}}\right|>z, \cdots,\left|\upsilon_{i_{t}}\right|>z\right)\\
\le &\sum^* P\left(\left|Z_{i_{1}}\right|>z-\epsilon_n(\log(N))^{-1/2}, \cdots,\left|Z_{i_{t}}\right|>z-\epsilon_n(\log(N))^{-1/2}\right)\\
&+\left(\begin{array}{c} |D_N|\\ t\end{array}\right)c_1t^{5/2}\exp\left(-\frac{n^{1/2}\epsilon_n}{c_2t^{3}(\log N)^{1/2}}\right)
\end{align*}
where $(Z_{i_1},\cdots,Z_{i_t})$ follows a multivariate normal distribution with mean zero and the same covariance matrix with $(\upsilon_{i_1},\cdots,\upsilon_{i_t})$. By the proof of Theorem 2 in \citet{feng2022test}, we have
\begin{align*}
\sum^* P\left(\left|Z_{i_{1}}\right|>z+\epsilon_n(\log(N))^{-1/2}, \cdots,\left|Z_{i_{t}}\right|>z+\epsilon_n(\log(N))^{-1/2}\right)\to \frac{1}{t !} \pi^{-t / 2} e^{-t y / 2}\\
\sum^* P\left(\left|Z_{i_{1}}\right|>z-\epsilon_n(\log(N))^{-1/2}, \cdots,\left|Z_{i_{t}}\right|>z-\epsilon_n(\log(N))^{-1/2}\right)\to \frac{1}{t !} \pi^{-t / 2} e^{-t y / 2}
\end{align*}
with $\epsilon_n\to 0$ and $N\to \infty$.
Additionally,
\begin{align*}
\left(\begin{array}{c} |D_N|\\ t\end{array}\right)c_1t^{5/2}\exp\left(-\frac{n^{1/2}\epsilon_n}{c_2t^{3}(\log N)^{1/2}}\right)\le C \left(\begin{array}{c} N\\ t\end{array}\right)t^{5/2}\exp\left(-\frac{n^{1/2}\epsilon_n}{c_2t^{3}(\log N)^{1/2}}\right)\to 0
\end{align*}
for $\epsilon_n\to 0$ sufficiently slow. Thus, we have
\begin{align*}
\sum^* P\left(\left|\upsilon_{i_{1}}\right|>z, \cdots,\left|\upsilon_{i_{t}}\right|>z\right)\to \frac{1}{t !} \pi^{-t / 2} e^{-t y / 2}.
\end{align*}
Then, by Bonferroni inequality,
$$
\sum_{t=1}^{2 k}(-1)^{t-1} \alpha_{t} \leq P\left(\max _{i \in D_{N}}\left|\upsilon_{i}\right|>z\right) \leq \sum_{t=1}^{2 k+1}(-1)^{t-1} \alpha_{t}
$$
for any $k \geq 1 .$  let $N \rightarrow \infty$, we have
$$
\begin{aligned}
\sum_{t=1}^{2 k}(-1)^{t-1} \frac{1}{t !}\left(\frac{1}{\sqrt{\pi}} e^{-x / 2}\right)^{t} & \leq \liminf _{N \rightarrow \infty} P\left(\max _{i \in D_{N}}\left|\upsilon_{i}\right|>z\right) \\
& \leq \limsup _{N \rightarrow \infty} P\left(\max _{i \in D_{N}}\left|\upsilon_{i}\right|>z\right) \leq \sum_{t=1}^{2 k+1}(-1)^{t-1} \frac{1}{t !}\left(\frac{1}{\sqrt{\pi}} e^{-x / 2}\right)^{t}
\end{aligned}
$$
for each $k \geq 1 .$ By letting $k \rightarrow \infty$ and using the Taylor expansion of the function $1-e^{-x}$, so we obtain the result.\hfill$\square$

\subsubsection{Proof of Theorem \ref{th22}}
By Lemma \ref{c7}, we have
\begin{align*}
P\left(|\hat{U}_{ij}(k)-E(\hat{U}_{ij}(k))|>t\right)\le 2e^{t^2/c}
\end{align*}
for some positive constant $c$. Taking the same procedure as Theorem \ref{th12}, we can also obtain the result. \hfill$\square$

\subsubsection{Proof of Theorem \ref{th23}}

According to Theorem 3 in \citet{feng2022whitenoise} and Assumption (A3), we can easily obtain the result. \hfill$\square$

\subsection{Proof of Theorems of Degenerate U-Statistics}
\subsubsection{Proof of Theorem \ref{th31}}

{

\begin{lemma}\lbl{winter_fire} For $N\geq 1$, let $\iota_N$ be positive integers with $\lim_{N\to\infty}\iota_N/N=1$. Let $Y_{iv}, i=1,\cdots,\iota_N, v=1,\cdots,M$ be $N(0,1)$-distributed random variables and $\cov(Y_{iv},Y_{is})=0$ for $v \not=s$. Let $\Y_i=(Y_{i1},\cdots,Y_{iM})^\top$ and $\Xi_{ij}=\cov(\Y_i,\Y_j)$. Assume $|\lambda_{max}(\Xi_{ij}\Xi_{ij}^\top)|\leq \delta_N^{2+2c}$ for $c>0$ and all $1\leq i<j \leq \iota_N$ , where  $\{\delta_N;\, N\geq 1\}$ are constants satisfying $0<\delta_N=o(1/\log N).$
Define $W_{i_k}=\sum_{v=1}^M\lambda_vY^2_{i_kv}$.
Given $x\in \mathbb{R}$, set $z=2  \lambda_1\log (N)+ \lambda_1\left(\mu_{1}-2\right) \log \log (N)+ \lambda_{1} y+o(1/\log(N)).$ Then, for any fixed $m\geq 1$, we have
\bea\lbl{cskh}
\Big(\frac{\Gamma\left(\mu_{1} / 2\right)}{ \kappa}\frac{N}{e^{-y/2}}\Big)^m\cdot P(W_{i_1}>z, \cdots, W_{i_m}>z)
\to 1
\eea
as $N\to\infty$ uniformly for all $1\leq i_1<\cdots < i_m\leq \iota_N.$
\end{lemma}
\noindent\textbf{Proof of Lemma \ref{winter_fire}} For $m=1$, (\ref{cskh}) is followed by Equation (6) in \citet{zolotarev1961concerning}. Assume Equation (\ref{cskh}) holds with $m=k-1$. We will prove it also holds with $m=k$.

Define $\Y_S=(\Y_{i_1}^\top,\cdots,\Y_{i_{m-1}}^\top)^\top$, $\bms_{i_mS}=\cov(\Y_{i_m},\Y_S)$ and $\bms_{Si_m}=\bms_{i_mS}^\top$. So $\Y_{i_m}=(\Y_{i_m}-\bms_{i_mS}\Y_S)+\bms_{i_mS}\Y_S\doteq U_Y+V_Y$ . Thus, by the conditional distribution of multivariate normal distributions, we have $\Y_{i_m}-\bms_{i_mS}\Y_S\sim N(\bm 0, \I_M-\bms_{i_mS}\bms_{Si_m})$ is independent of $\Y_S$.
Define $\A=\diag\{\lambda_1,\cdots,\lambda_M\}$. Thus, we have
\begin{align*}
&P(W_{i_1}>z, \cdots, W_{i_m}>z)\\
=&P\left(U_Y^\top \A U_Y+2U_Y^\top \A V_Y+V_Y^\top\A V_Y\ge z, \min_{1\le l\le m-1}\Y_{i_l}^\top \A\Y_{i_l}\ge z\right)\\
= &P\left(U_Y^\top \A U_Y+2U_Y^\top \A V_Y+V_Y^\top\A V_Y\ge z, 2U_Y^\top \A V_Y+V_Y^\top\A V_Y\le C\delta_N, \min_{1\le l\le m-1}\Y_{i_l}^\top \A\Y_{i_l}\ge z\right)\\
&+P\left(U_Y^\top \A U_Y+2U_Y^\top \A V_Y+V_Y^\top\A V_Y\ge z, 2U_Y^\top \A V_Y+V_Y^\top\A V_Y>C\delta_N, \min_{1\le l\le m-1}\Y_{i_l}^\top \A\Y_{i_l}\ge z\right)\\
\le & P\left(U_Y^\top \A U_Y\ge z-C\delta_N, \min_{1\le l\le m-1}\Y_{i_l}^\top \A\Y_{i_l}\ge z\right)+P\left(2U_Y^\top \A V_Y+V_Y^\top\A V_Y>C\delta_N\right)\\
\le & P\left(U_Y^\top \A U_Y\ge z-C\delta_N\right)P\left(\min_{1\le l\le m-1}\Y_{i_l}^\top \A\Y_{i_l}\ge z\right)+P\left(U_Y^\top \A V_Y>C\delta_N/4\right)+P\left(V_Y^\top\A V_Y>C\delta_N/2\right)
\end{align*}
By $U_Y\sim N(\bm 0, \I_M-\bms_{i_mS}\bms_{Si_m})$, we have $$U_Y^\top \A U_Y\sim \xi_U^\top (\I_M-\bms_{i_mS}\bms_{Si_m})^{1/2}\A(\I_M-\bms_{i_mS}\bms_{Si_m})^{1/2} \xi_U$$
where $\xi_U\sim N(\bm 0, \I_M)$. Define the eigenvalues of $(\I_M-\bms_{i_mS}\bms_{Si_m})^{1/2}\A(\I_M-\bms_{i_mS}\bms_{Si_m})^{1/2} $ are $\tilde \lambda_1\ge \tilde \lambda_2\ge\cdots\ge \tilde \lambda_{M}$ and $\tilde\Lambda, \tilde{\kappa}, \tilde \mu_1$ are the corresponding parameters as in Proposition 3.2 in \citet{drton2020high}. Because $\lambda_{max}(\bms_{i_mS}\bms_{Si_m})\le \delta_N^{2+2c}$, $\tilde \lambda_1=\lambda_1(1+o(\delta_N))$. So does $\tilde\Lambda, \tilde{\kappa}, \tilde \mu_1$.
So by Equation (6) in \citet{zolotarev1961concerning}, we have
\begin{align*}
P\left(U_Y^\top \A U_Y\ge z-C\delta_N\right)&\to\frac{\tilde{\kappa}}{\Gamma(\tilde{\mu}_1/2)}\left(\frac{z-C\delta_N+\tilde\Lambda}{2\tilde \lambda_1}\right)^{\tilde{\mu}_1/2-1}\exp\left(-\frac{z-C\delta_N+\tilde{\Lambda}}{2\tilde \lambda_1}\right) \\
&\to \frac{ \kappa}{\Gamma\left(\mu_{1} / 2\right)}\frac{e^{-y/2}}{N}
\end{align*}
by the assumption $\delta_N=o(1/\log N)$. So by the
\begin{align*}
P\left(U_Y^\top \A U_Y\ge z-C\delta_N\right)P\left(\min_{1\le l\le m-1}\Y_{i_l}^\top \A\Y_{i_l}\ge z\right)\to  \left(\frac{ \kappa}{\Gamma\left(\mu_{1} / 2\right)}\frac{e^{-y/2}}{N}\right)^m
\end{align*}
Next, we will show that $P\left(U_Y^\top \A V_Y>C\delta_N/4\right)+P\left(V_Y^\top\A V_Y>C\delta_N/2\right)=o(N^{-m})$.
Similarly,
\begin{align*}
V_Y^\top \A V_Y\sim \xi_V^\top \bms_{Si_m }\A \bms_{i_m S}\xi_V
\end{align*}
where $\xi_V\sim N(\bm 0, \I_{(m-1)M)}$. Define the eigenvalues of $\bms_{Si_m }\A \bms_{i_m S}$ are $\zeta_1,\cdots,\zeta_M$. So
\begin{align*}
V_Y^\top \A V_Y\sim \sum_{k=1}^{(m-1)M} \zeta_k \xi_k^2
\end{align*}
where $\xi_k$ are all independently distributed as $N(0,1)$. Thus, for small enough constant $\varpi$,
\begin{align*}
E(\exp(\varpi \delta_N^{-1-c} V_Y^\top \A V_Y))=&\prod_{k=1}^M Ee^{\varpi \delta_N^{-1-c} \zeta_k\xi_k^2}=\exp\left\{-\frac{1}{2}\sum_{k=1}^{(m-1)M}\log[1-2\varpi \delta_N^{-1-c} \zeta_k]\right\}\\
\le & \exp\left(2\varpi \delta_N^{-1-c} \sum_{k=1}^{(m-1)M}\zeta_k\right)
\end{align*}
In addition,
\begin{align*}
\sum_{k=1}^{(m-1)M}\zeta_k=\tr(\bms_{Si_m }\A \bms_{i_m S})\le \lambda_{max}(\bms_{Si_m }\bms_{i_m S})\tr(\A)\le \delta_N^{2+2c} \Lambda.
\end{align*}
So $E(\exp(\varpi \delta_N^{-1-c} V_Y^\top \A V_Y))\le \exp\left(2\varpi \delta_N^{1+c} \Lambda\right) \le C_2$ for some constant $C_2>0$. By the Markov inequality, we have
\begin{align*}
P\left(V_Y^\top\A V_Y>C\delta_N/2\right) &= P\left( \varpi \delta_N^{-1-c} V_Y^\top\A V_Y> C \varpi \delta_N^{-c} /2\right) \\
&\le \exp\left(-C\varpi\delta_N^{-c}/2\right)E(\exp(\varpi \delta_N^{-1-c} V_Y^\top \A V_Y))\\
&\le C_2 \exp\left(-C\varpi\delta_N^{-c}/2\right)=o(N^{-m})
\end{align*}
for large enough constant $C$. Similarly,
\begin{align*}
U_Y^\top \A V_Y\sim \sum_{k=1}^{M} \rho_k \xi_k\eta_k
\end{align*}
where $\eta_k$ are all independently distributed as $N(0,1)$, $\rho_1,\cdots,\rho_M$ are the singular values of $(\I_M-\bms_{i_mS}\bms_{Si_m})^{1/2}\A\bms_{i_mS}$. And then
\begin{align*}
 & P\left(U_Y^\top \A V_Y>C\delta_N/4\right) \\
 =  & P\left(\varpi\delta_N^{-1-c} U_Y^\top \A V_Y>C\varpi\delta_N^{-c}/4\right) \\
\le & \exp\left(-C\varpi\delta_N^{-c}/4\right)E(\exp(\varpi \delta_N^{-1-c} U_Y^\top \A V_Y))\\
\le & \exp\left(-C\varpi\delta_N^{-c}/4\right) \exp\left(-\frac{1}{2}\sum_{k=1}^M \log[1-\varpi^2\delta_N^{-2-2c}\rho_k^2]\right)\\
\le &\exp\left(-C\varpi\delta_N^{-c}/4\right) \exp\left(\varpi^2\delta_N^{-2-2c}\sum_{k=1}^M \rho_k^2\right)\\
\le &\exp\left(-C\varpi\delta_N^{-c}/4\right) \exp\left(\varpi^2\delta_N^{-2-2c}\tr(((\I_M-\bms_{i_mS}\bms_{Si_m})^{1/2}\A\bms_{i_mS})^2)\right)\\
\le &\exp\left(-C\varpi\delta_N^{-c}/4\right) \exp\left(\varpi^2\delta_N^{-2-2c}\lambda_{max}(\bms_{Si_m }\bms_{i_m S})\tr(\A^2)\right)\\
\le &\exp\left(-C\varpi\delta_N^{-c}/4\right) \exp\left(\varpi^2\tr(\A^2)\right)=o(N^{-m})
\end{align*}
for large enough constant $C>0$. Thus, we have
\begin{equation}
P(W_{i_1}>z, \cdots, W_{i_m}>z)\le \left(\frac{ \kappa}{\Gamma\left(\mu_{1} / 2\right)}\frac{e^{-y/2}}{N}\right)^m+o(N^{-m}). \label{w1}
\end{equation}
Further more,
\begin{align*}
&P(W_{i_1}>z, \cdots, W_{i_m}>z)\\
=&P\left(U_Y^\top \A U_Y+2U_Y^\top \A V_Y+V_Y^\top\A V_Y\ge z, \min_{1\le l\le m-1}\Y_{i_l}^\top \A\Y_{i_l}\ge z\right)\\
\ge & P\left(U_Y^\top \A U_Y\ge z+C\delta_N, \min_{1\le l\le m-1}\Y_{i_l}^\top \A\Y_{i_l}\ge z\right)-P\left(2U_Y^\top \A V_Y+V_Y^\top\A V_Y\le -C\delta_N\right)\\
\ge & P\left(U_Y^\top \A U_Y\ge z+C\delta_N\right)P\left(\min_{1\le l\le m-1}\Y_{i_l}^\top \A\Y_{i_l}\ge z\right)-P\left(2U_Y^\top \A V_Y+V_Y^\top\A V_Y\le -C\delta_N\right),
\end{align*}
where the first inequality is based on the fact that:
\begin{align*}
&P\left(\min_{1\le l\le m-1}\Y_{i_l}^\top \A\Y_{i_l}\ge z\right)-P\left(U_Y^\top \A U_Y+2U_Y^\top \A V_Y+V_Y^\top\A V_Y\ge z, \min_{1\le l\le m-1}\Y_{i_l}^\top \A\Y_{i_l}\ge z\right) \\
 = & P\left(U_Y^\top \A U_Y+2U_Y^\top \A V_Y+V_Y^\top\A V_Y < z, \min_{1\le l\le m-1}\Y_{i_l}^\top \A\Y_{i_l}\ge z\right)\\
 = & P\left(U_Y^\top \A U_Y+2U_Y^\top \A V_Y+V_Y^\top\A V_Y < z, 2U_Y^\top \A V_Y+V_Y^\top\A V_Y < - C\delta_N, \min_{1\le l\le m-1}\Y_{i_l}^\top \A\Y_{i_l}\ge z\right) \\
& + P\left(U_Y^\top \A U_Y+2U_Y^\top \A V_Y+V_Y^\top\A V_Y < z, 2U_Y^\top \A V_Y+V_Y^\top\A V_Y \ge - C\delta_N, \min_{1\le l\le m-1}\Y_{i_l}^\top \A\Y_{i_l}\ge z\right) \\
\le & P\left(2U_Y^\top \A V_Y+V_Y^\top\A V_Y < - C\delta_N \right) + P\left(U_Y^\top \A U_Y  < z + C\delta_N,\min_{1\le l\le m-1}\Y_{i_l}^\top \A\Y_{i_l}\ge z\right).
\end{align*}
Obviously,
\begin{align*}
P\left(U_Y^\top \A U_Y\ge z+C\delta_N\right)P\left(\min_{1\le l\le m-1}\Y_{i_l}^\top \A\Y_{i_l}\ge z\right)\to  \left(\frac{ \kappa}{\Gamma\left(\mu_{1} / 2\right)}\frac{e^{-y/2}}{N}\right)^m
\end{align*}
by $\delta_N=o(1/\log N)$. Next,
\begin{align*}
&P\left(2U_Y^\top \A V_Y+V_Y^\top\A V_Y\le -C\delta_N\right)\\
\le & P\left(U_Y^\top \A V_Y\le -C\delta_N/2\right)=P\left(U_Y^\top \A V_Y\ge C\delta_N/2\right)=o(N^{-m}).
\end{align*}
So
\begin{align}\label{w2}
&P(W_{i_1}>z, \cdots, W_{i_m}>z)\ge \left(\frac{ \kappa}{\Gamma\left(\mu_{1} / 2\right)}\frac{e^{-y/2}}{N}\right)^m+o(N^{-m}).
\end{align}
Then, we obtain the result by (\ref{w1}) and (\ref{w2}).   \hfill$\square$

\begin{lemma}\label{winter_cool}
Let $W_{i}=\sum_{v=1}^M\lambda_vY^2_{iv}$ where $Y_{iv}\stackrel{i.i.d}{\sim} N(0,1)$ for $v=1,\cdots,M$. Define $\Y_{i}=(Y_{i1},\cdots,Y_{i M})^\top$ and $\Y_S=(\Y_{i_1}^\top,\cdots,\Y_{i_{m-1}}^\top)^\top$ where $S=\{i_1,\cdots,i_{m-1}\}$ which for any $i,j\in S$, $\lambda_{max}(\Xi_{ij}\Xi_{ij}^\top)\le \delta_N^{2+2c_0}$ for some constant $c_0>0$. Let $\Xi_{ij}=\cov(\Y_{i},\Y_{j})$. Let $\Y_{i_{m}}$ satisfy $\max_{1\le j\le m-1}|\lambda_{max}(\Xi_{i_mj}\Xi_{i_mj}^\top)|>\delta_N^{2+2c_0}$ and $\lambda_{max}(\bms_{i_mS}\bms_{Si_m})\le \delta \in  (0,1)$ where $\bms_{i_mS}=\cov(\Y_{i_m},\Y_S)$ and $\bms_{Si_m}=\bms_{i_mS}^\top$. Then, we have
\begin{align}
P\left(W_{i_m}\ge z, \min_{1\le l\le m-1}W_{i_l}\ge z\right)\le C(\log N)^cN^{\varrho-m}
\end{align}
where $\varrho\in (0,1)$, $c, C$ are some positive constants.
\end{lemma}
\proof Define $\Y_{i_m}=(\Y_{i_m}-\bms_{i_mS}\Y_S)+\bms_{i_mS}\Y_S\doteq U_Y+V_Y$ . Thus, by the conditional distribution of multivariate normal distributions, we have $\Y_{i_m}-\bms_{i_mS}\Y_S\sim N(\bm 0, \I_M-\bms_{i_mS}\bms_{Si_m})$ is independent of $\Y_S$. Define $\A=\diag\{\lambda_1,\cdots,\lambda_M\}$. Thus, we have
\begin{align*}
&P\left(W_{i_m}\ge z, \min_{1\le l\le m-1}W_{i_l}\ge z\right)\\
=&P\left(\Y_{i_m}^\top \A\Y_{i_m}\ge z, \min_{1\le l\le m-1}\Y_{i_l}^\top \A\Y_{i_l}\ge z\right)\\
=&P\left(U_Y^\top \A U_Y+2U_Y^\top \A V_Y+V_Y^\top\A V_Y\ge z, \min_{1\le l\le m-1}\Y_{i_l}^\top \A\Y_{i_l}\ge z\right)\\
\le &P\left(2U_Y^\top \A U_Y+2V_Y^\top\A V_Y\ge z, \min_{1\le l\le m-1}\Y_{i_l}^\top \A\Y_{i_l}\ge z\right)\\
\le &P\left(U_Y^\top \A U_Y\ge \frac{1}{4}z, \min_{1\le l\le m-1}\Y_{i_l}^\top \A\Y_{i_l}\ge z\right)+P\left(V_Y^\top \A V_Y\ge \frac{1}{4}z, \min_{1\le l\le m-1}\Y_{i_l}^\top \A\Y_{i_l}\ge z\right)\\
=&P\left(U_Y^\top \A U_Y\ge \frac{1}{4}z\right)P\left( \min_{1\le l\le m-1}\Y_{i_l}^\top \A\Y_{i_l}\ge z\right)+P\left(V_Y^\top \A V_Y\ge \frac{1}{4}z, \min_{1\le l\le m-1}\Y_{i_l}^\top \A\Y_{i_l}\ge z\right)
\end{align*}
By Lemma \ref{winter_fire}, we have
\begin{align*}
&P\left(\min_{1\le l\le m-1}\Y_{i_l}^\top \A\Y_{i_l}\ge z\right)\leq (1+2\epsilon)^{m-1} \Big(\frac{ \kappa}{\Gamma\left(\mu_{1} / 2\right)}\frac{e^{-y/2}}{N}\Big)^{m-1}\le C N^{1-m}
\end{align*}
By $U_Y\sim N(\bm 0, \I_M-\bms_{i_mS}\bms_{Si_m})$, we have $$U_Y^\top \A U_Y\sim \xi_U^\top (\I_M-\bms_{i_mS}\bms_{Si_m})^{1/2}\A(\I_M-\bms_{i_mS}\bms_{Si_m})^{1/2} \xi_U$$
where $\xi_U\sim N(\bm 0, \I_M)$.
So
\begin{align*}
\tilde\lambda_1&\doteq\lambda_{max}(((\I_M-\bms_{i_mS}\bms_{Si_m})^{1/2}\A(\I_M-\bms_{i_mS}\bms_{Si_m})^{1/2})\\
&\ge \lambda_{max}(\A)\lambda_{\min}(\I_M-\bms_{i_mS}\bms_{Si_m})\ge \lambda_1(1-\delta)
\end{align*}
by the assumption.
By Equation (6) in \citet{zolotarev1961concerning}, we have
\begin{align*}
P\left(U_Y^\top \A U_Y\ge \frac{1}{4}z\right)\le& (1+\epsilon) \frac{\tilde{\kappa}}{\Gamma(\tilde{\mu}_1/2)}\left(\frac{z/4+\tilde\Lambda}{2\tilde \lambda_1}\right)^{\tilde{\mu}_1/2-1}\exp\left(-\frac{z/4+\tilde{\Lambda}}{2\tilde \lambda_1}\right)\\
\le & C\left(\frac{z/4+\tilde\Lambda}{2\lambda_1(1-\delta)}\right)^{\tilde{\mu}_1/2-1}\exp\left(-\frac{z/4+\tilde{\Lambda}}{2 \lambda_1(1-\delta)}\right)\\
\le & C (\log N)^c N^{-(1-\delta)/4}
\end{align*}
where $c=(\tilde{\mu}_1-\mu_1)/2$.
Thus, we have
\begin{align}
P\left(U_Y^\top \A U_Y\ge \frac{1}{4}z\right)P\left( \min_{1\le l\le m-1}\Y_{i_l}^\top \A\Y_{i_l}\ge z\right)\le C (\log N)^c N^{-m+(3+\delta)/4}
\end{align}
Define $\tilde\A=\diag\{\A,\cdots,\A\}\in \mathbb{R}^{(m-1)M\times(m-1)M}$. Next, we consider
\begin{align*}
&P\left(V_Y^\top \A V_Y\ge \frac{1}{4}z, \min_{1\le l\le m-1}\Y_{i_l}^\top \A\Y_{i_l}\ge z\right)\\
=&P\left(\Y_S^\top\bms_{Si_m} \A \bms_{i_mS}\Y_S\ge \frac{1}{4}z, \min_{1\le l\le m-1}\Y_{i_l}^\top \A\Y_{i_l}\ge z\right)\\
\le &P\left(\Y_S^\top\bms_{Si_m} \A \bms_{i_mS}\Y_S\ge \frac{1}{4}z, \Y_{S}^\top \tilde\A\Y_{S}\ge (m-1)z\right)\\
\le & P\left(\Y_S^\top((1-\epsilon) \tilde\A+\epsilon\bms_{Si_m} \A \bms_{i_mS})\Y_S\ge ((m-1)(1-\epsilon)+\frac{1}{4}\epsilon)z\right)\\
=& P\left(\xi_S^\top\bms_S^{1/2}((1-\epsilon) \tilde\A+\epsilon\bms_{Si_m} \A \bms_{i_mS})\bms_S^{1/2}\xi_S\ge ((1-\epsilon)m-1+\frac{5}{4}\epsilon)z\right)
\end{align*}
where $\xi_S \sim N(\bm 0, \I_{(m-1)M})$.
We have
\begin{align*}
\breve\lambda_1\doteq&\lambda_{max}\left(\bms_S^{1/2}((1-\epsilon)\tilde\A+\epsilon\bms_{Si_m} \A \bms_{i_mS})\bms_S^{1/2}\right)\\
\ge & (1-(m-1)\delta_N)\lambda_{max}\left((1-\epsilon)\tilde\A+\epsilon \bms_{Si_m} \A \bms_{i_mS}\right)\\
\ge &(1-(m-1)\delta_N)\left(\lambda_{max}\left((1-\epsilon)\tilde\A\right)-\lambda_{max}\left(\epsilon\bms_{Si_m} \A \bms_{i_mS}\right)\right)\\
\ge &\lambda_1(1-(m-1)\delta_N)(1-(1+\delta)\epsilon)\ge (1-\delta/2)\lambda_1
\end{align*}
for a small enough positive real number $\epsilon$.

By Equation (6) in \citet{zolotarev1961concerning}, we have
\begin{align*}
&P\left(\xi_S^\top\bms_S^{1/2}((1-\epsilon)\tilde\A+\epsilon\bms_{Si_m} \A \bms_{i_mS})\bms_S^{1/2}\xi_S\ge ((1-\epsilon)m-1+\frac{5}{4}\epsilon)z\right)\\
\le & (1+\epsilon) \frac{{\breve{\kappa}}}{\Gamma(\breve{\mu}_1/2)}\left(\frac{((1-\epsilon)m-1+\frac{5}{4}\epsilon)z+\breve\Lambda}{2\breve \lambda_1}\right)^{\breve{\mu}_1/2-1}\exp\left(-\frac{((1-\epsilon)m-1+\frac{5}{4}\epsilon)z+\breve{\Lambda}}{2\breve \lambda_1}\right)\\
\le & C\left(\frac{((1-\epsilon)m-1+\frac{5}{4}\epsilon)z+\breve\Lambda}{2\lambda_1(1-\delta/2)}\right)^{\breve{\mu}_1/2-1}\exp\left(-\frac{((1-\epsilon)m-1+\frac{5}{4}\epsilon)z+\breve{\Lambda}}{2 \lambda_1(1-\delta/2)}\right)\\
\le & C (\log N)^c N^{-\frac{(1-\epsilon)m-1+\frac{5}{4}\epsilon}{1-\delta/2}} \le C (\log N)^c N^{-m+(3+\delta)/4}
\end{align*}
by setting $\frac{(1-\epsilon)m-1+\frac{5}{4}\epsilon}{1-\delta/2}\ge m-(3+\delta)/4$. \hfill$\square$

}

{
\begin{lemma} \label{le67}
Define $\upsilon_t=\sum_{v=1}^{M} \lambda_{v}\left({(n-k_t)}^{-1 / 2} \sum_{i=1}^{{n-k_t}} \phi_v(Z_{i,t})\right)^{2}$ where $Z_{i,t}$ is the corresponding random variable in Condition (C6) with respect to $t\in \{1,\cdots,N\}$.
Let
$$
\beta_{s}=\sum^* P\left(\upsilon_{i_{1}}>z, \cdots,\upsilon_{i_{s}}>z\right)
$$
for $1 \leq s \leq N,$ where the sum runs over all $i_{1}<\cdots<i_{s}$ and $i_{1} \in D_{N}, \cdots, i_{s} \in D_{N} .$
Then,
$$
\lim _{N \rightarrow \infty} \beta_{s}=\frac{1}{s !}\left(\frac{  \kappa}{\Gamma\left(\mu_{1} / 2\right)} e^{-\frac{y}{2}}\right)^{-s}
$$
for each $s \geq 1 .$
\end{lemma}
\proof According to Theorem 1.1 in \citet{zaitsev1987gaussian}, we have
\begin{align*}
 &\sum^* P\left(\sum_{v=1}^M\lambda_vY_{i_1v}^2>z+\epsilon_n(\log(N))^{-1}, \cdots,\sum_{v=1}^M\lambda_vY_{i_sv}^2>z+\epsilon_n(\log(N))^{-1}\right)\\
 &-\left(\begin{array}{c} |D_N|\\ s\end{array}\right) c_1s^{5/2}\exp\left(-\frac{n^{1/2}\epsilon_n}{c_2s^{3}(\log N)}\right)\\
\le &\sum^* P\left(\upsilon_{i_1}>z, \cdots,\upsilon_{i_{s}}>z\right)\\
\le &\sum^* P\left(\sum_{v=1}^M\lambda_vY_{i_1v}^2>z-\epsilon_n(\log(N))^{-1}, \cdots,\sum_{v=1}^M\lambda_vY_{i_sv}^2>z-\epsilon_n(\log(N))^{-1}\right)\\
&+\left(\begin{array}{c} |D_N|\\ s\end{array}\right)c_1s^{5/2}\exp\left(-\frac{n^{1/2}\epsilon_n}{c_2s^{3}(\log N)}\right)
\end{align*}
where $(Y_{i_1 v},\cdots,Y_{i_s v})$ follows a multivariate normal distribution with mean zero and the same covariance matrix with $({(n-k_{i_1})}^{-1 / 2} \sum_{j=1}^{{n-k_{i_1}}} \phi_v(Z_{j,{i_1}}),\cdots,{(n-k_{i_s})}^{-1 / 2} \sum_{j=1}^{{n-k_{i_s}}} \phi_v(Z_{j,{i_s}}))$. By the condition $N=o(n^{\epsilon})$, there exist small enough $\epsilon_n\to 0$ satisfy
$$\left(\begin{array}{c} |D_N|\\ s\end{array}\right)c_1s^{5/2}\exp\left(-\frac{n^{1/2}\epsilon_n}{c_2s^{3}(\log N)}\right)\to 0.$$
Define $W_{i}=\sum_{v=1}^M\lambda_vY_{iv}^2$ and $\Y_i=(Y_{i1},\cdots,Y_{iM})$. So we only need to show that
$$\sum^* P\left(W_{i_1}>z, \cdots,W_{i_{s}}>z\right)\to\frac{1}{s !}\left(\frac{  \kappa}{\Gamma\left(\mu_{1} / 2\right)} e^{-\frac{y}{2}}\right)^{-s}.$$
Recalling $D_N:=\{1\leq i \leq N;\, |B_{N,i}|< N^{ \kappa}\}$, we write
\beaa
\big\{(i_1, \cdots, i_t)\in (D_N)^t;\, i_1<\cdots < i_t\big\}=F_t\cup G_t
\eeaa
where $\sigma_{i_ri_s}=\lambda_{max}(\Xi_{i_ri_s}\Xi_{i_ri_s}^\top)$ and $\Xi_{i_ri_s}=\cov(\Y_{i_r},\Y_{i_s})$.
\bea
&&F_t:=\big\{(i_1, \cdots, i_t)\in (D_N)^t;\, i_1<\cdots < i_t\ \mbox{and}\ |\sigma_{i_r i_s}|\leq \delta_N^{2+2c}\ \mbox{for all}\ 1\leq r< s\leq t\}; \nonumber\\
&& G_t:=\big\{(i_1, \cdots, i_t)\in (D_N)^t;\, i_1<\cdots < i_t\ \mbox{and}\ |\sigma_{i_r i_s}|> \delta_N^{2+2c}\ \mbox{for a pair}\ (i_r, i_s)\ \mbox{with} \nonumber\\
 && ~~~~~~~~~~~~~~~~~~~~~~~~~~~~~~~~~~~~~~~~~~~~~~~~~~~~~~~~~~~~~~~~~~~~~~~~~~~~~~~~~~~~~~~1\leq r< s\leq t\big\}. \nonumber\\
 && \lbl{urine_problem}
\eea
Now, think $D_N$ as graph with $|D_N|$ vertices. Keep in mind that $|D_N|\leq N$ and $|D_N|/ N\to 1$.  Any two different vertices from  them, say, $i$ and $j$  are connected if $|\sigma_{i j}|> \delta_N^{2+2c}.$ In this case we also say there is an edge between them. By the definition $D_N$, each vertex in the graph has at most $N^{\varsigma}$ neighbors. Replacing ``$n$", ``$q$" and ``$t$" in Lemma 7.1 in \citet{feng2022test} with ``$|D_N|$", ``$N^{\varsigma}$" and ``$t$", respectively, we have that $|G_t|\leq N^{t+\varsigma-1}$ for each $2\leq t\leq N$. Therefore $\binom{|D_N|}{t}\geq |F_t|\geq \binom{|D_N|}{t}-N^{t+\varsigma-1}$. Since $D_N/N\to 1$  and $\varsigma=\varsigma_N\to 0$ as $N\to\infty$, we know
\bea\lbl{hong_huo}
\lim_{N\to\infty}\frac{|F_t|}{N^t}=\frac{1}{t!}.
\eea
Here
\beaa
\beta_t&=&\sum_{(i_1, \cdots, i_t)\in F_t} P(W_{i_1}>z, \cdots, W_{i_t}>z)+\\
 &&  \sum_{(i_1, \cdots, i_t)\in G_t} P(W_{i_1}>z, \cdots, W_{i_t}>z).
\eeaa
From Lemma \ref{winter_fire} and \eqref{hong_huo} we have
\beaa
\sum_{(i_1, \cdots, i_t)\in F_t} P(W_{i_1}>z, \cdots, W_{i_t}>z)\to \frac{1}{t!} \Big(\frac{ \kappa}{\Gamma\left(\mu_{1} / 2\right)} e^{-\frac{y}{2}}\Big)^t
\eeaa
as $N\to\infty.$ As a consequence, it remains to show
\bea\lbl{wkcrby}
\sum_{(i_1, \cdots, i_t)\in G_t} P(W_{i_1}>z, \cdots, W_{i_t}>z) \to 0
\eea
as $N\to \infty$ for each $t\geq 2.$

Next, we will prove \eqref{wkcrby}. If $t=2$, the sum of probabilities in \eqref{wkcrby} is bounded by
$|G_2|\cdot \max_{1\leq i<j\leq N}P(W_i>z, W_j>z )$. By Lemma 7.1 in \citet{feng2022test}, $|G_2|\leq N^{\varsigma+1}$. Since $|\sigma_{ij}|\leq \varrho$, by Lemma \ref{winter_cool}, \bea\lbl{country_love}
P(W_i>z, W_j>z ) \leq  \frac{(\log N)^C}{N^{(5-\varrho)/4}}
\eea
uniformly for all $1\leq i < j \leq N$ as $N$ is sufficiently large, where $C>0$ is a constant not depending on $N.$ We then know \eqref{wkcrby} holds. So the remaining job is to show  \eqref{wkcrby} for $t\geq 3.$

Let $N \geq 2$ and $\left(\sigma_{i j}\right)_{N \times N}$ be a non-negative definite matrix. For $\delta_N>0$ and a set $A \subset\{1,2, \cdots, m\}$ with $2 \leq m \leq N$, define
$$
\wp(A)=\max \left\{|S| ; S \subset A \text { and } \max _{i \in S, j \in S, i \neq j}\left|\sigma_{i j}\right| \leq \delta_N^{2+2c}\right\}
$$
Easily, $\wp(A)$ takes possible values $0,2 \cdots,|A|$, where we regard $|\varnothing|=0 .$ If $\wp(A)=0$, then $\left|\sigma_{i j}\right|>\delta_N^{2+2c}$ for all $i \in A$ and $j \in A$.

Now we will look at $G_t$ closely. To do so, we classify $G_t$ into the following subsets
\beaa
G_{t,j}=\big\{(i_1, \cdots, i_t)\in G_t;\, \wp(\{i_1, \cdots, i_t\})=j\big\}
\eeaa
for $j=0, 2, \cdots, t-1.$ By the definition of $G_t$, we see $G_t=\cup G_{t,j}$ for $j=0, 2, \cdots, t-1.$  Since $t\geq 3$ is fixed, to show \eqref{wkcrby}, it suffices to prove
\bea\lbl{what_tie}
\sum_{(i_1, \cdots, i_t)\in G_{t,j}} P(W_{i_1}>z, \cdots, W_{i_t}>z)  \to 0
\eea
for any $j\in\{0, 2, \cdots, t-1\}$.

Assume $(i_1, \cdots, i_t)\in G_{t,0}$. This implies that $|\sigma_{i_r i_s}|> \delta_N^{2+2c}$ for all  $1\leq r<s\leq t$. Therefore, the subgraph $\{i_1, \cdots, i_t\}\in  G_t$ is a clique. Taking $n=|D_N|\leq N$, $t=t$ and $q=N^{\varsigma}$. Then by Lemma 7.1 in \citet{feng2022test}, $|G_{t,0}| \leq N^{1+\varsigma(t-1)}\leq N^{1+t\varsigma}$.  Thus,  the sum from \eqref{what_tie} is bounded by
\bea\lbl{Huang}
N^{1+t\varsigma}\cdot \max_{1\leq i< j\leq N}P(W_i>z, W_j>z)\leq N^{1+t\varsigma}\cdot\frac{(\log N)^C}{N^{(5-\varrho)/4}}\to 0
\eea
as $N\to\infty$ by using \eqref{country_love}. So \eqref{what_tie} holds with $j=0.$

Now we assume $(i_1, \cdots, i_t)\in G_{t,j}$ with $j\in \{2, \cdots, t-1\}$. By definition, there exits $S\subset \{i_1, \cdots, i_t\}$ such that $\max_{ i\in S, j\in S, i\ne j}|\sigma_{ij}|\leq \delta_N^{2+2c}$ and for each $k\in \{i_1, \cdots, i_t\}\backslash S$, there exists $i\in S$ satisfying
$|\sigma_{ik}| > \delta_N^{2+2c}$. Looking at the last statement we see two possibilities: (i) for each $k\in \{i_1, \cdots, i_t\}\backslash S$, there exist at least two indices, say, $i\in S$, $j\in S$ with $i \ne j$ satisfying
$|\sigma_{ik}| > \delta_N^{2+2c}$ and $|\sigma_{jk}| > \delta_N^{2+2c}$; (ii) there exists $k\in \{i_1, \cdots, i_t\}\backslash S$ such that $|\sigma_{ik}| > \delta_N^{2+2c}$ for an unique $i\in S$. However, for $(i_1, \cdots, i_t)\in G_{t,j}$, (i) and (ii) could happen at the same time for different $S$, say, (i) holds for $S_1$ and  (ii) holds for $S_2$ simultaneously. Thus, to differentiate the two cases, we introduce following two definitions. Set
\bea\lbl{comrade}
H_{t,j}&=&\big\{(i_1, \cdots, i_t)\in G_{t,j};\, \ \mbox{there exist}\ S \subset \{i_1, \cdots, i_t\}\ \mbox{with}\ |S|=j\ \mbox{and} \nonumber\\
&&\max_{ i\in S, j\in S, i\ne j}|\sigma_{ij}|\leq \delta_N^{2+2c}\ \mbox{such that}
\ \mbox{for any}\ k \in \{i_1, \cdots, i_t\}\backslash S\ \mbox{there exist}\ r\in S, s\in S,\nonumber\\
&& r \ne s\ \mbox{satisfying}\
  \min\{|\sigma_{kr}|,  |\sigma_{ks}|\}>\delta_N^{2+2c}\big\}. \lbl{sdrtn}
\eea
Replacing ``$n$", ``$q$" and ``$t$" in Lemma 7.1 in \citet{feng2022test} with ``$|D_N|$", ``$N^{\varsigma}$" and ``$t$", respectively, we have that $|H_{t,j}|\leq t^t\cdot N^{j-1+ (t-j+1)\varsigma}$ for each $t\geq 3$.
Again, set
\bea\lbl{Take_ball}
H_{t,j}'&=&\big\{(i_1, \cdots, i_t)\in G_{t,j};\, \ \mbox{for any}\ S \subset \{i_1, \cdots, i_t\}\ \mbox{with}\ |S|=j\ \mbox{and} \nonumber\\
&&\max_{i\in S, j\in S, i\ne j}|\sigma_{ij}|\leq \delta_N^{2+2c}\
\mbox{there exists}\ k \in \{i_1, \cdots, i_t\}\backslash S\ \mbox{such that}\ |\sigma_{kr}|>\delta_N^{2+2c}\nonumber\\
&& \mbox{for a unique}\ r \in S \}.
\eea
From Lemma 7.1 in \citet{feng2022test}, we see $|H_{t,j}'| \leq t^t\cdot N^{j+(t-j)\varsigma}.$ It is easy to see $G_{t,j}=H_{t,j}\cup H_{t,j}'$. Therefore, to show \eqref{what_tie}, we only need to prove
\bea\lbl{what_tie_1}
\sum_{(i_1, \cdots, i_t)\in H_{t,j}} P(W_{i_1}>z, \cdots, W_{i_t}>z) \to 0
\eea
 and
\bea\lbl{what_tie_2}
\sum_{(i_1, \cdots, i_t)\in H_{t,j}'} P(W_{i_1}>z, \cdots, W_{i_t}>z) \to 0
\eea
as $N\to\infty$ for $j=2, \cdots, t-1$. In fact, let $S$ be as in \eqref{sdrtn}, then by using Lemma \ref{winter_fire}, the probability in \eqref{what_tie_1} is bounded by $P(\cap_{l\in S}\{W_{l}>z\})\leq C\cdot N^{-j}$ uniformly for all $S$ as $N$ is sufficiently large, where $C$ is a constant not depending on $N$. Thus,
\beaa
\sum_{(i_1, \cdots, i_t)\in H_{t,j}} P(W_{i_1}>z, \cdots, W_{i_t}>z)
&\leq & t^t\cdot N^{j-1+ (t-j+1)\varsigma}\cdot \big(C\cdot N^{-j}\big)\\
& \leq & (Ct^t)\cdot N^{-1+t\varsigma}
\eeaa
as $N$ is sufficiently large. By assumption $\varsigma=\varsigma_N\to 0$, we then get \eqref{what_tie_1}.

Now we show \eqref{what_tie_2}. Recall the definition of $H_{t,j}'$. For  $(i_1, \cdots, i_t)\in H_{t,j}'$, pick $S \subset \{i_1, \cdots, i_t\}$ with $|S|=j$, $\max_{i\in S, j\in S, i\ne j}|\sigma_{ij}|\leq \delta_N^{2+2c}$ and $k \in \{i_1, \cdots, i_t\}\backslash S$ such that $\delta_N^{2+2c}<|\sigma_{kr}|\leq \varrho$ for a unique $r \in S.$  Then the probability from  \eqref{what_tie_2} is bounded by
\beaa
P\Big(W_{k}>z, \bigcap_{l\in S}\{W_{l}>z\}\Big)
\eeaa
for $2\leq j \leq t-1$.  Taking $m=j+1$ in Lemma \ref{winter_cool}, then  the probability above is dominated by
\beaa
\frac{2^{j+1}}{z}\cdot  \exp\Big\{-\frac{z^2}{2}\Big(j+\frac{1-\varrho}{4}\Big) \Big\}=O\Big(\frac{(\log N)^{c_{1}}}{N^{j+(1-\varrho)/4}}\Big)
\eeaa
for some constant $c_{1}$ not depending on $N$. As stated earlier, $|H_{t,j}'| \leq t^t\cdot N^{j+(t-j)\varsigma}.$   Multiplying the two quantities, since $\varsigma=\varsigma_N\to 0$, we see the sum from \eqref{what_tie_2} is of order $O(N^{-(1-\varrho)/8})$. Therefore  \eqref{what_tie_2} holds. We then have proved \eqref{what_tie} for any $j\in\{0, 2, \cdots, t-1\}$. The proof is completed. \hfill$\square$

}

\paragraph{Proof of Theorem \ref{th31}}
We proceed in two steps, proving first the case $m=2$ and then generalizing to $m \geq 2$. For notational convenience we introduce the constants $b_{1}:=\|h\|_{\infty}<\infty$ and $b_{2}:=$ $\sup _{v}\left\|\phi_{v}\right\|_{\infty}<\infty$.
Similar to the proof of Theorem \ref{th21}, we define $\{u_s\}_{s=1}^{N}=\{U_{ij}(k)\}_{1\le i,j\le p, 1\le k\le K}$ and $\{\X_{tijk}\}=\{(\varepsilon_{t,i},\varepsilon_{t+k,j})^\top\}_{1\le t\le n-k}$. So we rewrite $U_{ij}(k)$ in the following forms
\begin{align}
u_s=\frac{1}{C_{n-k_s}^m}\sum_{1\le t_1<t_2,\cdots,<t_m\le n-k_s}h(\X_{t_1,i_s j_s k_s},\cdots,\X_{t_m, i_s j_s k_s}).
\end{align}

\noindent{\bf Step I.} Suppose $m=2$. We start with the scenario that there are infinitely many nonzero eigenvalues. For a large enough integer $M$ to be specified later, we define the ``truncated" kernel of
$h_{2}\left(z_{1}, z_{2} ; \mathbb{P}_{Z}\right)$ as $h_{2, M}\left(z_{1}, z_{2} ; \mathbb{P}_{Z}\right)=\sum_{v=1}^{M} \lambda_{v} \phi_{v}\left(z_{1}\right) \phi_{v}\left(z_{2}\right),$ with corresponding U-statistic
$$
u_{M, s}:=\left(\begin{array}{c}
n - k_s \\
2
\end{array}\right)^{-1} \sum_{1 \leq i<j \leq n-k_s} h_{2, M}\left(Z_{i}, Z_{j} ; \mathbb{P}_{Z}\right)
$$
For simpler presentation, define $Y_{v, i}=\phi_{v}\left(Z_{i}\right)$ for all $v=1,2, \ldots$ and $i \in[n-k_s] .$ In view of the expansions of $h_{2, M}(\cdot)$ and $h_{2}(\cdot), u_{M, s}$ and $u_{s}$ can be written as
$$
\begin{aligned}
 u_{M, s} &=\frac{1}{{n-k_s}-1}\left\{\sum_{v=1}^{M} \lambda_{v}\left({(n-k_s)}^{-1 / 2} \sum_{i=1}^{{n-k_s}} Y_{v, i}\right)^{2}-\sum_{v=1}^{M} \lambda_{v}\left(\frac{\sum_{i=1}^{{n-k_s}} Y_{v, i}^{2}}{{n-k_s}}\right)\right\} \\
 u_{s} &=\frac{1}{{n-k_s}-1}\left\{\sum_{v=1}^{\infty} \lambda_{v}\left({(n-k_s)}^{-1 / 2} \sum_{i=1}^{{n-k_s}} Y_{v, i}\right)^{2}-\sum_{v=1}^{\infty} \lambda_{v}\left(\frac{\sum_{i=1}^{{n-k_s}} Y_{v, i}^{2}}{{n-k_s}}\right)\right\}
\end{aligned}
$$
Define $M=[n^{(1-3\theta)/5}]. $ By the definition of $\theta$, there exist a positive absolute constant $C_{\theta}$ such that $\sum_{v=M+1}^\infty \lambda_v \le C_{\theta}n^{-\theta}$ for all sufficiently large $n$. Thus, for any $\epsilon>0$,
we have
\begin{align*}
P\left(\max_{1\le s\le N}(n - k_s -1)|u_{M,s}-u_s|\ge\epsilon\right)\le& NP\left((n - k_s-1)|u_{M,s}-u_s|\ge\epsilon\right)\\
\le &2Ne^{1/12}\exp\left(-\frac{\epsilon}{12b_2^2\sum_{v=M+1}^\infty \lambda_v}\right)\\
\le &2Ne^{1/12}\exp\left(-\frac{\epsilon n^{\theta}}{12b_2^2C_{\theta}}\right)\to 0
\end{align*}
by $\log N=o(n^{\theta})$. Here the second inequality are followed by (A.9) in \citet{drton2020high}.
Thus, by
\begin{align*}\label{u1}
  &\left|\max_{1\le s\le N}(n -k_s-1)u_s-\max_{1\le s\le N}(n -k_s-1)u_{M,s}\right|\le \max_{1\le s\le N}(n - k_s-1)|u_s-u_{M,s}|\to 0
\end{align*}
we only need to show that
\begin{align}
& \mathbb{P}\left\{\max_{1\le s\le N}(n - k_s -1)u_{M,s}-2 \lambda_1 \log (N)- \lambda_1\left(\mu_{1}-2\right) \log \log (N)+\Lambda \leq \lambda_{1}y\right\}\nonumber\\
&\to \exp \left\{-\frac{  \kappa}{\Gamma\left(\mu_{1} / 2\right)} \exp \left(-\frac{y}{2}\right)\right\}
\end{align}
Define
\begin{align*}
 \tilde u_{M, s} &=\frac{1}{{n-k_s}-1}\left\{\sum_{v=1}^{M} \lambda_{v}\left({(n-k_s)}^{-1 / 2} \sum_{i=1}^{{n-k_s}} Y_{v, i}\right)^{2}-\sum_{v=1}^{M} \lambda_{v}\right\}.
\end{align*}
Thus, for any $\epsilon>0$,
we have
\begin{align*}
P\left(\max_{1\le s\le N}(n-k_s-1)|u_{M,s}-\tilde u_{M, s}|\ge\epsilon\right)\le& NP\left((n-k_s-1)|u_{M,s}-\tilde u_{M, s}|\ge\epsilon\right)\\
\le &NP\left(\left|\sum_{v=1}^M\lambda_v\frac{\sum_{i=1}^{n-k_s}(Y_{v,i}^2-1)}{n-k_s}\right|\ge \epsilon\right)\\
\le &2N\exp\left(-\frac{(n-k_s)\epsilon^2}{48\Lambda^2(b_2^2+1)^2}\right)\to 0
\end{align*}
by $\log N=o(n^{\theta})$. Here the second inequality are followed by (A.10) in \citet{drton2020high}.
Thus, by
\begin{align*}\label{u1}
  &\left|\max_{1\le s\le N}(n - k_s -1)\tilde u_{M, s}-\max_{1\le s\le N}(n-k_s-1)u_{M,s}\right|\le \max_{1\le s\le N}(n-k_s-1)|\tilde u_{M, s}-u_{M,s}|\to 0
\end{align*}
we only need to show that
\begin{align}
& \mathbb{P}\left\{\max_{1\le s\le N}(n-k_s-1)\tilde u_{M, s}-2 \lambda_1 \log (N)- \lambda_1\left(\mu_{1}-2\right) \log \log (N)+\Lambda \leq \lambda_{1} y\right\}\nonumber\\
&\to \exp \left\{-\frac{  \kappa}{\Gamma\left(\mu_{1} / 2\right)} \exp \left(-\frac{y}{2}\right)\right\}
\end{align}
Define $\upsilon_s=\sum_{v=1}^{M} \lambda_{v}\left({(n-k_s)}^{-1 / 2} \sum_{i=1}^{{n-k_s}} Y_{v, i}\right)^{2}=(n-k_s-1)\tilde u_{M, s}+\sum_{v=1}^M \lambda_v$. By the definition of $\Lambda$, we have $\Lambda-\sum_{v=1}^M \lambda_v=\sum_{v=M+1}^\infty \lambda_v=O(n^{-\theta})$, thus, we only need to show that
\begin{align}
& \mathbb{P}\left\{\max_{1\le t\le N}\upsilon_t-2 \lambda_1\log (N)- \lambda_1\left(\mu_{1}-2\right) \log \log (N) \leq \lambda_{1} y\right\}\nonumber\\
&\to \exp \left\{-\frac{  \kappa}{\Gamma\left(\mu_{1} / 2\right)} \exp \left(-\frac{y}{2}\right)\right\}
\end{align}
Define $z=2  \lambda_1\log (N)+ \lambda_1\left(\mu_{1}-2\right) \log \log (N)+ \lambda_{1} y$.
%

By Theorem 4.1 in \citet{drton2020high} and $\log N=o(n^{\theta})$, we have
\begin{align*}
P(\upsilon_i \geq z)=\frac{ \kappa}{\Gamma\left(\mu_{1} / 2\right)}\left(\frac{z}{2 \lambda_{1}}\right)^{\mu_{1} / 2-1} \exp \left(-\frac{z}{2 \lambda_{1}}\right)\{1+o(1)\}\sim \frac{ \kappa}{\Gamma\left(\mu_{1} / 2\right)}\frac{e^{-y/2}}{N}
\end{align*}
Thus,
$$
P\left(\max _{i \in C_{N}}\upsilon_{i}>z\right) \leq\left|C_{N}\right| \cdot P(\upsilon_i \geq z) \rightarrow 0
$$
$\operatorname{asp} \rightarrow \infty .$ Set $D_{N}:=\left\{1 \leq i \leq N ;\left|B_{N, i}\right|<N^{\varsigma}\right\} .$ By assumption, $\left|D_{N}\right| / N \rightarrow 1$ as $N \rightarrow \infty$
Easily,
$$
\begin{aligned}
P\left(\max _{i \in D_{N}}\upsilon_{i}>z\right) & \leq P\left(\max _{1 \leq i \leq N}\upsilon_{i}>z\right) \\
& \leq P\left(\max _{i \in D_{N}}\upsilon_{i}>z\right)+P\left(\max _{i \in C_{N}}\upsilon_{i}>z\right)
\end{aligned}
$$
Therefore, to prove Theorem \ref{th31}, it is enough to show
$$
\lim _{N \rightarrow \infty} P\left(\max _{i \in D_{N}}\upsilon_{i}>z\right)=1-\exp \left\{-\frac{ \kappa}{\Gamma\left(\mu_{1} / 2\right)} \exp \left(-\frac{y}{2}\right)\right\}
$$
as $N \rightarrow \infty$. Define
$$
\beta_{t}=\sum^* P\left(\upsilon_{i_{1}}>z, \cdots,\upsilon_{i_{t}}>z\right)
$$
for $1 \leq t \leq N,$ where the sum runs over all $i_{1}<\cdots<i_{t}$ and $i_{1} \in D_{N}, \cdots, i_{t} \in D_{N} .$
By Lemma \ref{le67},
$$
\lim _{N \rightarrow \infty} \beta_{t}=\frac{1}{t !}\left(\frac{  \kappa}{\Gamma\left(\mu_{1} / 2\right)} e^{-\frac{y}{2}}\right)^{-t}
$$
for each $t \geq 1 .$
Then, by Bonferroni inequality,
$$
\sum_{t=1}^{2 k}(-1)^{t-1} \beta_{t} \leq P\left(\max _{i \in D_{N}}\upsilon_{i}>z\right) \leq \sum_{t=1}^{2 k+1}(-1)^{t-1} \beta_{t}
$$
for any $k \geq 1 .$  let $N \rightarrow \infty$, we have
$$
\begin{aligned}
\sum_{t=1}^{2 k}(-1)^{t-1} \frac{1}{t !}\left(\frac{  \kappa}{\Gamma\left(\mu_{1} / 2\right)} e^{-\frac{y}{2}}\right)^{t} & \leq \liminf _{N \rightarrow \infty} P\left(\max _{i \in D_{N}}\upsilon_{i}>z\right) \\
& \leq \limsup _{N \rightarrow \infty} P\left(\max _{i \in D_{N}}\upsilon_{i}>z\right) \leq \sum_{t=1}^{2 k+1}(-1)^{t-1} \frac{1}{t !}\left(\frac{  \kappa}{\Gamma\left(\mu_{1} / 2\right)} e^{-\frac{y}{2}}\right)^{t}
\end{aligned}
$$
for each $k \geq 1 .$ By letting $k \rightarrow \infty$ and using the Taylor expansion of the function $1-e^{-x}$, so we obtain the result.

\noindent{\bf Step II.}
For $m\ge 2$, by the Hoeffding decomposition, we have
\begin{align*}
u_s=C_m^2 H_{n-k_s}^{(2)}(\cdot, \mathbb{P}_{Z_s})+\sum_{\ell=3}^m C_m^\ell H_{n-k_s}^{(\ell)}(\cdot, \mathbb{P}_{Z_s})
\end{align*}
where for any measure $\mathbb{P}_{Z_s}$ and kernel $h$, $H_{n-k_s}^{(\ell)}\left(\cdot ; \mathbb{P}_{Z_s}\right)$ is the U-statistic based on the completely degenerate kernel $h^{(\ell)}\left(\cdot ; \mathbb{P}_{Z_s}\right)$ from (\ref{hl}) :
$$
H_{n-k_s}^{(\ell)}\left(\cdot ; \mathbb{P}_{Z_s}\right):=\left(\begin{array}{c}
n-k_s \\
\ell
\end{array}\right)^{-1} \sum_{1 \leq i_{1}<i_{2}<\cdots<i_{\ell} \leq n-k_s} h^{(\ell)}\left(Z_{i_{1}}, \ldots, Z_{i_{\ell}} ; \mathbb{P}_{Z_s}\right).
$$
To prove the result, we only need to show that $ \max_{1\le s\le N}(n-k_s-1)H_{n-k_s}^{(\ell)}(\cdot, \mathbb{P}_{Z_s})=o(1)$ for $\ell\ge 3$.

By Proposition 2.3(c) in \citet{arcones1993limit}, there exist positive constant $C_1, C_2$ such that for all $\epsilon_n>0$,
\begin{equation}
P\left((n-k_s)^{\ell/2}| H_{n-k_s}^{(\ell)}(\cdot, \mathbb{P}_{Z_s})|\ge \epsilon_n\right) \le C_1\exp\left(-C_2\left(\frac{\epsilon_n}{2^\ell b_1}\right)^{2/\ell}\right)
\end{equation}
So, for any $\epsilon_1>0$,
\begin{align*}
P\left(\max_{1\le s\le N} (n - k_s -1)H_{n-k_s}^{(\ell)}(\cdot, \mathbb{P}_{Z_s})\ge \epsilon_1\right)\le &NP\left((n-k_s-1)H_{n-k_s}^{(\ell)}(\cdot, \mathbb{P}_{Z_s})\ge \epsilon_1\right)\\
\le &C_1N\exp\left(-C_2\left(\frac{(n-k_s)^{\ell/2-1}\epsilon_1}{2^\ell b_1}\right)^{2/\ell}\right)\to 0
\end{align*}
by the condition $\log N=o(n^{\theta})$. Here we complete the proof. \hfill$\square$

\subsubsection{Proof of Theorem \ref{th32}}
The proof is similar to the proof of Theorem 4.3 in \citet{drton2020high}. So we omit it here.

\subsubsection{Proof of Theorem \ref{th33}}
According to Theorem 3 in \citet{feng2022whitenoise} and Assumption (A4), we can easily obtain the result. \hfill$\square$

{\color{black}
\subsection{Proof of Theorem S.1}
First, we restate the following lemma in \citet{arr1989}.
\begin{lemma}\label{chenstein}
Let $I$ be an index set and $\left\{B_\alpha, \alpha \in I\right\}$ be a set of subsets of I; that is, $B_\alpha \subset I$ for each $\alpha \in I$. Let also $\left\{\eta_\alpha, \alpha \in I\right\}$ be random variables. For a given $t \in \mathcal{R}$, set $\lambda=\sum_{\alpha \in I} \operatorname{pr}\left(\eta_\alpha>t\right)$. Then
$$
\left|\operatorname{pr}\left(\max _{\alpha \in I} \eta_\alpha \leq t\right)-e^{-\lambda}\right| \leq \min \left(1, \lambda^{-1}\right)\left(b_1+b_2+b_3\right),
$$
\end{lemma}
where
$$
\begin{aligned}
b_1 & \equiv \sum_{\alpha \in I} \sum_{\beta \in B_\alpha} \operatorname{pr}\left(\eta_\alpha>t\right) \operatorname{pr}\left(\eta_\beta>t\right), b_2 \equiv \sum_{\alpha \in I} \sum_{\beta \neq \alpha, \beta \in B_\alpha} \operatorname{pr}\left(\eta_\alpha>t, \eta_\beta>t\right), \\
b_3 & \equiv \sum_{\alpha \in I} E\left|\operatorname{pr}\left\{\eta_\alpha>t \mid \sigma\left(\eta_\beta, \beta \notin B_\alpha\right)\right\}-\operatorname{pr}\left(\eta_\alpha>t\right)\right|
\end{aligned}
$$
where $\sigma\left(\eta_\beta, \beta \notin B_\alpha\right)$ is the $\sigma$-algebra generated by $\left\{\eta_\beta, \beta \notin B_\alpha\right\}$. In particular, if $\eta_\alpha$ is independent of $\left\{\eta_\beta, \beta \notin B_\alpha\right\}$ for each $\alpha,$ then $b_3=0$.

Next, we adopt Lemma \ref{chenstein} to prove Theorem S1.

\proof In Lemma \ref{chenstein}, let $I=\{(i,j,k): 1\le i,j\le p, 1\le k\le K\}$. For $u=\{(i,j,k)\in I\}$, set $B_u=\{(l,m,q)\in I:\{i,j\}\cap\{l,m\}\not=\emptyset\}$, $\eta_u=\{|\psi_{ij}(k)|\}$, $\psi_{ij}(k)=\sqrt{\frac{5(n-k+1)}{2}}|\Xi_{ij}(k)|$ and $A_u=\{\psi_{ij}(k)>t_y\}$, $t_y=2\log N-\log\log N+y$. By Lemma \ref{chenstein}, we have $b_3=0$ by the independence assumption. Thus,
\begin{align*}
|P\left(\Xi_n^2\le t_y\right)-e^{-\lambda_n}|\le b_{1,n}+b_{2,n}
\end{align*}
where $\lambda_n=NP(\psi_{ij}(k)>t_y)$ and
\begin{align*}
b_{1,n}\le 2K^2p^3P^2(\psi_{ij}(k)>t_y)=O(p^{-1})\to 0
\end{align*}
since
\begin{align*}
P(\psi_{ij}(k)>t_y)&=P(|N(0,1)|>t_y)(1+o(1))\\
&=2(1-\Phi(t))(1+o(1)) =\sqrt{\frac{2}{\pi}}\frac{e^{-t^2/2}}{t}(1+o(1))
\end{align*}
by Theorem 2.1 in Chatterjee (2021). Additionally, by Lemma C4 in Han et al. (2017), we have $\Xi_{ij}(k)$ is independent of $\Xi_{is}(l)$ if $i\not=j\not=s$ or $k\not=l$. So
\begin{align*}
b_{2,n}&\le 2K^2p^3P(\psi_{ij}(k)>t_y,\psi_{is}(l)>t_y)+K^2pP(\psi_{ii}(k)>t_y,\psi_{ii}(l)>t_y)\\
&\le 2K^2p^3P^2(\psi_{ij}(k)>t_y)+K^2pP(\psi_{ii}(k)>t_y)=O(p^{-1})\to 0.
\end{align*}
Obviously, we have $\lambda_n\to \pi^{-1 / 2} \exp (-y / 2)$ as $p\to \infty$. So we obtain the result. \hfill$\Box$

}



\end{document}